\documentclass{article}
\usepackage{amsmath,amssymb,bm,yhmath}
\usepackage[all]{xy}
\usepackage{caption,graphicx}
\font\ibf=cmbxti10

\usepackage{multicol}

\title{All finitely presented groups are QSF}
\author{Valentin {\sc Po\'enaru}\footnote{Professor Emeritus, Universit\'e Paris-Sud, UMR 8628 du CNRS, Math\'ematiques, B\^atiment 425, 91405 Orsay Cedex, France. e-mail: valpoe@hotmail.com}}
\date{\sl (April 2014)}

\begin{document}

\maketitle

\vglue 1cm

\setcounter{section}{-1}
\section{Introduction}\label{sec0}
\setcounter{equation}{0}

This is the third and last of the trilogy of papers of which the first two are \cite{29} and \cite{39}, leading to the following final result.

\bigskip

\noindent {\bf Theorem A.} {\it All the finitely presented groups $\Gamma$ have the QSF property.}

\bigskip

Remember that the property in question has been introduced by S. Brick, M. Mihalik and J. Stallings (see \cite{3}, \cite{35}) and some general comments concerning it may be found in the introduction to \cite{39}. If we invoke some results of L. Funar and D. Otera \cite{7}, \cite{14}, \cite{41}, then there is also another way for stating Theorem A which some readers may find more congenial, namely

\bigskip

\noindent {\bf Theorem A$\bm '$.} (Alternative form of Theorem A.) {\it For any finitely presented group $\Gamma$ we may find a smooth closed manifold $M$ (of some high dimension), such the $\pi_1 M = \Gamma$ and that the universal covering space $\widetilde M$ is {\ibf geometrically simply connected} (GSC).}

\bigskip

Remember that GSC means that there is a handlebody decomposition s.t. the $1$-handles are in cancelling position with the $2$-handles, see here also \cite{30} and \cite {40}.

\smallskip

The rest of the present introduction is a brief survey of the proof of Theorem A, modulo the papers \cite{29}, \cite{39}, of which some tidbits will be reminded too.

\smallskip

In \cite{29}, for each $\Gamma$ we have constructed a presentation $\Gamma = \pi_1 M(\Gamma)$ along the following lines; $M(\Gamma)$ is a compact $3$-manifold, with {\ibf singularities}, and in \cite{39} we have introduced a certain $(N+4)$-dimensional cell-complex, with large $N$, called $S_u \, \widetilde M (\Gamma)$. Very roughly speaking, $S_u \, \widetilde M(\Gamma)$ is an infinitely foamy, high dimensional thickening of the universal covering space $\widetilde M(\Gamma)$. Actually, as explained in \cite{39}, the ``$S_u$'' is a functor. Here comes now our

\bigskip

\noindent {\bf Theorem B.} (The main result of \cite{39}, recalled here.) {\it The $S_u \, \widetilde M(\Gamma)$ is geometrically simply connected (GSC).}

\bigskip

The paper \cite{39} gives the full proof of Theorem B, relying strongly on \cite{29}. The paper \cite{42} is a coda to our trilogy.

\smallskip

There are very good reasons to work, not with the usual $2^{\rm d}$ presentations for $\Gamma$, but with $3^{\rm d}$ presentations. Let us say we have a presentation $M(\Gamma)$ which is a singular handlebody of some not yet determined dimension. Now, in order to get the local finiteness in \cite{29}, the first paper of the present trilogy, it was necessary that the handles of index one and two, and these are the ones which are really relevant in geometric group theory, should have co-cores of positive dimensions, allowing us to corral at infinity various unwanted infinite accumulations. This excludes the mundane dimension two for our presentations of $\Gamma$.

\smallskip

Next, the technology of \cite{39}, the second paper in the trilogy requires exploring every nook and hook of $M(\Gamma)$ with some dense subcomplexes and their zipping. That technology cannot work nicely if $\dim M(\Gamma) \geq 4$. So, eventually, the $\dim M(\Gamma) = 3$ get-forced on us.

\smallskip

The exact geometry of $S_u \, \widetilde M(\Gamma)$, explicitly explained in \cite{39}, will be very important for us in this paper. There is to begin with at $3^{\rm d}$ level, a first cell-complex $\Theta^3 (fX^2)$, then a $4^{\rm d}$ thickening of it $\Theta^4 (\Theta^3 (fX^2) , {\mathcal R})$, the notation $\Theta^4 (\ldots , {\mathcal R})$ being here like in \cite{8},  \cite{19},  \cite{36} and finally our $S_u \, \widetilde M(\Gamma)$ is, in a first approximation, but only in a first approximation,
\begin{equation}
\label{eq0.1}
S_u \, \widetilde M(\Gamma) = \Theta^4 (\Theta^3 (fX^2) , {\mathcal R}) \times B^N \, .
\end{equation}
All the three objects above, $\Theta^3$, $\Theta^4$, $S_u \, \widetilde M(\Gamma)$ are $\Gamma$-dependent. At a first, simple-minded level, all these three objects would be non locally finite, but this is certainly not something which we could live with. So, this lack of local finiteness is something which will have to be taken care of, requiring a certain amount of technology.

\smallskip

Except in the very special case when $\Gamma = \pi_1 M^3$, where $M^3$ is a smooth closed $3$-manifold, the $\Theta^3$ is never smooth, but if it would not be for that looming non local finiteness, the $\Theta^4$ and $S_u$ would be.

\smallskip

Local finiteness is realized by surging out the locus of non-local-finiteness and then, making up for this deletion, by the addition of a system of compensating $2$-handles of appropriate dimension. This will create singularities, i.e. non-manifold points. The singular locus certainly contains the attaching zones of the compensating $2$-handles and more. In the case of $\Theta^3$ there are other singularities too, while for $\Theta^4$ and $S_u \, \widetilde M(\Gamma)$ there are no others. But since there are singularities, we only have cell-complexes. For $\Theta^4$ or $S_u \, \widetilde M(\Gamma)$, the {\ibf correct} definition takes the following form

\bigskip

\noindent (0.2) \quad $\{$a non-compact smooth part of dimension four, respectively $N+4\} + \{$infinitely many compensating $2$-handles, also of dimension four or $N+4\}$. 

\bigskip

With this, the correct definition of the $S_u \, \widetilde M(\Gamma)$ which occurs in Theorem B is not (\ref{eq0.1}), but the following

\bigskip

\noindent (0.3) \quad $S_u \, \widetilde M(\Gamma) \equiv \{$The smooth part of the cell-complex $\Theta^4 (\Theta^3 (fX^2) , {\mathcal R})$, which is a smooth non-compact $(N+4)$-manifold, with very large boundary$\} \times B^N + \sum \, \{$compensating $2$-handles of dimension $N+4\}$.

\bigskip

Of course, one may ask, why not thicken to even higher dimensions and instead of a cell-complex like in (0.2), get a smooth manifold. The answer is that, in order to get from $S_u \, \widetilde M(\Gamma) \in {\rm GSC}$ to $\Gamma \in {\rm QSF}$, we need our $\Theta^4$ and $\Theta^3$ above, which certainly are singular. And, because of this, we need a {\ibf singular} $S_u \, \widetilde M(\Gamma)$, defined like in (0.3).

\smallskip

The group $\Gamma$ acts freely on each of the three objects $\Theta^3$, $\Theta^4$ and $S_u \, \widetilde M(\Gamma)$, once they are correctly defined, in the style of (0.2) or something more complicated, not to be described here, for $\Theta^3$. Unfortunately, none of the three actions above is cocompact.

\smallskip

But then, it turns out that there is a $\Gamma$-invariant subcomplex
\setcounter{equation}{3}
\begin{equation}
\label{eq0.4}
\Theta^3 (\mbox{co-compact}) \subset \Theta^3 (fX^1)
\end{equation}
which is {\ibf co}-compact. It occurs at the end of the following $\Gamma$-equivariant process
\begin{equation}
\label{eq0.5}
\Theta^3 (fX^2) \underset{\rm THE \ MULTI\mbox{-}GAME}{=\!\!=\!\!=\!\!=\!\!=\!\!=\!\!=\!\!=\!\!=\!\!=\!\!=\!\!=\!\!=\!\!=\!\!=\!\!=\!\!\Longrightarrow} \Theta^3 ({\rm new}) \underset{\rm collapse}{-\!\!\!-\!\!\!-\!\!\!-\!\!\!-\!\!\!-\!\!\!\longrightarrow} \Theta^3(\mbox{co-compact}) \, ,
\end{equation}
where the double arrow consists, in succession, of a PROPER, infinite $3^{\rm d}$ Whitehead dilatation, a PROPER addition of an infinite system of $3$-handles, followed by the cancellation of these $3$-handles with a PROPER system of $2$-handles, pre-existing in $\Theta^3 (fX^2)$. These handles to be cancelled are completely disjoined from the compensating $2$-handles which make good for the surging out of the non-local finiteness locus.

\bigskip

\noindent {\bf Remark.} The attaching zone of the compensating $2$-handles are far from the place where the deleted locus was. When we define correctly the $\Theta^3 (fX^2)$, something which takes a form analogous with (0.2), but non singular

\bigskip

\noindent $\Theta^3 (fX^2) \ (\mbox{correctly defined}) = \{$a $3^{\rm d}$ cell-complex which is a non-compact {\ibf singular} manifold, with undrawable singularities of the type described in \cite{8}, \cite{19}, \cite{36}$\} + \sum \, \{$compensating $2$-handles of dimension three$\}$,

\bigskip

\noindent then the double arrow, which we call the MULTI-GAME, stays far from the compensating $2$-handles. This means that we also have now a $\Theta^4 ({\rm new})$ defined like $\Theta^4 (\Theta^3 (fX^2) , {\mathcal R})$ with $\Theta^3 (fX^2)$ replaced by $\Theta^3 ({\rm new})$, and an $S_u ({\rm new})$. One of the effects of the the multi-game under discussion now, is to change the infinitely generated $\pi_2 \, \Theta^3 (fX^2)$ into a finitely generated $\pi_2 \, \Theta^3 ({\rm new})$. There is also here the following little fact

\bigskip

\noindent {\bf Lemma C.} {\it Because $S_u \, \widetilde M(\Gamma)$ is GSC, the $(N+4)$-dimensional $S_u ({\rm new})$ is also GSC.}

\bigskip

I will explain now the notion of {\ibf Dehn exhaustibility}, in the framework of {\ibf pure} $p$-dimensional complexes, denoted by $M^p , K^p , \ldots$. By definition, a pure $p$-dimensional simplicial complex $M^p$ is such that the maximum possible dimension of any simplex is $p$ and any simplex $\sigma$ of dimension $q < p$ is face of a $p$-dimensional complex.

\bigskip

\noindent {\bf Definition (0.6).} A pure $p$-complex $M^p$ is Dehn-exhaustible iff for any compact $k \overset{i}{\subset} M^p$ there is a compact, simply-connected pure $K^p$ which is abstract (i.e. not necessarily a subcomplex of $M^p$) and which comes with a commutative diagram
\setcounter{equation}{6}
\begin{equation}
\label{eq0.7}
\xymatrix{
k \ar@{-}[rr]^{j} \ar[dr]_-i &&K^p \ar[dl]^-{g}  \\ 
&M^p
}
\end{equation}
where $j$ is an inclusion, $g$ a simplicial {\ibf immersion} and where the following Dehn-type condition is fulfilled, for the set of double points $M_2 (g) \subset K^p$
\begin{equation}
\label{eq0.8}
i(k) \cap M_2(g) = \emptyset \, .
\end{equation}

\bigskip

If, in this context, $M^p$ is a smooth $p$-manifold, we may as well require that $K^p$ be a smooth $p$-manifold too and $g$ a smooth immersion. It is in this smooth connection, for $p=3$, that this concept first occurred in my old papers \cite{23}, \cite{24}, \cite{25} and, independently, in the work of A. Casson \cite{9} too. It is those old papers which motivated S. Brick, M. Mihalik and J. Stallings to introduce the concept QSF. In \cite{3} one also finds the following

\bigskip

\noindent {\bf Variant of Dehn's Lemma.} {\it Let $W^3$ be a smooth open $3$-manifold which is Dehn-exhaustible (which certainly implies that $\pi_1 W^3 = 0$). Then $W^3$ admits an exhaustion by compact codimension zero simply-connected submanifolds. Hence, we also have that $\pi_1^{\infty} W^3 = 0$.}

\bigskip

The proof follows the same pattern as for the classical Dehn's lemma. Our proof of Theorem A never makes use of this variant of Dehn's lemma, which I only mentioned here as a historical illustration. Actually our D.E. (Dehn-exhaustibility) implies QSF but, when it comes to groups $\Gamma$, but while QSF is presentation independent, DE is not.

\smallskip 

With a little additional work, from \cite{23}, \cite{24}, \cite{26} one can extract a proof of the following fact, which should be kept in mind for what will follow afterwards.

\bigskip

\noindent {\bf Proposition D.} {\it Let $V^p$ be a smooth open $p$-manifold, such that there exists some $m \in Z_+$ with the property that $V^p \times B^m$ is GSC. Then $V^p$ is DE.}

\bigskip

Our next lemma is now

\bigskip

\noindent {\bf Lemma E.} {\it The fact that $S_u ({\rm new})$ is GSC implies that $\Theta^4 ({\rm new})$ is Dehn-exhaustible, in the context of pure $4$-complexes.}

\bigskip

The proof is a relatively easy modification of the proof of proposition {\rm D} which, as we have said can be done like in \cite{23}, \cite{24}, \cite{26}.

\smallskip

By more or less similar, but harder arguments, because the situation is now more singular, one can prove

\bigskip

\noindent {\bf Lemma F.} {\it The fact that $\Theta^4 ({\rm new})$ is DE implies that $\Theta^3 ({\rm new})$ is also Dehn-exhaustible.}

\bigskip

In the proof of Lemma F, the Dehn-exhaustibility of $\Theta^4 ({\rm new})$ replaces the GSC property of $V^p \times B^m$ from the context of Proposition D. Similarly, the canonical retraction
$$
\Theta^4 ({\rm new}) \equiv \Theta^4 (\Theta^3 ({\rm new}) , {\mathcal R}) \overset{r}{-\!\!\!-\!\!\!-\!\!\!\longrightarrow} \Theta^3 ({\rm new}) \, ,
$$
plays in the proof of our lemma E the same role as the projection $V^p \times B^m \overset{\pi}{-\!\!\!-\!\!\!-\!\!\!\longrightarrow} V^p$, in the context of Proposition D.

\smallskip

The final step in our proof of Theorem A is now the following.

\bigskip

\noindent {\bf Lemma G.} {\it Using the fact that $\Theta^3 ({\rm new})$ is DE and making also use of a complete knowledge of the structure of the collapse $\Theta^3 ({\rm new}) \longrightarrow \Theta^3 (\mbox{\rm co-compact})$ from {\rm (0.5)}, one can show that $\Theta^3 (\mbox{\rm co-compact})$ is QSF.}

\bigskip

Since there is a free co-compact action
$$
\Gamma \times \Theta^3 (\mbox{co-compact}) \longrightarrow \Theta^3 (\mbox{co-compact}) \, ,
$$
our Lemma G implies the Theorem A.

\smallskip

Of course, in the next pages, this fast overview of the proof of Theorem A will be developed with full details.

\smallskip

Finally, there is also a CODA to the trilogy, namely the paper \cite{42}, to be very soon available too.

\vglue 1cm

Thanks are due to David Gabai and Louis Funar for very helpful conversations. I also wish to thank, once more the IHES for its constant friendly help, and last but not least, my many thanks are due to C\'ecile Gourgues for the typing of this paper, and to Marie-Claude Vergne for the drawings.

\newpage

\section{The game}\label{sec1}
\setcounter{equation}{0}

We start by  reviewing the geometrical objects which the present paper will have to deal with. All these objects have been already introduced in \cite{39}, a paper of which the present one is a direct continuation. In terms of this \cite{39}, we will be here constantly in the context of the Variant II, and we will repeat right now the little exposition from the Complement (6.21.5) in \cite{39}. All the references to numbers between prentices, until further notice, will refer to \cite{39}.

\smallskip

One starts with the $\Theta^3 (fX^2)$ from (2.12). This object, as such, contains already all the fins $F_{\pm}$ (minus their rims, as it will turn out), has the $\partial \Sigma (\infty)^{\wedge}$ ($\supset$ rims of fins) deleted, AND IT FAILS to be locally finite at the $p_{\infty\infty} (S)$'s (see here (1.15.0)). Next, as part of the big passage from Variant I to Variant II, one adds to the $\Theta^3 (fX^2)$ the $\underset{R_0}{\sum} \, {\rm int} \, R_0 \times [0,\infty)$, far from the $p_{\infty\infty} (S)$'s. We will review now, completely, how in Section VI of \cite{37} one perform the change from the Variant I to Variant II, and see here also (6.21) in \cite{37}. We start, like in (6.18), with
\setcounter{equation}{-1}
\begin{equation}
\label{eq1.0}
\xymatrix{
\left(\underset{R_0}{\sum} \, R_0 , \underset{R_0}{\sum} \, \partial R_0 \right) \ar[rr]_{\varphi} &&\left(\Sigma (\infty)_*^{\wedge} , \partial \Sigma (\infty)_*^{\wedge} \right) \\ 
\underset{R_0}{\sum} \ {\rm int} \, R_0 \ar[u] \ar[rr] &&{\rm int} \left( \ring\Sigma (\infty)_* \cup {\rm fins} \right), \ar[u]
}
\end{equation}
where ${\rm int} \left( \ring\Sigma (\infty)_* \cup {\rm fins} \right)$ is defined like in (6.8.1) \cite{39}, with $p_{\infty\infty} ({\rm all}) \times [-\varepsilon, \varepsilon]$ deleted. In VI \cite{39} it was essential to work with $S'_u (M(\Gamma) - H)_{\rm II} = S'_u (\widetilde M (\Gamma) - H)_{\rm II} \diagup \Gamma$ and in order to define it, we {\ibf had} to start from
$$
\Theta^3 (fX^2 - H)'_{\rm II} =\{\mbox{the $\Theta^3 (fX^2-H)'$ from (4.13.1) in \cite{39}}\} \cup \sum_R {\rm int} \, R_0 \times [0,\infty) \, ,
$$
where the two pieces are glued along ${\rm int} \left( \ring\Sigma (\infty)_* \cup {\rm fins} \right)$.

\smallskip

As a preliminary for proving that $S_u \, \widetilde M (\Gamma)_{\rm II} \in {\rm GSC}$, it was shown in Section VI of \cite{39} that $S'_u \, \widetilde M (\Gamma)_{\rm II} \in {\rm GSC}$. The context $S'_u$ was essential there, for proving the compactness lemma. In the present paper we start directly from the fact that $S_u \, \widetilde M (\Gamma)_{\rm II} \in {\rm GSC}$, and the context $S'_u$ is, by now, a mere intermediary tool which we will forget about. So, without loosing the all-important GSC feature, we can proceed now slightly differently than above. Like in (6.21.5) in \cite{39} which supersedes the (6.18), we will start by extending the range of $\underset{R_0}{\sum} \ {\rm int} \, R_0$ in (\ref{eq1.0}), from ${\rm int} \left( \ring\Sigma (\infty)_* \cup {\rm fins} \right)$ to
$$
\ring\Sigma (\infty)^{\wedge} \equiv \Biggl\{{\rm int} \left( \ring\Sigma (\infty)_* \cup {\rm fins} \right), \ \mbox{with all the contribution of $p_{\infty\infty} ({\rm proper})$ restored back$\Biggl\}$} \eqno (1.0.1)
$$
$$
\supsetneqq {\rm int} \left( \ring\Sigma (\infty)_* \cup {\rm fins} \right).
$$

With this, we define now the presently useful
\newpage
\begin{equation}
\label{eq1.1}
\Theta^3 (fX^2)_{\rm II} \equiv \Biggl[\{\Theta^3 (fX^2) \ \mbox{(from (2.12) \cite{39}) with the contribution of $p_{\infty\infty}(S)$ deleted}\} \ \underset{\overbrace{\mbox{\footnotesize$\ring\Sigma(\infty)^{\wedge}$}}}{\cup} 
\end{equation}
$$
\underset{R_0}{\sum} \ {\rm int} \, R_0 \times [0,\infty)\Biggl] + \, \Bigl\{\mbox{the compensating $2$-handles} \ \underset{P_{\infty\infty} (S)}{\sum} \, D^2 (p_{\infty\infty} (S)) \times \left[ -\frac\varepsilon4 , \frac\varepsilon4 \right]\Bigl\} \, .
$$

The piece $[\ldots]$ in (1.1) will be denoted by $[\Theta^3]_{\rm II}$. With ${\rm int} \, \Sigma (\infty)$ defined like in (2.13.1) \cite{39}, i.e. with the contribution of $p_{\infty\infty} (S)$ deleted, we have now
$$
\bigcup_{\overbrace{\mbox{\footnotesize$\ring\Sigma(\infty)^{\wedge}$}}} \ \sum_{R_0} {\rm int} \, R_0 \times [0,\infty) = \bigcup_{\overbrace{\mbox{\footnotesize${\rm int} (\Sigma (\infty)$ (2.13.1))}}} \ \left( {\rm int} \, \Sigma (\infty) \right) \times [0,\infty) \, .
\eqno ({\rm 1.1.bis})
$$
Next, we go $4$-dimensional and introduce the cell-complex
\begin{equation}
\label{eq1.2}
\Theta^4 (\Theta^3 (fX^2), {\mathcal R})_{\rm II} \equiv \Theta^4 ([\Theta^3]_{\rm II} , {\mathcal R}) \ (\mbox{which is smooth}) \ + 
\end{equation}
$$
+ \sum_{p_{\infty\infty} (S)} D^2 (p_{\infty\infty} (S)) \times \left[ -\frac\varepsilon4 , \frac\varepsilon4 \right] \times I \underset{\pi_{4,3}}{-\!\!\!-\!\!\!-\!\!\!-\!\!\!-\!\!\!\longrightarrow} \Theta^3 (fX^2)_{\rm II} \, .
$$
Here $\pi_{4,3} \mid \Theta^4 ([\Theta^3]_{\rm II} , {\mathcal R}) = \{$the natural retraction on $[\Theta^3]_{\rm II}$, of which $\Theta^4 (\ldots)$ is a smooth regular neighbour\-hood$\}$, and here also $\pi_{4,3} \mid \{\mbox{$2$-handle}\}$ is the obvious projection
$$
D^2 \times \left[ -\frac\varepsilon4 , \frac\varepsilon4 \right] \times I \longrightarrow D^2 \times \left[ -\frac\varepsilon4 , \frac\varepsilon4 \right] \, .
$$

Finally, we go high-dimensional (i.e. $(N+4)$-dimensional, with $N$ high) and introduce there a cell-complex
\begin{equation}
\label{eq1.3}
S_u \, \widetilde M(\Gamma)_{\rm II} \equiv \Theta^4 ([\Theta^3]_{\rm II} , {\mathcal R}) \times B^N + \sum_{p_{\infty\infty} (S)} D^2 (p_{\infty\infty} (S)) \times \left[ -\frac\varepsilon4 , \frac\varepsilon4 \right] \times I \times \frac12 \, B^N
\end{equation}
$$
\underset{\pi_{N+4,4}}{-\!\!\!-\!\!\!-\!\!\!-\!\!\!-\!\!\!-\!\!\!\longrightarrow} \ \Theta^4 (\Theta^3 (fX^2),{\mathcal R})_{\rm II} \, .
$$

Very importantly, while $\Theta^4 ([\Theta^3]_{\rm II} , {\mathcal R})$, where ${\mathcal R}$ is a desingularization, like in \cite{8}, \cite{21}, is ${\mathcal R}$-dependent, this dependence gets washed away when one goes from $\Theta^4 ([\Theta^3]_{\rm II} , {\mathcal R})$ to $S_u \, \widetilde M(\Gamma)_{\rm II}$. So, just like it was the case for $\Theta^3 (fX^2)_{\rm II}$, the $S_u \, \widetilde M(\Gamma)_{\rm II}$ admits now a free $\Gamma$-action which is co-compact ($=$ with compact fundamental domain).

\smallskip

I will restate now the main result of \cite{39}, namely

\bigskip

\noindent {\bf The statement (1.3.1).} {\it The $(N+4)$-dimensional cell-complex $S_u \, \widetilde M(\Gamma)_{\rm II}$, which fails to be smooth exactly along the $\underset{p_{\infty\infty} (S)}{\sum} C (p_{\infty\infty} (S)) \times \left[ -\frac\varepsilon4 , \frac\varepsilon4 \right] \times I \times \frac12 \, B^N$, is GSC.}

\bigskip

In \cite{39}, the present statement (1.3.1) had appeared as point 2) in the GSC Theorem 2.3.

\smallskip

From now on, the numbers of our formulae will no longer refer to \cite{39}, unless explicitly said so. This was already the case with (\ref{eq1.1}) to (\ref{eq1.3}).

\bigskip

\noindent {\bf A Remark.} Notice the sequence of increases and dimensions, throughout this series of papers:
$$
\left(X^2 \overset{\rm zipping}{-\!\!\!-\!\!\!-\!\!\!-\!\!\!-\!\!\!-\!\!\!\longrightarrow} \, fX^2 \right) \Longrightarrow \Theta^3 (fX^2) \Longrightarrow \Theta^4 (\Theta^3 , {\mathcal R}) \Longrightarrow S_u (\dim = N+4) \, .
$$

The zipping is best dealt with in $2^{\rm d}$, but in order to get to the all-important GSC feature, we need to go high-dimensional. The {\ibf geometric realization} of the zipping, our key to GSC, takes place essentially in the {\ibf supplementary dimensions} (those which are in addition to four).

\smallskip

Since the Variant I from \cite{39} will never any longer occur in this present paper, the subscript ``II'' for the objects defines in (\ref{eq1.1}) to (\ref{eq1.3}) above, may often be dropped.

\smallskip

We present now the {\ibf elementary game}, a transformation conceived  a priori at the level of (\ref{eq1.1}) (then at the other two levels above too). This is a semi-local process, generically labelled by an $\{$ideal Hole$\} \subset \bigcup \mbox{limit walls} \equiv \Sigma_1 (\infty)$ (see (1.14) in \cite{39}) $=$
$$
= \sum S_{\infty}^2 ({\rm BLUE}) \cup \sum (S^1 \times I)_{\infty} ({\rm RED}) \cup \sum {\rm Hex}_{\infty} ({\rm BLACK}) \subset \widetilde M (\Gamma) \, ,
$$
OR in the degenerate cases by an arc (which could possibly be reduced to a single point contained in the intersection of two limit walls (of different colours)). Contrary to the ideal Holes which correspond to exactly one GAME, these arcs can correspond to several such, possibly infinitely many. The BLUE, RED, BLACK elementary games will always be localized inside the part of $fX^2$ restricted to some handle of $\widetilde M (\Gamma)$, explicitly: a $h^0 ({\rm BLUE})$, a $h^1 ({\rm RED})$ plus the adjacent $h^0$'s (now RED/ BLUE), or finally $h^2 ({\rm BLACK})$ and the adjacent $h^0 , h^1$'s.

\smallskip

At the bottom of the geometric structure coming with an elementary game, we always find a $2$-cell called Sq like ``Square'', see here the formulae (\ref{eq1.4}), (\ref{eq1.21}), (\ref{eq1.24}), and also the figures 1.2, 1.6. The Sq is, according to the case, a piece of some compact wall $W({\rm BLUE}) , W({\rm RED}) , W({\rm BLACK})$. So much for the COLOURS attached to the elementary games.

\smallskip

We will start with the easiest, paradigmatical BLUE case and, in the simplest of the BLUE variants one considers first the $U^2(B)$ from formula (\ref{eq1.4}) below. Eventually this should be part of $fX^2$, with $(X^2,f)$ like in (1.1) from \cite{39} but, for simplicity's sake think of it now as living in $R^3$. Here it is
\begin{equation}
\label{eq1.4}
U^2 (B) = \left\{ \underbrace{[-1 \leq x \leq 1 , -1 \leq y \leq 1 , z=0]}_{\mbox{\footnotesize call this Sq, like ``square''}} \, \cup \, \partial \, {\rm Sq} \times [0 \leq z \leq N + \varepsilon_1] \right\} \cup
\end{equation}
$\cup \ \{$infinitely many $2$-handles (i.e. here $2$-cells), parallel to Sq and being BLUE, like it, namely the ${\rm Sq} \times \{ z_1 \} , {\rm Sq} \times \{z_2 \} , \ldots$, where $0< z_1 < z_2 < z_3 \ldots < N$ and $\underset{n = \infty}{\lim} z_n = N \}$.

\smallskip

\noindent In this simplest of the BLUE variants, the $[-1 \leq x \leq 1] \times [0 \leq z \leq N+\varepsilon_1] \times \{ y = \pm 1 \}$, respectively $[-1 \leq y \leq 1] \times [0 \leq z \leq N+\varepsilon_1] \times \{ x = \pm 1 \}$ are (pieces of) BLACK, respectively RED walls. In the non-generic, more complicated variants, the BLACK (or RED) walls might be replaced by an infinite BLACK/BLUE (or RED/BLUE) staircase, stretching through $[0 \leq z < N)$. At $z = N$ we have, generically, an ideal Hole $\subset \, S_{\infty}^2$. Figure 1.1 suggests, schematically, what we are talking about here.

$$
\includegraphics[width=12cm]{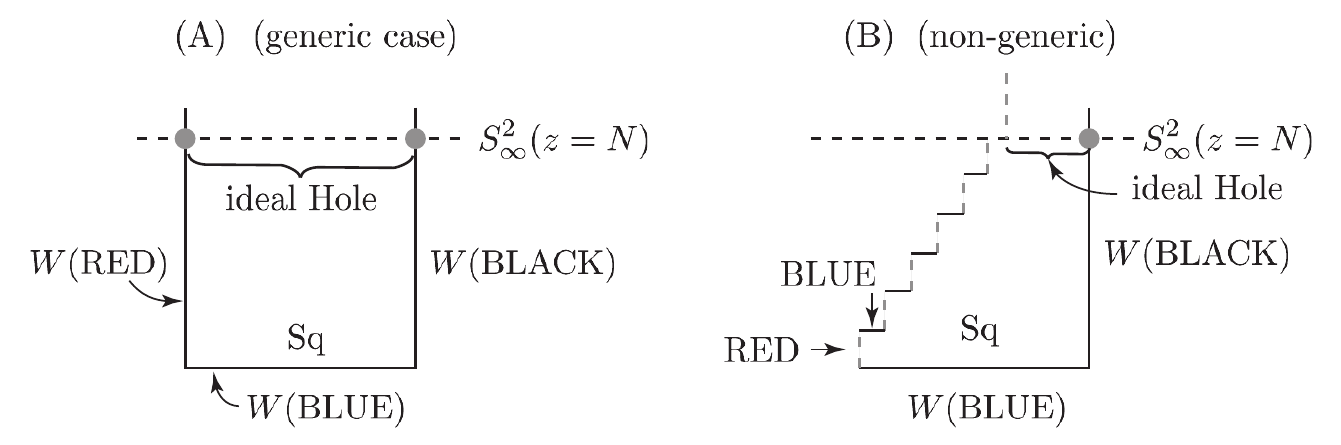}
$$
\label{fig1.1}

\centerline {\bf Figure 1.1.} 

\smallskip

\begin{quote}
Schematical representations of the generic $U^2({\rm BLUE})$ and of one of its variants. There is here, for instance, an additional variant where the straight $W({\rm BLACK})$ in (B) is replaced by another BLUE/RED infinite staircase and where the ideal hole which is squashed is reduced to an ideal arc contained in $S_{\infty}^1 = S_{\infty}^2 \cap (S^1 \times I)_{\infty}$.
\end{quote}

\bigskip

\noindent {\bf Remark.} In (\ref{eq1.4}) and also in the other similar formulae, the Sq instead of being a square, could be a polygon with more than four sides. \hfill $\Box$

\bigskip

What follows next, is a complement to the formula (\ref{eq1.4}), and it concerns the {\ibf ``lateral walls''} piece of (\ref{eq1.4}), by  which we mean the piece $\partial \, {\rm Sq} \times [0 \leq z < \ldots]$.

\smallskip

Specifically for the BLUE case, the following will happen

\bigskip

\noindent (1.4.1) \quad When the $U^2 (B)$ corresponds to an ideal (BLUE) Hole, and not to some arc in $S_{\infty}^1$, then the lateral part of (\ref{eq1.4}), even when it is an infinite staircase, exists already at the level $X^2$ without us having to go to $fX^2$. 

\bigskip

This certainly concerns the two drawings in Figure~1.1. \hfill $\Box$

\bigskip

When we move from $fX^2$ to the $\Theta^3 (fX^2)$ (\ref{eq1.1}), then $U^2(B)$ is to be replaced by the $U^3(B)$ below, essentially its regular neighbourhood, and here $0 < \varepsilon \ll \varepsilon_1$:
\begin{equation}
\label{eq1.5}
U^3(B) = U^2(B) \times [-\varepsilon,\varepsilon] - \{\partial \, {\rm Sq} \times [(z=N) \times \varepsilon]\} \, , \ \mbox{occurring as $S_{\infty}^1$ in Figure 1.2.}
\end{equation}
Here the factor $[-\varepsilon,\varepsilon]$ is supposed to be such that the $+\varepsilon$ is pointing towards the interior of the Sq. The deleted part of the formula is in $\partial \Sigma (\infty)$, with a $\partial \Sigma (\infty)$ like in (2.13.1) from \cite{39} and, very importantly, the spots via which our $U^3(B)$ communicates with the outside world are exactly the following ones:
\begin{equation}
\label{eq1.6}
\{\mbox{the outer $\varepsilon$ side of} \ U^2 \times [-\varepsilon,\varepsilon]\} \cup \{{\rm the} \ z > N \} \, ,
\end{equation}
to which we have to add the following item too

\bigskip

\noindent (1.6.1) \quad We are now in the context (\ref{eq1.1}) with $\underset{R_0}{\sum} \, {\rm int} \, R_0 \times [0,\infty)$ resting, among other things, on
$$
\Sigma_1 (\infty) \cap U^3({\rm BLUE}) \, .
$$

We will not add this kind of contribution to our $U^3 ({\rm COLOUR})$, it will never touch the bowls ${\mathcal B}$, and it will not interfere with the various constructions in the present section, which will have as their climax the MAIN MULTIGAME LEMMA 1.5.

\smallskip

Notice that, in (\ref{eq1.5}), the ${\rm Sq} \times (z=N)$, resting on the $S_{\infty}^1 \equiv \partial \, {\rm Sq} \times [(z=N) \times (-\varepsilon)] \subset S_{\infty}^2$ is an ideal Hole of BLUE colour. We will embellish $U^3(B)$ with a PROPER hypersurface
\begin{equation}
\label{eq1.7}
{\mathcal B} (\mbox{like ``BOWL''}) = \{\mbox{a copy of $R^2$ PROPERLY embedded inside ${\rm int} \, U^3(B) \subset U^3(B)$,}
\end{equation}
$$
\mbox{resting, at infinity, on $S_{\infty}^1$}\} \, .
$$

The ${\mathcal B}$ will be ``sent to infinity'' by adding to $U^3(B)$ a copy of ${\mathcal B} \times [0,\infty)$ along ${\mathcal B} = {\mathcal B} \times \{0\}$. Figure 1.2 suggests the embellished $U^3(B)$. The position of the $\overset{\infty}{\underset{n=1}{\sum}} \partial H_n^3$, attaching zones of the $3$-handles $H_1^3 , H_2^3, H_3^3,\ldots$ which the Figure 1.2 suggests us to attached to $U^3({\rm BLUE})$, should be slightly changed, with respect to what we see in the drawings, by letting the $\partial H^3$'s climb at least partially on ${\mathcal B} \times [0,\infty)$ so that we should fulfill the following condition
\begin{equation}
\label{eq1.8}
\lim_{n=\infty} \partial H_n^3 \subset ({\mathcal B} \times \{\infty\}) \cup S_{\infty}^2
\end{equation}
making the embedding $\overset{\infty}{\underset{1}{\sum}} \, \partial H_n^3 \subset U^3(B) \cup {\mathcal B} \times [0,\infty)$ PROPER.

$$
\includegraphics[width=12cm]{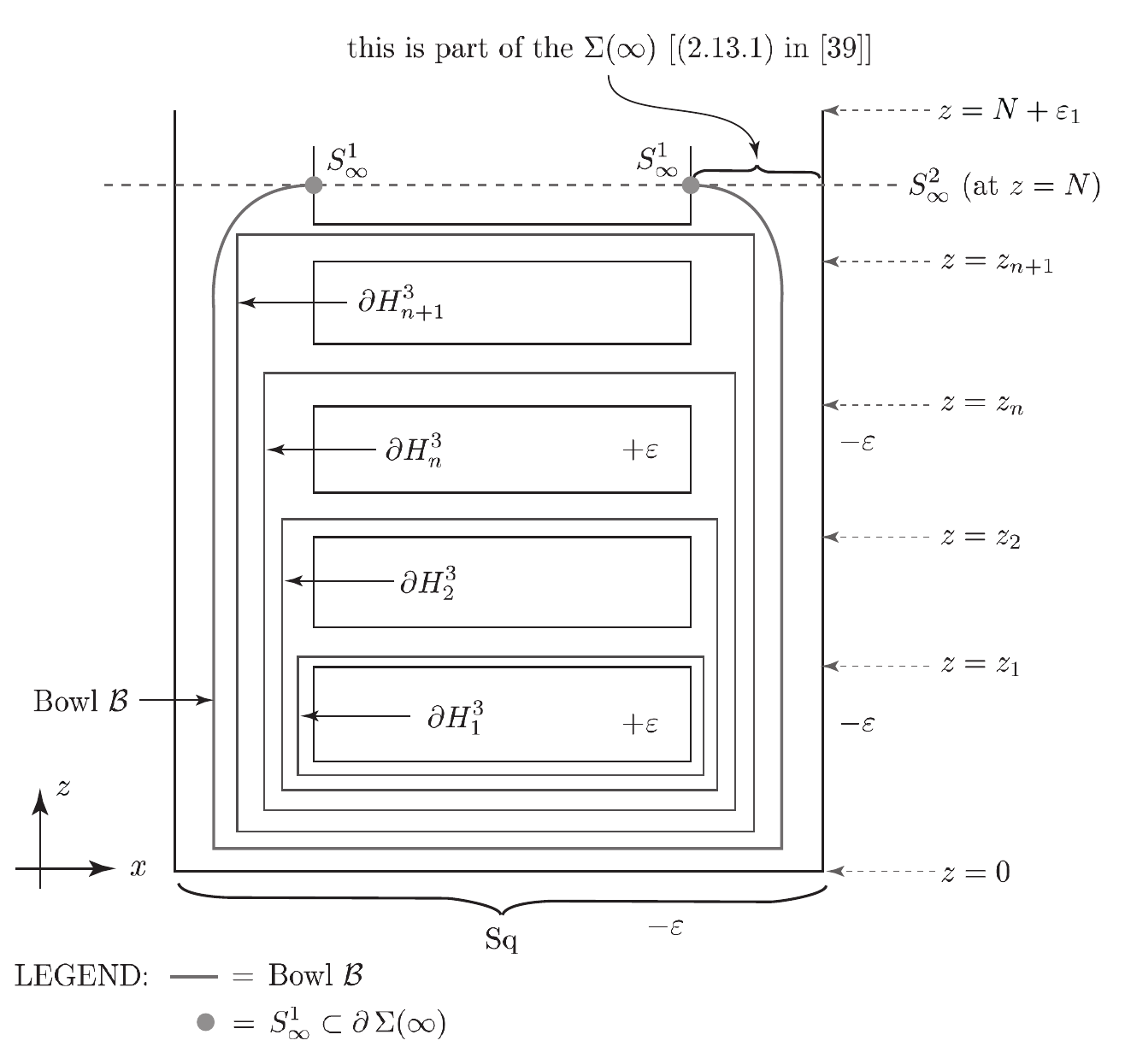}
$$
\label{fig1.2}
\centerline {\bf Figure 1.2.} 

\smallskip

\begin{quote}
The $U^3({\rm Blue})$ embellished with the BOWL ${\mathcal B}$ and with infinitely many $3$-handle attaching spheres $\partial H_1^3 , \partial H_2^3 , \ldots$ which accumulate on ${\mathcal B} \cup S_{\infty}^2$.
\end{quote}

\bigskip

With $U^3 = U^3(B)$ like above, we introduce now the BLUE transformation which {\ibf is} our elementary BLUE game
\begin{equation}
\label{eq1.9}
U^3 \overset{\rm BLUE}{\underset{\rm transformation}{=\!\!=\!\!=\!\!=\!\!=\!\!=\!\!=\!\!=\!\!=\!\!=\!\!=\!\!=\!\!=\!\!\Longrightarrow}} U^3({\rm new}) \equiv \Bigl\{ U^3(B) \cup {\mathcal B} \cup [0,\infty) \, \mbox{with all the 2-handles}
\end{equation}
$$
{\rm Sq} \times \{z_1\} , {\rm Sq} \times \{ z_2 \} , \ldots \ \mbox{\ibf deleted}\Bigl\} \, .
$$

This transformation does not touch to (\ref{eq1.6}) and so, when $U^3({\rm BLUE})$ is part of a larger (singular) $3^{\rm d}$ object, like the $\Theta^3 (fX^2)$ (\ref{eq1.1}) for instance, let us call this $X^3 \supset U^3({\rm BLUE})$, then one can go from the semilocal (\ref{eq1.9}) to a more global BLUE transformation
\begin{equation}
\label{eq1.10}
X^3 \equiv X^3 ({\rm old}) \overset{\rm BLUE}{\underset{\rm transformation}{=\!\!=\!\!=\!\!=\!\!=\!\!=\!\!=\!\!=\!\!=\!\!=\!\!=\!\!=\!\!=\!\!\Longrightarrow}} X^3({\rm new}) \, .
\end{equation}

\bigskip

\noindent {\bf BLUE Lemma 1.1.} {\it In the context of {\rm (1.10)}, assume that $\Theta^4 (X^3 ({\rm old}) , {\mathcal R}) \times B^N$ is GSC, then the $\Theta^4 (X^3 ({\rm new})$, ${\mathcal R}) \times B^N$ is also GSC.}

\bigskip

It should be understood here that the $\Theta^4 \times B^N$ in the statement above may be read like the $S_u$ in the formula (\ref{eq1.3}) and, anyway, in this context we will always have things like
\begin{equation}
\label{eq1.11}
U^3 ({\rm BLUE}) \cap \left\{ D^2 (p_{\infty\infty} (S)) \times \left[ - \frac\varepsilon4 , \frac\varepsilon4 \right] \ \mbox{in (\ref{eq1.1})} \right\} = \emptyset \, .
\end{equation}

\bigskip

\noindent {\bf Proof of the BLUE Lemma.} Because of (\ref{eq1.8}) we can use the $\underset{n=1}{\overset{\infty}{\sum}} \partial H_n^3$ as a recipee for attaching a PROPER infinite system of $(N+4)$-dimensional handles of index $\lambda = 3$, which we call $\underset{1}{\overset{\infty}{\sum}} \, H_n^3$, to $\Theta^4 (X^3 ({\rm old}) , {\mathcal R}) \times B^N$.  We get then
\begin{equation}
\label{eq1.12}
(\Theta^4 (X^3 ({\rm old}) , {\mathcal R}) \times B^N) + \sum_{n=1}^{\infty} H_n^3 \in {\rm GSC} \, .
\end{equation}
The $3$-handles above are in cancelling position with the $2$-handles of $\Theta^4 (U^3 ({\rm BLUE}) , {\mathcal R}) \times B^N \subset \Theta^4 (X^3 ({\rm old})$, ${\mathcal R}) \times B^N$. Actually, Figure 1.2 tells us that the geometric intersection  matrix is $\partial H_i^3$. ${\rm Sq} \times \{ z_j \} = \delta_{ij}$. It follows that we have a diffeomorphism
$$
\Theta^4 (X^3 ({\rm new}) , {\mathcal R}) \times B^N \underset{\rm DIFF}{=} \Theta^4 (X^3 ({\rm old}) , {\mathcal R}) \times B^N + \sum_1^{\infty} H_n^3 \, ,
$$
which combined with (\ref{eq1.12}) yields our desired conclusion. \hfill $\Box$

\bigskip

Notice that the presence of $\underset{\overbrace{\mbox{\footnotesize${\rm int} \, \Sigma (\infty)$}}}{\bigcup} \, ({\rm int} \, \Sigma (\infty)) \times [0,\infty)$ neither interferes with the action in this lemma, nor changes its conclusions.

\smallskip

We move now to the RED elementary games. The formula (\ref{eq1.4}) is to be replaced by now by the (\ref{eq1.13}) below, which superficially may look just like one of the variants of (\ref{eq1.4}).
\begin{equation}
\label{eq1.13}
U({\rm RED}) = \{{\rm Sq} \, \cup \, [-1 \leq y \leq 1, x=\pm 1 , 0 \leq z \leq N+\varepsilon_1] \, \cup \, [-1 \leq x \leq 1 , y = \pm 1 , 0 \leq z < N \} + \sum_{n=1}^{\infty} {\rm Sq} \times z_n \, .
\end{equation}

\bigskip

\noindent This may again have variants where the RED ideal Hole is replaced by an ideal arc and (1.4.1) is, generally speaking, violated now. In the generic case, explicitly written down in (\ref{eq1.13}), the lateral pieces $(-1 \leq y \leq 1 , x = \pm 1 , 0 \leq z \leq N+\varepsilon_1)$ are vertical piece of BLACK walls $(W_{(\infty)} ({\rm BLACK}))$, while the $(-1 \leq x \leq 1 , y = \pm 1 , 0 < z < N)$ are BLUE/RED infinite staircases stretching inside $[0 \leq z < N)$.

\bigskip

IF the (1.4.1) would hold in our RED context too, but generically speaking it does not, then the lateral walls in our formula (\ref{eq1.13}) would be like in the Figure 1.3, and make sense already at the level of $X^2$. Now when we move from $X^2$ to $fX^2$, then the Figure 1.3 should be completed with the items below.

$$
\includegraphics[width=12cm]{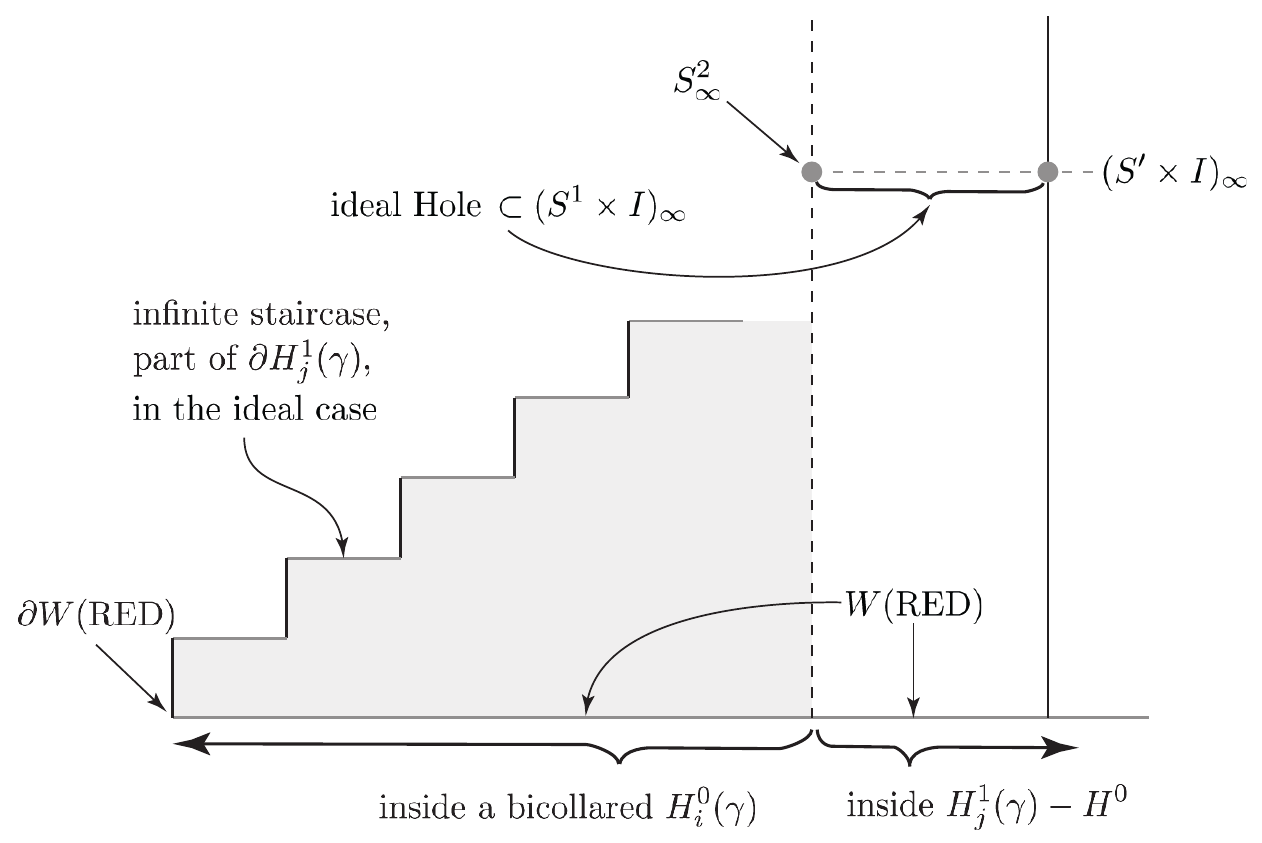}
$$
\label{fig1.3}

\centerline {\bf Figure 1.3.} 

\smallskip

\begin{quote}
Very schematical view of the lateral walls for $U^2({\rm RED})$ in the ideal, highly non-generic case when the (1.4.1) would be verified in the RED situation. The bicollared handles $H_j^1 (\gamma)$, $H_i^0 (\gamma)$ correspond to some bona-fide handles $h_j^1 , h_i^0 \subset \widetilde M(\Gamma)$. Here the $W({\rm RED})$ is the innermost RED level of our $H_j^1 (\gamma)$. Inside the infinite staircase, the vertical arcs are BLUE and the horizontal ones are RED.
\end{quote}

\bigskip

\noindent (1.13.1) \quad The infinitely many 2-handles $\underset{1}{\overset{\infty}{\sum}} \, {\rm Sq} \times z_n$. These can be the continuation of the horizontal red steps of our staircase OR also traces of other $W({\rm RED})$'s coming from the other side of the staircase and crossing it.

\bigskip

\noindent (1.13.2) \quad A lot of BLUE WALLS going through the shaded area, continuations of the BLUE horizontal steps, and others too. It should be understood that, {\ibf before} any RED game can start, these BLUE pieces have to be demolished by other preliminary games, possibly infinitely many of them. We will call this the {\ibf preliminary cleaning} operation, which corresponds to the shaded area from the Figure~1.3. All this was in the ideal case.

\bigskip

In the generic real life case of the RED game, the (1.4.1) is violated. Then, the clean situation depicted in the Figure 1.3 is to be changed as follows.

\smallskip

Corresponding to $h_j^1$ there are now infinitely many bicollared handles $H_j^1 (\gamma_1) , H_j^1 (\gamma_2) , \ldots$ each attached to some bicollared $H_i^0 (\gamma_1) , H_i^0 (\gamma_2) , \ldots$ (and to $H_i^0 (\gamma_1)^* , H_i^0 (\gamma_2)^* , \ldots$, at the other end too). The $H_j^1 (\gamma_k)$'s come with disjoined $\partial H_j^1 (\gamma_k)$'s which when $k_2 > k_1$ come closer and closer to $S_{\infty}^2$ and which are such that $\underset{k=\infty}{\lim} \, \partial H_j^1 (\gamma_k) \subset S_{\infty}^2$. Each $H_j^1 (\gamma_k)$ has its own innermost $W_k ({\rm RED})$, and these come closer and closer to $(S^1 \times I)_{\infty}$ when $k_2 > k_1$, so that we also get
$$
\lim_{k=\infty} W_k ({\rm RED}) \ (= \, \mbox{innermost NATURAL level of} \ H_j^1 (\gamma_k)) = (S^1 \times I)_{\infty} \, .
$$

These $W_k ({\rm RED})$'s are 2-by-2 disjoined, with their $\partial W_k$'s coming closer and closer to $S_{\infty}^2$, as $k$ increases and we have $\underset{k=\infty}{\lim} \, \partial W_k \subset S_{\infty}^2$.

\smallskip

So far this is {\ibf NOT} yet a $U^2 ({\rm RED})$, but out of the infinite maze of $H_i^0 (\gamma_k)$'s and $H_j^1 (\gamma_k)$'s, with $k \to \infty$, which we have described, one can extract a clean $U^2({\rm RED})$, on the lines of (\ref{eq1.13}), (1.13.1), (1.13.2), by proceeding as follows.

\smallskip

Consider, to begin with, the first two $W_k({\rm RED})$'s, these are the $W_1$ and $W_2$ in the simplest pristine case. Inside $h_i^0$ one can find (inside the corresponding complete $fX^2$ picture) a finite RED/BLUE staircase $A_1$ which has the following features

\medskip

1) The $A_1$, which might start with a collar of $\partial W_1$ inside $W_1$, joins $\partial W_1$ to $\partial W_2$; similarly the $A_1^*$ in $(h_i^0)^*$.

\medskip

2) At the level of $h_i^0 \cup h_j^1 \cup (h_i^0)^*$ the embedded surface
$$
A_1 \cup W_1 \cup A_1^* \cup W_1^*
$$
encloses a space homeomorphic to $(S^1 \times I) \times [0,1]$ inside which the pieces of $fX^2$ which may be found, are of the following kinds:

\medskip

2.1) Pieces of $W({\rm BLUE})$'s, to be killed by a preliminary cleaning operation like in (1.13.2), before any RED game can start.

\medskip

2.2) Pieces of $W_{(\infty)} ({\rm BLACK})$'s. These are actually necessary for the preliminary cleaning above. Out of them, on par with the $A_1, A_2,\ldots$ we start building the other part of the lateral walls of the $U^2 (B)$'s of the preliminary cleaning, possibly infinite RED/BLACK staircases.

\bigskip

Next we go to $W_2$ and $W_3$, for which we find a finite RED/BLUE staircase $A_2$, analogous to $A_1$. This continues indefinitely, until we build the
$$
\{(-1 \leq x \leq 1 , y = \pm 1, 0 \leq z < N) \subset \{\mbox{lateral surface of} \ U^2 ({\rm RED}) / (\ref{eq1.13})\}\} 
$$
$$
= (A_1 \cup A_2 \cup A_3 \cup \ldots) + (A_1^* \cup A_2^* \cup A_3^* \cup \ldots) \, ,
$$
which is the $[-1 \leq x \leq 1 , y = \pm 1 , 0 \leq z < N]$ in (\ref{eq1.13}), the $[-1 \leq y \leq 1 , x = \pm 1 , 0 \leq z \leq N + \varepsilon_1]$ being provided, in the clean pristine case, by the $W_{(\infty)} ({\rm BLACK})$ resting on $W_1$, like the $W_{\pm}$ in Figure 3.1 below.

\smallskip

With all this, there is now an area ${\mathcal A} \subset h_i^0$, contained between $\partial H_j^1 (\gamma_1) \cap h_i^0$ and our newly created clean $U^2 ({\rm RED})$.

\smallskip

I make now the following.

\bigskip

\noindent {\bf Claim (1.14).} One can break $fX^2 \cap {\mathcal A}$ into infinitely many $U^2(B)$'s and $U^2 ({\rm RED})$'s, each of them corresponding to an arc in $S_{\infty}^1$, and here the same given arc may parametrize several, possibly infinitely many, such $U^2$'s. With this, adjacent to our clean $U^2 ({\rm RED})$, inside ${\mathcal A}$ there are infinitely many elementary BLUE and RED games to be played. These additional games are independent from the ``main'' RED game which corresponds to the clean $U^2 ({\rm RED})$ which we have just constructed. \hfill $\Box$

\bigskip

The paradigm for the RED game is given by the next

\bigskip

\noindent {\bf RED Lemma 1.2.} {\it We will state our lemma for the clean $U^2 ({\rm RED})$ (let's say the one in {\rm (\ref{eq1.13})}), but it is valid for the other {\rm RED} games produced by {\rm (1.14)} too.}
\begin{enumerate}
\item[1)] {\it Everything said, in the context of the {\rm BLUE} elementary game, from {\rm (\ref{eq1.5})} up to {\rm (\ref{eq1.10})} included, with the analogue of the Figure {\rm 1.2} included too, remains valid for the {\rm RED} elementary game. But notice here that, for the {\rm (\ref{eq1.6})} to be valid in the {\rm RED} case, the preliminary cleaning operation which kills unwanted pieces of $W({\rm BLUE})$'s, is necessary. We also have again the analogue of {\rm (1.6.1)}, of course. The infinite symphony of games has to be played in a precise order, rather than all simultaneously.}
\item[2)] {\it The analogue of Lemma {\rm 1.1} is valid for the elementary {\rm RED} games too.}
\end{enumerate}

\bigskip

We move finally to the BLACK elementary game. We will have now  two kinds  of complication: certainly, like before, the kind of complication we had when going from BLUE to RED, i.e. the issue of going from an infinite messy picture to the single clean $U^2$. But now, on top of that, we also have the immortal singularities $S \subset \Theta^3 (fX^2)$ too, making that the analogue of $U^3$, even without the bowls ${\mathcal B}$, fails to be now a smooth 3-manifold. Figure 1.4 which should be compared to a detail of the Figure 1.1 in \cite{39}, displays a piece of $W({\rm BLACK})$. With the modification with respect to the Figure 1.2, which the Figure 1.6 below may suggest, this piece of $W({\rm BLACK})$ will generate the Sq for the $U^2 ({\rm BLACK})$, which is still to come.

\smallskip

Generically, each $W({\rm BLACK})$ contains exactly two $p_{\infty\infty} (S)$'s, but in order to simplify the exposition we will pretend that there is exactly one.

$$
\includegraphics[width=15cm]{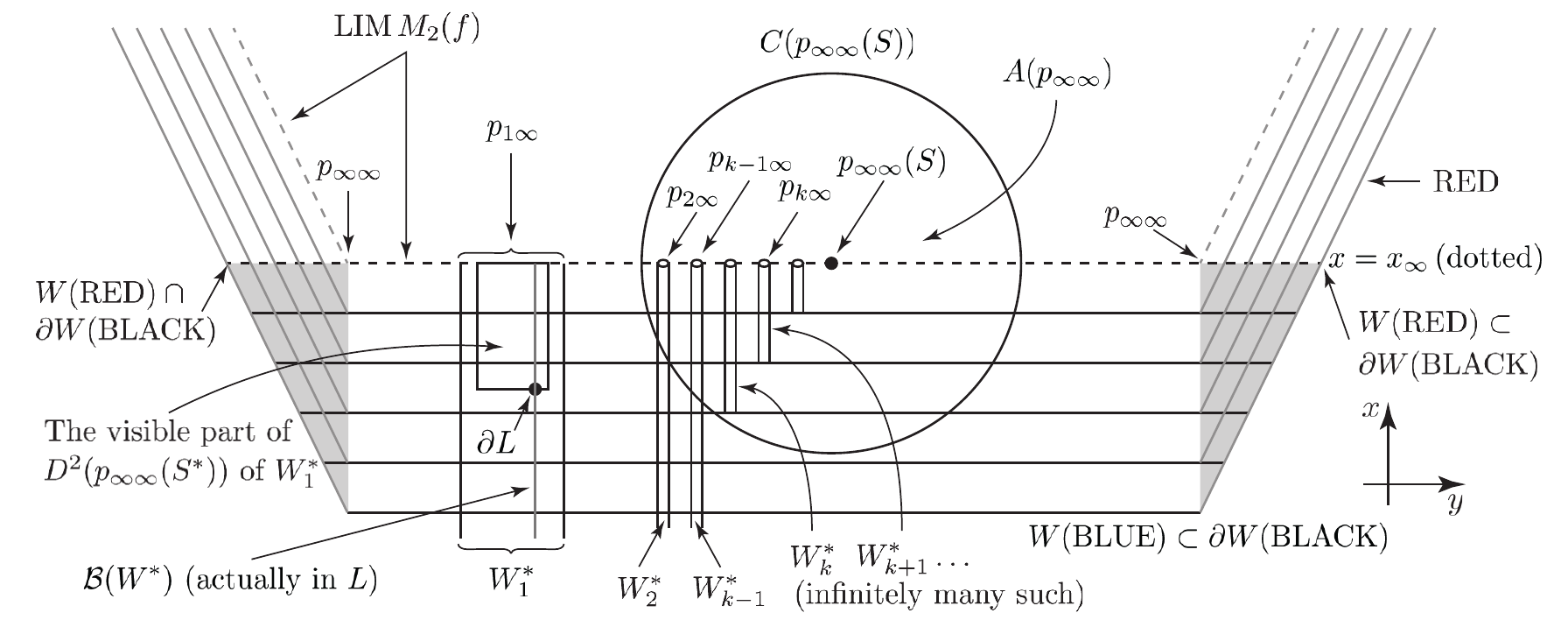}
$$
\label{fig1.4}

\centerline {\bf Figure 1.4.} 

\smallskip

\begin{quote}
We see here a typical $W({\rm BLACK})$, assumed complete and only such will take active part in the BLACK games. The $A(p_{\infty\infty})$ (like in the Figure 2.2 from \cite{39}) is the $\{$piece of $W({\rm BLACK})$ inside the circle $C(p_{\infty\infty} (S))\} - \{ p_{\infty\infty} (S)\}$. The $p_{\infty\infty}$'s are immortal singularities of $fX^2$. The parts of the dual $W({\rm BLACK})_n^*$'s which are beyond these $p_{n\infty}$'s do no longer interact with our $W({\rm BLACK})$ (at the present level $fX^2$) and they will be ignored for a while but see then the Figure 1.5 too.  At $x=x_{\infty}$ our present space $\Theta^3 (fX^2)_{\rm II}$ is traintrack and the coordinate half-line $x > x_{\infty}$ bifurcates into an $x(W)$, coordinate of our drawing, and an $x(W^*)$ invisible here. This being said, the present $D^2 (p_{\infty\infty} (S^*))$ and ${\mathcal B} (W_1^*)$, are the same as in the Figure 1.6 below, which lives in the plane $(x = x_{\infty} , y , z)$.
\end{quote}

\bigskip

There is a collar of $\partial W({\rm BLACK})$ inside $W({\rm BLACK})$, which we will denote by $[\partial W({\rm BLACK}) , {\rm LIM} \, M_2 (f)$ $\cap \, W({\rm BLACK})]$, and here ``${\rm LIM} \, M_2 (f) \cap W({\rm BLACK})$'' means the dotted hexagon with six vertices $p_{\infty\infty}$ housed inside $W({\rm BLACK})$. The collar in question has three kinds of parts, and they should be readable in our Figure 1.4, giving a decomposition
\setcounter{equation}{14}
\begin{equation}
\label{eq1.15}
[\partial W({\rm BLACK}) , {\rm LIM} \, M_2 (f) \cap W({\rm BLACK})] = \{\mbox{(three) purely RED parts}\} 
\end{equation}
$$
\cup \, \{\mbox{(three) purely BLUE parts}\} \cup \{\mbox{(six) mixed, shaded parts}\} \, .
$$
So, we will focus now on a (not yet explicitly defined) elementary BLACK game, when our present $W({\rm BLACK}$, complete) is to play the role of Sq, in the not yet written down analogues of the formula (\ref{eq1.4}) and (\ref{eq1.13}). Our $U^3 ({\rm BLACK})$ will be defined by a third formula
\begin{equation}
\label{eq1.16}
U^2 ({\rm BLACK}) = ({\rm Sq} \cup \partial \, {\rm Sq} \times [0 \leq z < N)) + \sum_{n=1}^{\infty} {\rm Sq} \times z_n \, ,
\end{equation}
with an ideal BLACK Hole living at $z=N$ and where the following things should happen. The $\partial \, {\rm Sq} \times [0 \leq z < N]$ is now the union of six infinite staircases BLACK/BLUE and BLACK/RED and then also, in order to take care of the immortal singularities visible in Figure 1.4, we split away piece $A(p_{\infty\infty})$ from the rest of $W({\rm BLACK})$, and define
\begin{equation}
\label{eq1.17}
U^3 \mid {\rm Sq} \, ({\rm BLACK}) = \left(W({\rm BLACK}) \times [-\varepsilon,\varepsilon] - A(p_{\infty\infty})(S) \times [-\varepsilon,\varepsilon] \right) \cup
\end{equation}
$$
\underset{\overbrace{\mbox{\footnotesize$C(p_{\infty\infty})(S)$}}}{\cup} D^2 (p_{\infty\infty}(S)) \times \left[ - \frac\varepsilon4 , \frac\varepsilon4 \right] \subset U^3 ({\rm BLACK}) \, .
$$

Notice that this corresponds to what $W({\rm BLACK})$ anyway becomes via the basic step (1.15) of \cite{39}. There is actually also a clear $2^{\rm d}$ counterpart to (\ref{eq1.17}), and {\ibf that} is the Sq occurring in (\ref{eq1.16}). The reader should also be warned that, in the BLACK case, the passage from $U^2$ to $U^3$ is less simple-minded than in the BLUE or RED cases, involving among other things, (deletions) $+$ (additions) $+$ (splittings), to be described below.

\smallskip

But let us assume temporarily, that we are in the ideal case when, in the style of (1.4.1), we are in the possession of an (\ref{eq1.16}), which, ideally, pre-exists at level $X^2$.

\smallskip

Our discussion is at level $fX^2$, hence the preliminary cleaning steps mentioned below. Also we are now without any other piece of unwanted infinite staircase in the way. We are still not ready for the BLACK game. Some preliminary cleaning is necessary first. In order, this is:
\begin{enumerate}
\item[i)] Via BLUE games kill the unwanted pieces of BLUE walls in the purely BLUE and the mixed pieces of (\ref{eq1.15}).
\item[ii)] Then, via RED games kill the remaining unwanted pieces of RED walls inside the purely RED AND the mixed pieces of our same (\ref{eq1.15}).
\end{enumerate}

\smallskip

Forgetting for the time being about the bowls ${\mathcal B}$, we will define our $U^3 ({\rm BLACK}) \subset \Theta^3 (fX^2)_{\rm II}$ as a {\ibf smooth} $3^{\rm d}$ branch of a larger $3^{\rm d}$ train-track manifold. Let us say that $U^3 ({\rm BLACK})$ is defined by a formula like (\ref{eq1.5}), where now $A(p_{\infty\infty}) \times [-\varepsilon \leq z \leq \varepsilon]$ is {\ibf deleted} and where the compensating 2-handle $D^2 (p_{\infty\infty} (S)) \times \left[ - \frac\varepsilon4 , \frac\varepsilon4 \right]$ is {\ibf added} instead. To this description, we also give the following modulations.

\bigskip

\noindent (1.17.1) \quad Notice, to begin with, that among the infinitely many effective intersections $W^* \cap W({\rm BLACK})$ (all of them stopping at their corresponding immortal singularity, see Figure 1.4) all except finitely many are completely inside $A(p_{\infty\infty} (S))$. Let us say, and this will be now a conventional notation, which will simplify the exposition, that the $W_1^* , W_2^* , \ldots$ which are dual and transversal to our $W = W({\rm BLACK})$, divide into three disjoined categories, as follows:
\begin{enumerate}
\item[a)] The $W_1^*$, but there can be finitely many such, which does not touch $A(p_{\infty\infty}) \supset C (p_{\infty\infty} (S))$.
\item[b)] The $W_2^* + W_3^* + \ldots + W_{k-1}^*$ which do touch $A(p_{\infty\infty})$ and which also cross $C(p_{\infty\infty}(S))$. See for all this Figure 1.4.
\item[c)] The infinite rest, i.e. $W_k^* + W_{k+1}^* + \ldots$ which are such that $W^* \cap W \subset A(p_{\infty\infty})$, never making it to $\partial W$.
\end{enumerate}

\bigskip

With this, at level $\Theta^3 (fX^2)_{\rm II}$ (or $\Theta^3 (fX^2)_{\rm I}$), the hypersurface $C(p_{\infty\infty} (S)) \times (-\infty < z < +\infty)$, see here Figure 1.4, {\ibf splits} (abstractly), each of the $W_2^* , W_3^* , \ldots , W_{k-1}^*$ into a piece $W^* (A(p_{\infty\infty}))$ which does not intersect with $U^3 ({\rm BLACK})$, and a piece $W^* (\mbox{non-$A(p_{\infty\infty})$})$, which does. The factor $-\infty < z < +\infty$ occurring in the splitting surface above, goes transversally through the plane of Figure 1.4; it is {\ibf NOT} the $z$-coordinate in the Figure 1.6 below. The splitting surface certainly goes through $\Sigma (\infty)$. But we are only focusing on the effect of this {\ibf abstract} splitting on the $W^*$'s. Its interaction with $\underset{R_0}{\sum} \, {\rm int} \, R_0 \times [0,\infty)$ can be safely ignored, it is without consequence on our conclusions.

\bigskip

\noindent (1.17.2) \quad When it will come to the bowls ${\mathcal B} \subset U^3 ({\rm BLACK})$, the idea is now the following. The ${\mathcal B} + \underset{1}{\overset{\infty}{\sum}} \, \partial H_n^3$ will {\ibf NOT} use the $A(p_{\infty\infty})$ but they will {\ibf use the 2-handle} $D^2 (p_{\infty\infty}(S))$ instead, as the Figure 1.6 suggests us to do. But before this idea can actually be implemented, we need some additional steps.

\bigskip

\noindent (1.17.3) \quad This is a reminder: our $\Theta^3 (fX^2)_{\rm II}$ is a train-track manifold which contains, among others, branches at
$$
\sum_{p_{\infty\infty} (S)} C (p_{\infty\infty} (S)) \times [-\varepsilon , \varepsilon] \, .
$$
Here, our $U^3 ({\rm BLACK})$ as defined so far, uses exactly two branches out of the three possible ones (see (\ref{eq1.17})).

\bigskip

\noindent (1.17.4) \quad In the conditions of (1.17.1) and of the Figure 1.4, we impose the following things. At the immortal singularities $p_{2\infty} , p_{3\infty}, \ldots $, created at the b) $+$ c), the $W = W({\rm BLACK})$ is {\ibf overflowing}, while the $W_2^* + W_3^* + \ldots + W_{k-1}^* + W_k^* + \ldots$ are all {\ibf subdued}. At the $p_{1\infty}$, created by a) (and in real life this corresponds not just to one, but to finitely many immortal singularities), $W_1^*$ overflows and $W({\rm BLACK})$ is subdued. [The notions of ``overflowing'' and ``subdued'' have been defined in \cite{39}; see, in particular, formulae (2.10.1) and (2.10.2) and the claim here is that our present (1.17.4) is compatible with (2.10.1), (2.10.2) in \cite{39}. We do not care if they are not implemented by the specific trick from the Figure 5.2 in \cite{39}; that was just an illustration.]

\bigskip

The next item is a consequence of the present one.

\bigskip

\noindent (1.17.5) \quad The various, infinitely many thickened disks $A(p_{\infty\infty}) \subset \Theta^3 (fX^2)_{\rm II}$ are 2-by-2 disjoined. See, at this point, the Figures 1.4 and 1.6, and also the Figures 3.4, 3.5 in Section III, which complete them. Figure 1.6 illustrates well the stated fact.

\bigskip

Since our $W = W({\rm BLACK})$, which is concerned by the Figures 1.4 and 1.6 is subdued with respect to the overflowing $W_1^*$ and also overflowing with respect to the subdued $W_2^* + W_3^* + \ldots$, the $A(p_{\infty\infty}) + D^2 (p_{\infty\infty} (S))$ of $W_1^*$ occur in the Figures 1.4 $+$ 1.6, while those of $W_2^* + W_3^* +\ldots$ do not.

\smallskip

Provided now that the preliminary cleanings mentioned at i) $+$ ii) above have been performed, here is the list of spots where the $U^3 ({\rm BLACK})$ at least as defined so far, communicates with the outside world and to this, the analogue of (1.6.1) is to be added too.

\bigskip

\noindent (1.18) \quad $\{$The SPLITTING SURFACE $C(p_{\infty\infty} (S) \times [-\varepsilon , \varepsilon]$, via which our $U^3 ({\rm BLACK})$ communicates with the deleted $A(p_{\infty\infty})\} + \{$just like in (\ref{eq1.6}), the other $-\varepsilon$ side. {\ibf But} the piece $\{ z > N \}$ which had occurred in (\ref{eq1.6}) is now unexistant; from the viewpoint of our $U^3 ({\rm BLACK})$ the $z=N$ is at infinity$\} + \{$on both $\pm \, \varepsilon$ sides, at the level of Sq itself, our $U^3 ({\rm BLACK})$ is in contact with $W_1^* , W_2^*$ (non $A(p_{\infty\infty})) , \ldots , W_{k-1}^* ({\rm non} \ A(p_{\infty\infty}))$. Here, the decomposition $W^* = W^* (A(p_{\infty\infty})) \cup W^* ({\rm non} \ A(p_{\infty\infty}))$ is defined in the formula (1.17.1). The arcs of type $[\alpha,\beta]$ or $[\gamma,\delta]$ from the Figure 1.6, when on the $A(p_{\infty\infty})$ side, are communications of $A(p_{\infty\infty})$ with the outside world, and {\ibf not} communications of $U^3 ({\rm BLACK})$. The $[\delta , \alpha]$, $[\gamma , \beta]$ on the $W_1^*$ side {\ibf are} communications of $U^3 ({\rm BLACK})$. We have, for {\ibf their} rectangle $[\alpha \, \beta \, \gamma \, \delta]$
$$
[\alpha \, \beta \, \gamma \, \delta] = U^3 ({\rm BLACK})_W \cap \{{\rm The} \ A(p_{\infty\infty})_{W_1^*}, \mbox{which is deleted from} \ U^3 ({\rm BLACK})_{W_1^*} \} \, .
$$
Outside of Sq, and again in both sides $\pm \, \varepsilon$, our $U^3 ({\rm BLACK})$ also communicates with $\underset{n=k}{\overset{\infty}{\sum}} \, W_n^*$ too; see here the legend of Figure 1.5.$\}$. This ENDS formula (1.18).

\bigskip

We will give now a more detailed description of the interactions $U^3 ({\rm BLACK}) \cap W^*$.

\smallskip

Notice, to begin with, that starting with our $W = W_0 \equiv \{$our $W({\rm BLACK})\}$, there is a whole infinite family of parallel walls $W({\rm BLACK \ complete})$, parallel to $W_0$ and converging to the ideal BLACK Hole, namely
\setcounter{equation}{18}
\begin{equation}
\label{eq1.19}
W_1 , W_2 , W_3 , \ldots \ \mbox{{\ibf and} the} \ {\rm Sq} \times z_n \ \mbox{of our} \ U^3 ({\rm BLACK}) \ \mbox{is a (thickened piece of)} \ W_n \, .
\end{equation}

Our $W_1^* , W_2^* , W_3^* , \ldots$ are dual not only to our initial $W_0$, but to all the other $W_1 , W_2 , W_3 , \ldots$ too. Each of our $W_1^* , W_2^* , \ldots , W_{k-1}^*$ is getting zipped, at the level of (\ref{eq1.16}), with the rest of $U^2 ({\rm BLACK})$, along a zipping path (which when considered with time ordering reversed) is starting at $p_{1\infty}$ (Figure 1.4) or at $(W_2^* + W_3^* + \ldots + W_{k-1}^*) \cap C (p_{\infty\infty} (S))$ and involving the $W_{1<i<k}^* - A (p_{\infty\infty})$. In terms of the notations of Figure 1.5, these paths go first to some $S_{\ell \leq p}$ and then further to $S_{\ell}^*$.

$$
\includegraphics[width=9cm]{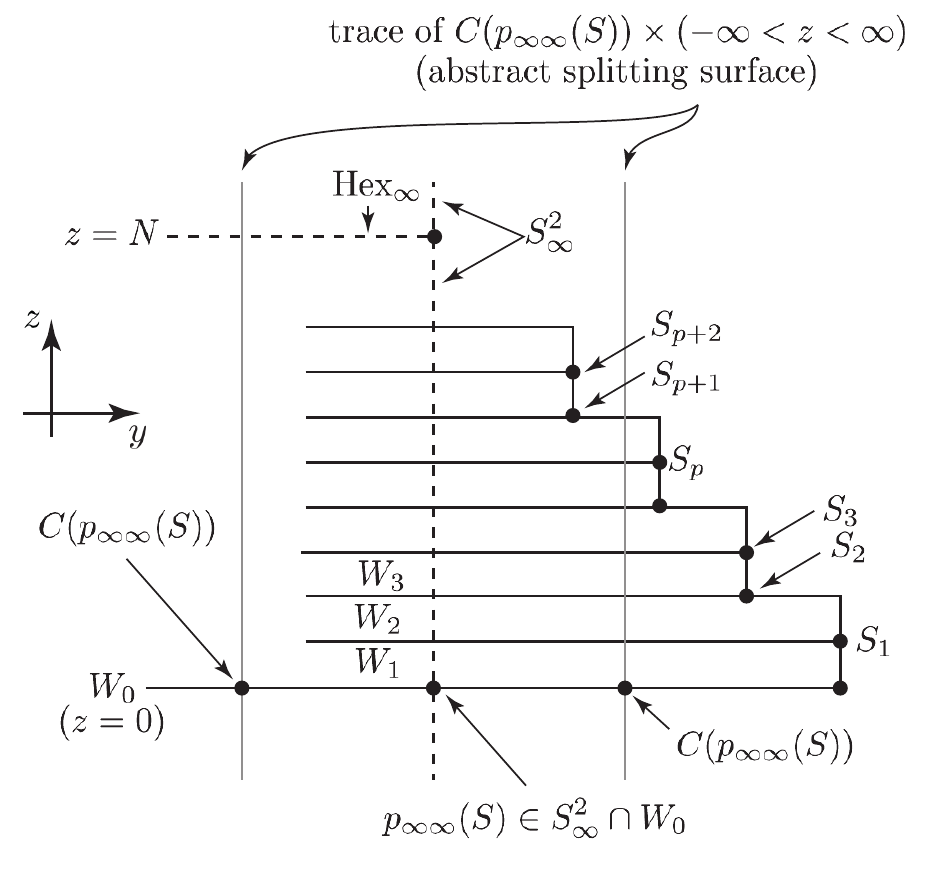}
$$
\label{fig1.5}

\centerline {\bf Figure 1.5.} 

\smallskip

\begin{quote}
In this figure, which is in the style of Figure 1.3, we have suggested a BLUE/BLACK infinite staircase, part of $\partial \, {\rm Sq} \times [0 \leq z < N]$ in (\ref{eq1.16}). We have suggested as fat points, mortal singularities which occur normally in the zipping of $X^2 \to fX^2$, at intermediary stages of the zipping in question. The $S_1 , S_2 , \ldots , S_p$ involve $\{W_1 , W_2 , \ldots\}$ and $\{W_1^* , \ldots , W_{k-1}^*\}$, while the $S_{p+1} , S_{p+2} , \ldots$ involve $\{W_1 , W_2 , \ldots\}$ and $\{W_k^* , W_{k+1}^* , \ldots\}$. When we go to the complete $fX^2$, each mortal singularity $S_n$ is replaced after a short zipping, by an immortal singularity which we call $S_n^*$, and which involves the same pair $(W,W^*)$. In the present figure, the horizontal walls are BLACK, while the vertical ones are BLUE.
\end{quote}

\bigskip

At this point, we will make the following CHANGES concerning the definition of $U^3({\rm BLACK})$, as presented so far. These changes will complete and/or supersede when necessary, the (1.18) above.

\bigskip

\noindent (1.20.1) \quad In the context of (1.18), when it comes to the communication with the outer world along the internal $+\varepsilon$ side, $U^3 ({\rm BLACK})$ communicates exactly via {\ibf ALL} the $W_1^* , W_2^* , \ldots$. For the $W_1^* , W_2^* , \ldots , W_{k-1}^*$, this communication starts on the Sq from (\ref{eq1.16}) $+$ (\ref{eq1.17}), and then continues along $\partial \, {\rm Sq} \times [0 \leq z < N]$, until we reach an intermediary singularity $S_p$, (see here the Figure 1.5) and next, along the $W_p \sim {\rm Sq} \times z_{n(p)}$, to the final immortal singularity $S_p^*$. Of course, ``$S_p$'' may read here $S_1 , S_2 , \ldots , S_p$. As the Figure 1.5 suggests, the $W_k^* , W_{k+1}^* , \ldots$ also have contacts with the 2-handles ${\rm Sq} \times z_n$ along arcs of type  $[S_p , S_p^*]$. The point is that, with these things our 2-handles $\overset{\infty}{\underset{1}{\sum}} \, {\rm Sq} \times z_n$ are not clean as they should be, when we will want to apply the handle-cancellations from the BLUE or RED games (see Lemma 1.1), in the BLACK case. This motivates the next change.

\bigskip

\noindent (1.20.2) \quad In order to free the 2-handles $\overset{\infty}{\underset{1}{\sum}} \, {\rm Sq} \times z_n$ from (1.20.1) the immortal singularities $S_1^* , S_2^* , \ldots$ living in their middle, we unzip $W_n$ and $W^*$ along $[S_n , S_n^*]$ creating a mortal singularity at $S_n$. The $\Theta^4 (\Theta^3 (fX^2) , {\mathcal R})_{\rm II}$ in (\ref{eq1.2}) does not feel the difference. We get now a {\ibf singular} $U^3 ({\rm BLACK})$ even before the bowls ${\mathcal B}$ are thrown into the game. At each $S_n$, our $U^3 ({\rm BLACK})$ has now a mortal singularity with one branch $\{$the previous smooth $U^3 ({\rm BLACK})\}$ and another very small $W_n^*$ branch, via the outer boundary of which our new $U^3 ({\rm BLACK})_{\rm singular}$ also communicates with the outer world. But it does not communicate, any longer with the outer world through the 2-handles themselves. Figure 1.6 is the analogue of the Figure 1.2 for our final $U^3 ({\rm BLACK})$.
$$
\includegraphics[width=148mm]{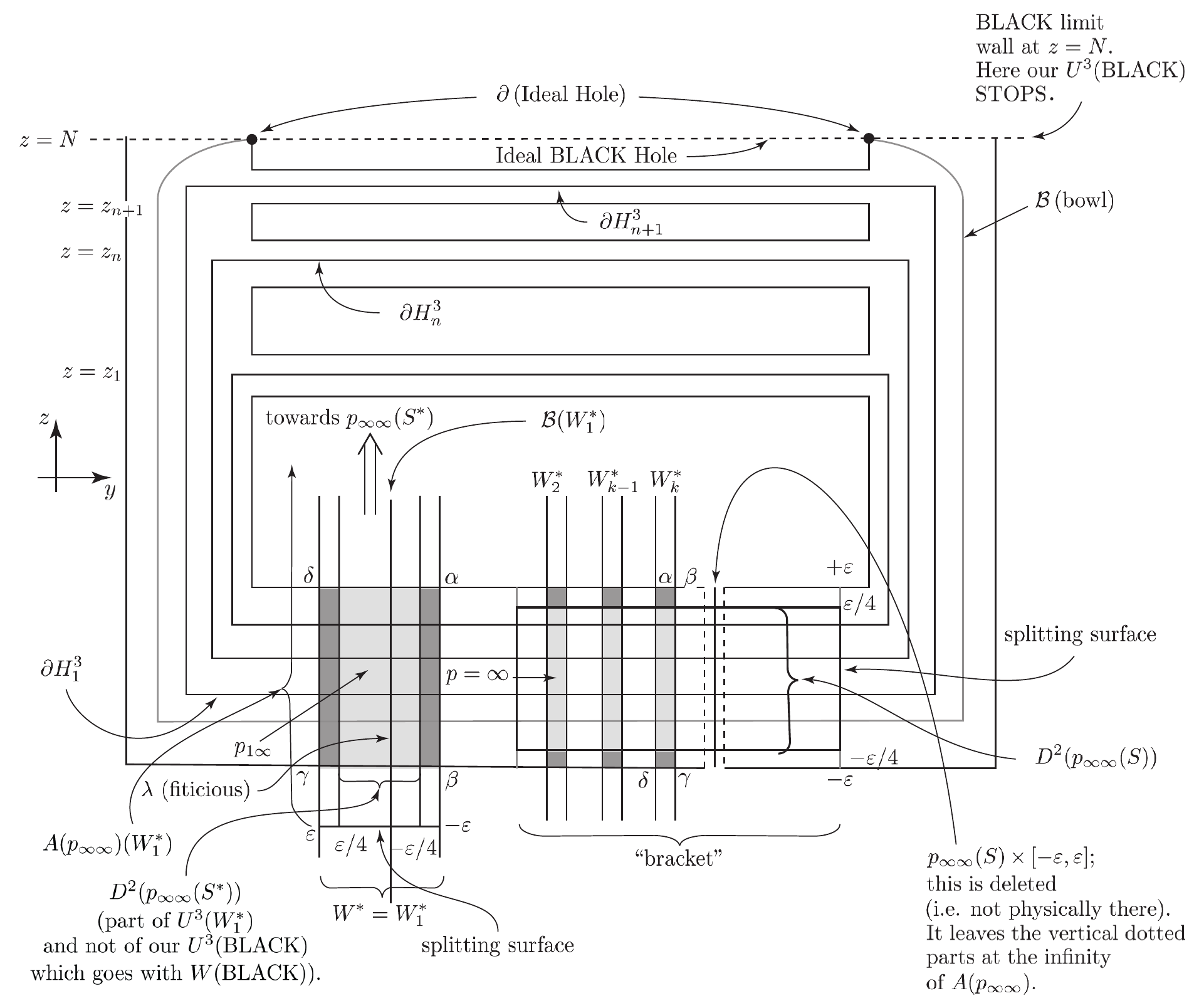}
$$
\label{fig1.6}
\centerline {\bf Figure 1.6.} 

\smallskip

\begin{quote}
The $U^3 ({\rm BLACK})$, with ${\rm Sq} \subset W({\rm BLACK})$. We are here at $x=x_{\infty}$. Here, at the ``bracket'', $C(p_{\infty\infty} (S)) \times [-\varepsilon,\varepsilon]$ splits the $\Theta^3 (fX^2)_{\rm II} \mid \{$our $W({\rm BLACK})\}$ into three branches: The outer part, which belongs to ${\rm Sq} \subset U^3 ({\rm BLACK})$, the $A(p_{\infty\infty}) \times [-\varepsilon,\varepsilon]$ part  (corresponding to our $W({\rm BLACK})$) and then  also the 2-handle $D^2 (p_{\infty\infty}(S)) \times \left[ -\frac\varepsilon4 , \frac\varepsilon4 \right]$. Remember here that the $\Theta^3 (fX^2)_{\rm II}$ is a train track manifold. The point ``$\lambda$'' is ficticious and it has been drawn in only for explanatory purposes. At our present level $x=x_{\infty}$, the ${\mathcal B} (W)$ goes through
$$
\{A(p_{\infty\infty}) (W_1^*) + D^2 (p_{\infty\infty}(S)) \} \, ,
$$
while the ${\mathcal B} (W_1^*)$ goes through $D^2 (p_{\infty\infty}(S)) (W_1^*)$, which is disjoined from the $\{\ldots\}$ above. The coordinates of $\lambda$ are
$$
\lambda = (x=x_{\infty} , y_0 , z_0) \equiv \{\mbox{the $z$ of $W({\rm BLACK})$ in the Figure 1.4, i.e. our $W$\}} \, .
$$
Our present $\lambda$ is a reminder of the physical point $\partial L$ from the Figure 1.4, which comes with the coordinates
$$
\partial L = \{ x  = x_0 \, (\partial L) < x_{\infty} , y_0 , z_0 \} \, .
$$
The only $A(p_{\infty\infty})$ with which our present $U^3 ({\rm BLACK})$ has an intersection, is the
$$
\{{\rm piece} \ [\alpha \, \beta \, \gamma \, \delta] \ \mbox{(L.H.S. of our figure)} \ \subset A(p_{\infty\infty}) (W_1^*)\} \subset U^3 ({\rm BLACK}) \ ({\rm of} \ W ({\rm BLACK})) \, .
$$

\smallskip

\begin{multicols}{2}
$\includegraphics[width=35mm]{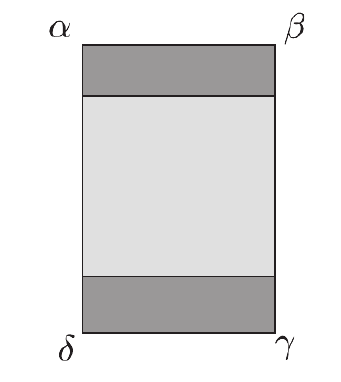}$ \\ Legend: This is an immortal singularity $S(W({\rm overflowing}) \cap W^* ({\rm subdued}))$ of $\Theta^3 (fX^2)_{\rm II}$. Along the simply shaded area, the $W^*$ looks superposed with the $D^2(W) \times \left[ -\frac\varepsilon4 , \frac\varepsilon4 \right]$, but this is just optical illusion, the $S$'s are always inside the $A(p_{\infty\infty})$'s, disjoined from the compensating 2-handle $D^2 (p_{\infty\infty})$. Along the two sides marked $[\alpha \, \beta], [\gamma \, \delta]$, it is $W^*$ which continues, while along $[\beta \, \gamma] , [\delta \, \alpha]$ it is the wall $W$.
\end{multicols}

\smallskip

In the present figure we are at $x=x_{\infty}$, reason for seeing the arc $p_{\infty\infty} (S) \times [-\varepsilon , \varepsilon]$ which, of course, is not physically present. The shaded areas correspond to the immortal singularities $p_{1\infty} , p_{2\infty}, \ldots$ from Figure 1.4.
\end{quote}

\bigskip

\noindent {\bf Comments concerning the Figure 1.6.} The figure in question lives at ($x=x_{\infty} , y ,z$). At the point marked $\lambda$, there is no actual contact ${\mathcal B} (W) \cap {\mathcal B} (W_1^*)$. We have there:
$$
{\mathcal B} (W) \subset A(p_{\infty\infty}) (W_1^*) \subset U^3 (W({\rm BLACK}))
$$
and
$$
{\mathcal B} (W_1^*) \subset D^2 (p_{\infty\infty} (S)(W_1^*)) \subset U^3 (W_1^* ({\rm BLACK})) \, .
$$
As far as $W({\rm BLACK})$ and its $U^3 ({\rm BLACK})$ are concerned, all the contacts $\left( {\mathcal B} + \underset{n=1}{\overset{\infty}{\sum}} \, \partial H^3_n \right) \cap A(p_{\infty\infty})$ have been transformed on the compensating handle $D^2 (p_{\infty\infty} (S)) \subset U^3 ({\rm BLACK})$. Also
$$
\underset{{\rm our} \ U^3 ({\rm BLACK})}{\underbrace{U^3 (W({\rm BLACK}))}} \, \cap \, A(p_{\infty\infty}) \ ({\rm of} \ W({\rm BLACK})) = \emptyset
$$
while
$$
U^3 (W({\rm BLACK})) \cap A(p_{\infty\infty}) \ ({\rm of} \ W_1^*) \ne \emptyset \, .
$$
What we see in the Figure 1.6 is a train-track and, importantly
$$
\left\{\mbox{all the shaded contribution $(S)$ of} \ \sum_1^{\infty} W_n^* \right\} \cap D^2 (p_{\infty\infty} (S)) = \emptyset \, ,
$$
and this equality concerns, of course $x=x_{\infty}$. The simple shading corresponds to $D^2 (p_{\infty\infty} (S))$ and its superposition with the shaded $S$'s is just an optical illusion. Also, still at $x=x_{\infty}$, the $\biggl\{\mbox{shaded contribution of}$ $\underset{n=2}{\overset{\infty}{\sum}} \, W_n^* \biggl\} \subset A(p_{\infty\infty})(W({\rm BLACK}))$, while the $\{$shaded contribution of $W_1^*\}$ is outside of $A(p_{\infty\infty})(W({\rm BLACK}))$, but inside $A(p_{\infty\infty})(W_1^*) \cap U^3 (W({\rm BLACK}))$.

\smallskip

One might have also noticed, already, that the notations of $W$ versus $W_1^*$ and $W_{i \geq 2}^*$ versus $W$, are symmetrical.

\smallskip

The additions of 3-handles, followed immediately by a cancellation of $\lambda=2$ and $\lambda=3$ handles demanded by the BLACK game for $U^3(W)$ are in no way disturbed by the dual $U^3(W^*)$'s. See also what is said below concerning (\ref{eq1.21}). So, our elementary BLACK games for $W$ and $W^*$ can be played in any order.

\smallskip

When we go outside of $x=x_{\infty}$, where our drawing lines then we find that (see here the Figure 1.4)
$$
\underset{{\mbox{\scriptsize attaching zone of the 2-handle} \atop D^2 (p_{\infty\infty}(S))({\rm of} \, W({\rm BLACK}))}}{\underbrace{\left( C(p_{\infty\infty}) \times \left[ - \frac\varepsilon4,\frac\varepsilon4\right] \right)}} \, \cap \, \sum_{i=2}^{k-1} W_i^* (\mbox{non-$A(p_{\infty\infty})$}) \ne \emptyset \, .
$$
So, far from $x=x_{\infty}$, we also find contacts
\setcounter{equation}{20}
\begin{equation}
\label{eq1.21}
\left( {\mathcal B} + \sum_{n=1}^{\infty} \partial H_n^3 \right) \cap \left[ W_1^* + \sum_{i=2}^{k-1} W_i^* ({\rm non} \ A(p_{\infty\infty})) \cap {\rm Sq} \right] \ne \emptyset \, ,
\end{equation}
which I claim to be {\ibf harmless}.

\bigskip

\noindent {\bf Remarks.} There is another alternative, avoiding the step (1.20.2). It is to insist that the cocores of the 2-handles do not touch the $W^*$'s and to work only with the cocores, not with the whole handles. At least for expository purposes, I thought the chosen variant, i.e. using (1.20.2) is smoother. \hfill $\Box$

\bigskip

Everything said so far was in the ideal case when, in the style of (1.4.1) and of the Figure 1.3, our present Figure 1.5 (and its undrawn RED zippings), concerns a single bicollared handles $H_j^2 (\gamma)$ of which $W_0$ is the unique $W(\mbox{BLACK, complete})$ (see here (1.13) in \cite{39}) and the staircase is part of
$$
\partial H_j^2 (\gamma) \cap \{{\rm adjacent} \ H_i^0 (\gamma)\}
$$
(in the RED sibling of Figure 1.5 this is then rather $\partial H_j^2 (\gamma) \cap \{{\rm adjacent} \ H_k^1 (\gamma)\}$) {\ibf AND} when no other bicollared handles perturb the clean picture which leads to Figure 1.6.

\smallskip

In the real life case, there are actually infinitely many $H_j^2 (\gamma_1) , H_j^2 (\gamma_2) , \ldots$ attached to $H_{i_1}^0 (\gamma_1) , H_{i_2}^0 (\gamma_2) , \ldots$ (three of them for each $\gamma_i$ and to $H_k^1 (\gamma_1) , H_k^1 (\gamma_2),\ldots$ (again three of them). We will denote by $W_0 (\gamma_{\ell})$ the unique $W(\mbox{BLACK complete})$ of $H_j^2 (\gamma_{\ell})$. Anyway, out of this infinite maze we have to extract now a clean picture. It may be assumed, without loss of generality that the location of the $W_0 (\gamma_1) ,W_0 (\gamma_2) , \ldots$ is such that there exists a unique BLACK limit wall to which they come closer and closer as $\gamma_1 < \gamma_2 < \ldots$ and converge to it when $\gamma_{\ell} \to \infty$. There is no harm in imposing this as a condition going with the (1.13) in \cite{39}, when considered at the target. Figure 1.7 should illustrate this.

$$
\includegraphics[width=14cm]{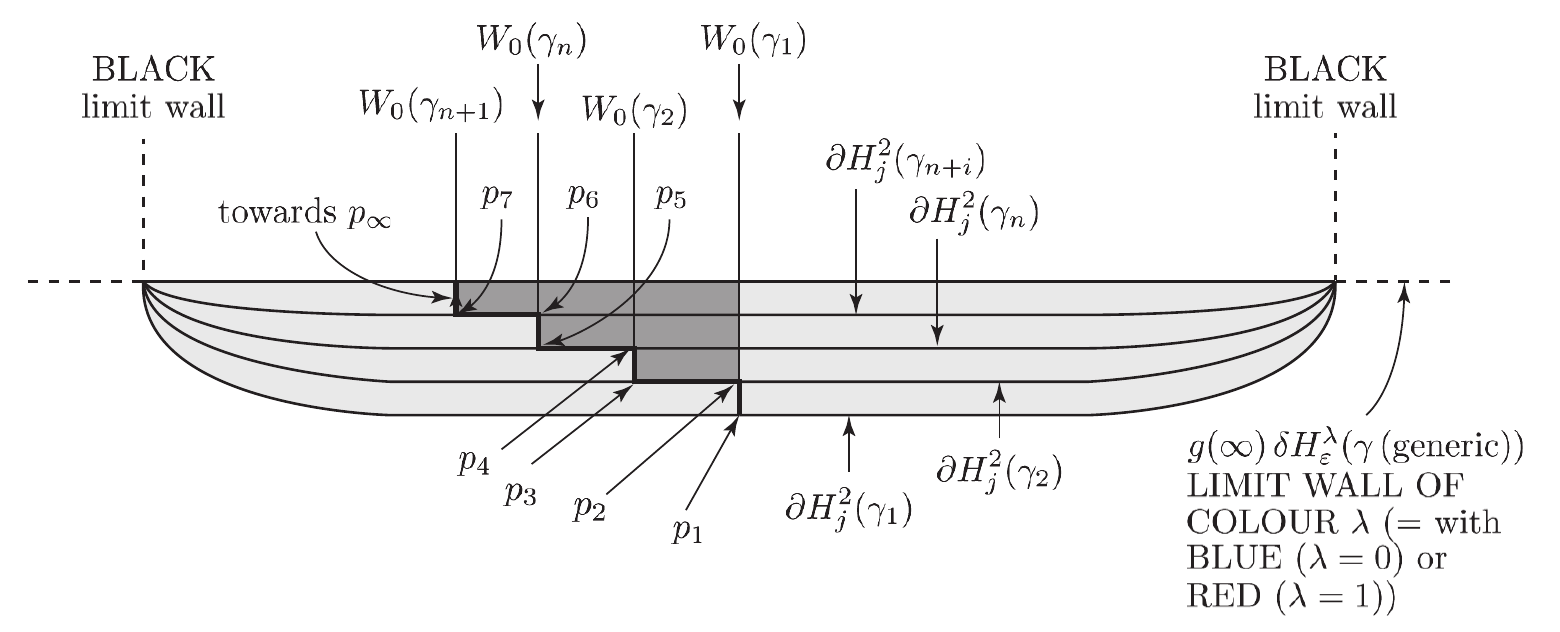}
$$
\label{fig1.7}

\centerline {\bf Figure 1.7.} 

\smallskip

\begin{quote}
We see here, at the target $\widetilde M(\Gamma)$, in the style of the Figure 2.2 from \cite{31}, the $g(\infty)$-image (see here formula (1.6) in \cite{39}) of a generic bicollared handle to which $H_j^2 (\gamma)$ is attached. We may assume that, in terms of (\ref{eq1.15}), for $W({\rm BLACK}) = W_0 (\gamma_0)$, this figure is a slice through the pure BLUE or pure RED part of the collar $[\partial W({\rm BLACK}) , {\rm LIM} \, M_2(f) \cap W({\rm BLACK})]$. The BLACK limit wall to the left of the figure {\ibf is} the ideal BLACK Hole of the clean $U^3({\rm BLACK})$ which gets created here.
\end{quote}

\bigskip

This same figure should suggest how, out of $\partial H_j^2 (\gamma_1) \cup W_0 (\gamma_1)$, $\partial H_j^2 (\gamma_2) \cup W_0 (\gamma_2) , \ldots$ we can extract a clean $U^2 ({\rm BLACK})$ (\ref{eq1.16}). We take as Sq the $\{W_0 (\gamma_1)$ modified like in (\ref{eq1.17})$\}$, as $\partial \, {\rm Sq} \times [0 \leq z < N]$ the infinite staircase suggested by $[p_1 , p_2 , p_3, p_4, p_5 , p_6 , p_7 , \ldots]$ and as 2-handles ${\rm Sq} \times z_n$ the sequence of $W_0 (\gamma_2) , W_0 (\gamma_3) , \ldots$ (starting at the staircase which occurs in fat lines). To be really OK, this still needs a preliminary cleaning of the doubly shaded region from Figure 1.7. This can be done by preliminary GAMES of Colour $\lambda$ ($< {\rm BLACK} \, (\lambda = 2)$) AND by collapsing away of pieces of $W(\mbox{BLACK NOT complete})$'s.

\smallskip

Also, independently of the main elementary Game which comes with the newly created clean $U^2({\rm BLACK})$, the simply shaded areas in the Figure 1.7 correspond to other, BLUE and RED games parametrized by ideal arcs. In the Figure~1.7, {\ibf we} have created an infinite staircase and selected appropriate pieces of successive $W(\mbox{BLACK complete})$'s, so as to create an $U^3 ({\rm BLACK})$ like in the Figure 1.6. But then, pieces of $W(\mbox{BLACK complete})$ are complete now between successive $\partial H (\gamma_n)$'s, not part of the construction above. They will have to be killed by additional, degenerate BLACK games, parametrized now by ideal arcs.

\bigskip

\noindent {\bf The BLACK Lemma 1.3.} {\it The manipulations above, concerning the Figure {\rm 1.7}, create $\infty + 1$ {\rm BLACK} games, a main one which corresponds to an ideal, $2$-dimensional {\rm BLACK} hole, and the infinitely many additional ones, each corresponding to some ideal arc. 

\smallskip

When all the $h_j^2$'s are taken into account, this is an exhaustive list of all the {\rm BLACK} games. The main black game(s) come with a figure like {\rm 1.6}, and get(s) a treatment like in {\rm (1.20.1)} $+$ {\rm (1.20.2)}. The additional {\rm BLACK} games come with a much simpler figure, very much like the Figure {\rm 1.2.}

\smallskip

The analogue of Lemma {\rm 1.1} remains valid for all these elementary {\rm BLACK} games.}

\bigskip

We consider now the bicollared handles $H_i^{\lambda} (\gamma) \subset Y(\infty)$ and also the
\begin{equation}
\label{eq1.22}
H_i^{\lambda} \equiv\{\mbox{the common $g(\infty)$-image of all the $H_i^{\lambda} (\gamma)$'s}\} \subset \widetilde M (\Gamma) \, .
\end{equation}

The $fX^2 \cap H_i^{\lambda}$, or rather their restriction to individual walls $W$ are the so-called ``complete figures'' of \cite{39}, and $H_i^{\lambda}$ corresponds to a $h_i^{\lambda} \subset \widetilde M (\Gamma)$.

\bigskip

\noindent {\bf Lemma 1.4.} {\it We can choose a unique compact wall, in each $H_i^{\lambda}$ above}
$$
W_i ({\rm COLOUR} \ \lambda) \subset fX^2 \cap H_i^{\lambda} \, , \quad \mbox{s.t.}
$$

1) {\it If $\partial H_{\alpha}^1 = H_{\beta}^0 - H_{\gamma}^0$, then $W_{\alpha} ({\rm RED})$ makes it all the way from $W_{\beta} ({\rm BLUE})$ to $W_{\gamma} ({\rm BLUE})$. We will denote, from now on by ``$\,W_{\alpha} ({\rm RED})$'', the $W_{\alpha} ({\rm RED})$ truncated by $B_{\beta}^3 ,B_{\gamma}^3$, the two $3$-balls bounded by the $2$-spheres $W_{\beta} ({\rm BLUE})$, $W_{\gamma} ({\rm BLUE})$ in their respective $H_{\beta}^0 , H_{\gamma}^0$'s. The truncated ``$\,W_{\alpha} ({\rm RED})$'' together with the $W_{\beta} ({\rm BLUE})$, $W_{\gamma} ({\rm BLUE})$ determines a $1$-handle $(D^2 \times I)_{\alpha} \subset H_{\alpha}^1$, which is attached to $B_{\beta}^3 ,B_{\gamma}^3$.}

\medskip

2) {\it For a given $H_i^2$, each of the infinitely many $X^2 \mid H_i^2 (\gamma)$'s contains exctly one $W_i (\mbox{\rm BLACK, complete})_{\gamma}$ and one of them (exactly)  will be our chosen $W_i({\rm BLACK})$. It will be assumed, again, that when
$$
\partial H_i^2 = \sum_{\alpha} H_{\alpha}^0 \cup \sum_{\beta} H_{\beta}^1 \, ,
$$
then $\partial W_j ({\rm BLACK})$ makes it all the way to
$$
\sum B_{\alpha}^3 \cup \sum (D^2 \times I)_{\beta}
$$
and, from now on we will denote by $W_j ({\rm BLACK})$, the $W_j ({\rm BLACK})$ truncated by the $3^{\rm d}$ object written above.}

\medskip

3) {\it Our choice can be made equivariantly, meaning that for each $x \in \Gamma$, when $x H_i^{\lambda} = H_j^{\lambda} (= H_{xi}^{\lambda})$ then we have
$$
x \, W_i ({\rm colour} \ \lambda) = W_j ({\rm colour} \ \lambda) \, ,
$$
an equivariance which should hold both for the untruncated and the truncated $W$'s.}

\medskip

4) {\it The $B_{\alpha}^3$ and $(D^2 \times I)_{\beta}$'s cut out of $fX^2$ finite complexes, where not only the compact $W$'s contribute, but the $W_{\infty}$'s too, call these $fX^2 \mid B_{\alpha}^3$, $fX^2 \mid (D^2 \times I)_{\beta}$. We can introduce the the {\ibf locally finite} simply-connected complex}
\begin{equation}
\label{eq1.23}
Y^2 \equiv \sum_{H_{\alpha}^0} \left( fX^2 \mid B_{\alpha}^3 \right) \cup \sum_{H_{\beta}^1} \left( fX^2 \mid (D^2 \times I)_{\beta} \right) \cup
\end{equation}
{\it $\cup \ \biggl\{ \underset{H_j^2}{\sum} \, W_j ({\rm BLACK})$, where in the cases when $W_i ({\rm BLACK}) = W_j ({\rm BLACK})^*$ (i.e. they come from two $H_i^2 , H_j^2$ like in the Figure {\rm 1.5} from {\rm \cite{39}}), then they are naturally zipped together, from the mortal singularity occurring on $\partial W_i ({\rm BLACK}) \, \cap \, \partial W_j ({\rm BLACK}) \, \cap \, \{$one of the attached $W_{\alpha} ({\rm BLUE}) = \partial B_{\alpha}^3\}$, to the corresponding immortal singularity $S(i,j) \in {\rm Sing} \, Y^2 \biggl\}$.}

\medskip

5) {\it There is a free action $\Gamma \times Y^2 \to Y^2$, which is co-compact, coming with $\pi_1 (Y^2 / \Gamma) = \Gamma$. But we certainly do {\ibf NOT} claim that $Y^2$ is {\rm QSF} and so, we cannot deduce that $\Gamma \in {\rm QSF}$ from things said so far.}

\bigskip

Notice that, with the zipping part of (\ref{eq1.23}) we get a natural inclusion $Y^2 \subset fX^2$. In terms of $X^2$ and/or of $fX^2$, our $Y^2$ is a union of walls (or pieces of walls) $W({\rm BLUE}), W({\rm RED}) , W({\rm BLACK})$ and $W_{(\infty)}({\rm BLACK})$, these last ones being caught inside the $(D^2 \times I)_{\beta}$'s or the $B_{\alpha}^3$'s. We also introduce the $3^{\rm d}$ cell-complex
\begin{equation}
\label{eq1.24}
\Theta^3 (Y^2) \equiv \Theta^3 (fX^2)_{\rm II} \mid Y^2 \subset \Theta^3 (fX^2)_{\rm II} \, ,
\end{equation} 
coming with a PROPERLY embedded (branched) surface
\begin{equation}
\label{eq1.25}
\left( {\rm int} \ \Sigma (\infty) \right) \cap \Theta^3 (Y^2) \subset \Theta^3 (Y^2) \, ,
\end{equation} 
with ${\rm int} \ \Sigma (\infty)$ like in (1.1.bis) above, i.e. $\{$the ${\rm int} \ \Sigma (\infty)$ for the $\Sigma (\infty)$ from (2.13.1) in \cite{39}, with all the contribution of the $p_{\infty\infty} (S)$'s removed, while the one of the $p_{\infty\infty} ({\rm proper})$ is left in place$\}$.

\smallskip

We will simplify the notations from (\ref{eq1.25}) into
\begin{equation}
\label{eq1.26}
\ring\Sigma (\infty) \equiv \left( {\rm int} \ \Sigma (\infty) \right) \cap \Theta^3 (Y) \, , \ \mbox{from now on.}
\end{equation}
Finally I will introduce the following subcomplex of the $\Theta^3 (fX^2)_{\rm II}$, namely
\begin{equation}
\label{eq1.27}
\Theta^3 ({\rm provisional}) \equiv \Theta^3 (Y) \underset{\overbrace{\mbox{\footnotesize$\ring\Sigma (\infty)$}}}{\cup} \ring\Sigma (\infty) \times [0,\infty) \, .
\end{equation}

\bigskip

\noindent {\bf The main multi-game Lemma 1.5.} {\it There exists an infinite sequence of elementary games, which we will call the {\ibf multi-game}}
\begin{equation}
\label{eq1.28}
\Theta^3 ({\rm old}) \equiv \Theta^3 (fX^2)_{\rm II} \overset{\mbox{\footnotesize MULTI-GAME}}{=\!\!=\!\!=\!\!=\!\!=\!\!=\!\!=\!\!=\!\!=\!\!=\!\!=\!\!=\!\!=\!\!=\!\!=\!\!\!\Longrightarrow} \Theta^3 ({\rm new}) \, ,
\end{equation}
{\it with features to be described below. It should be understood, to begin with, that as a consequence of BLUE, RED and BLACK Lemmas above, we have that the $S_u ({\rm new}) \equiv \Theta^4 (\Theta^3 ({\rm new}) , {\mathcal R}) \times B^N$ is {\rm G.S.C.}}

\medskip

0) {\it The multi-game in {\rm (\ref{eq1.28})} does not touch the $D^2 (p_{\infty\infty} (S))$'s. Also, one should read the definition of the $(N+4)$-dimensional cell-complex $S_u ({\rm new})$ above, like in the formula {\rm (\ref{eq1.3})}, namely as
$$
\{\mbox{a smooth $(N+4)$-manifold}\} + \sum_{p_{\infty\infty} (S)} \{\mbox{compensating $2$-handles of dimension $N+4$}\} \, .
$$

The MULTI-GAME leaves us with
$$
\pi_1 \, \Theta^3 ({\rm new}) = \pi_1 \, \Theta^3 ({\rm old}) = 0 \, .
$$

Out of the infinitely generated $\pi_2 \, fX^2 = \pi_2 \, \Theta^3 ({\rm old})$, the MULTI-GAME leaves only a finitely generated $\pi_2 \, \Theta^3 ({\rm new})$ alive.

\smallskip

Moreover, we will have a collapse}
\begin{equation}
\label{eq1.29}
\Theta^3 ({\rm new}) \overset{\mbox{\footnotesize collapse}}{-\!\!\!-\!\!\!-\!\!\!-\!\!\!-\!\!\!-\!\!\!-\!\!\!\longrightarrow} \Theta^3 ({\rm provisional}) \, .
\end{equation}
{\it The {\rm (\ref{eq1.28})}, the $S_u ({\rm new})$ and the {\rm (\ref{eq1.29})} are all $\Gamma$-equivariant and
$$
S_u ({\rm new}) = (S_u ({\rm new}) \diagup \Gamma)^{\sim} \, .
$$
Achieving all these things said above was, actually, the whole aim of our Multi-Game.}

\medskip

1) {\it The multi-game does not only delete things from the $\Theta^3 (fX^2)_{\rm II}$, it also adds the Bowls $\underset{n}{\sum} \, {\mathcal B}_n \times [0,\infty)$, one of them for each individual elementary game.

\smallskip

There is a PROPER map, which is injective, on each individual ${\mathcal B}$, call it
\begin{equation}
\label{eq1.30}
\sum_n {\mathcal B}_n \overset{\mathcal J}{-\!\!\!-\!\!\!-\!\!\!-\!\!\!\longrightarrow} \Theta^3 ({\rm new}) \supset \ring\Sigma (\infty) \times [0,\infty) \, ,
\end{equation}
where the $\ring\Sigma (\infty) \times [0,\infty)$ is like in {\rm (\ref{eq1.26})} above and, in the context of {\rm (\ref{eq1.30})} we also have
$$
{\mathcal J} \left( \sum_n {\mathcal B}_n \right) \cap \ring\Sigma (\infty) \times [0,\infty) = \emptyset \, . \eqno (1.30.1)
$$

The ${\mathcal J}{\mathcal B}_n$'s are $2$-by-$2$ disjoined {\ibf except} for the fact that, connected to those immortal singularities which survive at the level $Y^2$ (see {\rm (\ref{eq1.23})}), this may force transversal intersection lines where two ${\mathcal B}$'s cut through each other,}
$$
L = {\mathcal B} ({\rm BLACK}) \cap {\mathcal B} ({\rm BLACK}) \, .
$$

2) {\it There is a natural inclusion
\begin{equation}
\label{eq1.31}
\Theta^3 ({\rm provisional}) \subset \Theta^3 ({\rm new}) - \sum_n {\mathcal B}_n \times [0,\infty) \, ,
\end{equation}
the $\Theta^3 ({\rm new})$ is $\Gamma$-equivariant and so are $\Theta^3 ({\rm provisional})$ and the map  {\rm (\ref{eq1.31})}.}

\medskip

3) {\it Inside $\Theta^3 ({\rm provisional})$ lives another cell-complex, staying away from those things at the infinity of $\Theta^3 ({\rm provisional})$ which prevent the action 
$\Gamma \times \Theta^3 ({\rm provisional}) \longrightarrow \Theta^3 ({\rm provisional})$ from being co-compact. Let us call this cell-complex, which is a {\rm good approximation} of $\Theta^3 ({\rm provisional})$ and which will be made explicit much later,
\begin{equation}
\label{eq1.32}
\Theta^3 (\mbox{\rm co-compact}) \subset \Theta^3 ({\rm provisional}) \, .
\end{equation}
The $\Theta^3 (\mbox{\rm co-compact})$ inherits a free action from the free action of $\Gamma$ on $\Theta^3 ({\rm provisional})$ and this $\Gamma \times \Theta^3 (\mbox{\rm co-compact}) \longrightarrow \Theta^3 (\mbox{\rm co-compact})$ is now {\ibf co-compact} (i.e. it has a compact fundamental domain).

\smallskip

Moreover, the inclusion {\rm (\ref{eq1.32})} is itself $\Gamma$-equivariant.}

\medskip

4) {\it There is a big collapse
\begin{equation}
\label{eq1.33}
\Theta^3({\rm new}) \overset{\mbox{\footnotesize{\rm collapse (1.29)}}}{-\!\!\!-\!\!\!-\!\!\!-\!\!\!-\!\!\!-\!\!\!-\!\!\!-\!\!\!-\!\!\!-\!\!\!\longrightarrow} \Theta^3 ({\rm provisional}) \longrightarrow \Theta^3 (\mbox{\rm co-compact}) \, .
\end{equation}
The big collapse {\rm (\ref{eq1.33})} is itself $\Gamma$-equivariant.}

\medskip

5) {\it The multi-game {\rm (\ref{eq1.28})} can be played in such a way that, in the context of {\rm (\ref{eq1.23})} we should find a natural isomorphism
\begin{equation}
\label{eq1.34}
{\rm Sing} \, Y^2 \approx {\rm Sing} \, \widetilde M (\Gamma) \, .
\end{equation}
This means the following. When a given immortal singularity $\overline S \subset {\rm Sing} \, \widetilde M (\Gamma)$ breaks into a double infinity of immortal singularities $S \subset {\rm Sing} \, \Theta^3 (fX^2)_{\rm (I \, or \, II)}$, then out of all these $S$'s, {\ibf one and exactly one} remains alive at the level ${\rm Sing} \, Y^2 \approx {\rm Sing} \, \Theta^3 (Y^2)$.

\smallskip

We also find that, at the level of the immortal singularities, we have}
\begin{equation}
\label{eq1.35}
{\rm Sing} \, \Theta^3 (Y^2) = {\rm Sing} \, \Theta^3 (\mbox{co-compact}) \, .
\end{equation}

\bigskip

[The ${\mathcal B} \times \{0\}$'s are NOT counted among the immortal singularities ${\rm Sing} (\ldots)$ nor are the $\ring\Sigma (\infty) \times \{0\} = \ring\Sigma (\infty)$'s.]

\smallskip

The proof of Lemma 1.5 will occupy the Section III of the present paper.

\newpage

\section{From GSC to Dehn-Exhaustibility}\label{sec2}
\setcounter{equation}{0}

In this section our concern is to show that from the fact that
\setcounter{equation}{-1}
\begin{equation}
\label{eq2.0}
S_u ({\rm new}) = \Theta^4 (\Theta^3 ({\rm new}) , {\mathcal R}) \times B^N
\end{equation}
is GSC we can deduce that $\Theta^3 ({\rm new})$ itself is {\ibf Dehn-exhaustible}, a property which is stronger than QSF, in the sense that DE $\Longrightarrow$ QSF. Dehn-exhaustibility, which will be formally defined below, comes with two variants $4^{\rm d}$ Dehn-exhaustibility and $3^{\rm d}$ Dehn-exhaustibility. The DE notion has its roots in my old papers \cite{23}, \cite{24}, \cite{25} as well as in the related, but independent work of A. Casson \cite{9}. It may well have provided the inspiration for introducing the QSF \cite{3}, \cite {35}.

\bigskip

\noindent {\bf Comments.} 

\smallskip

A) The formula (\ref{eq2.0}) is schematical. One has actually to proceed like in (\ref{eq1.1}), (\ref{eq1.2}), (\ref{eq1.3}), i.e. delete the $p_{\infty\infty} (S)$'s and add compensatory 2-handles. In none of the dimensions which are involved in (\ref{eq2.0}), three, four and large $N+4$ do we have smooth manifolds, only cell-complexes.

\medskip

B) The action of $\Gamma$ on $\Theta^3 ({\rm new})$ is not co-compact. But what we will show in the next section, is that the big collapse in (\ref{eq1.33}) is nice enough so as to make possible the implication
$$
\Theta^3 ({\rm new}) \in {\rm DE} \Longrightarrow \Theta^3 (\mbox{co-compact}) \in {\rm QSF} \, .
$$

Since $\Gamma$ has a free co-compact action on $\Theta^3 (\mbox{co-compact})$, this will imply then that $\Gamma \in {\rm QSF}$. \hfill $\Box$

\bigskip

Now, just like the $\Theta^3 (fX^2)_{\rm II}$ in (\ref{eq1.1}), our $\Theta^3 ({\rm new})$ has the following general structure
\begin{equation}
\label{eq2.1}
\Theta^3 ({\rm new}) = \underbrace{\left[ \Theta^3 ({\rm new}) (\mbox{where the $\underset{p_{\infty\infty} (S)}{\sum} p_{\infty\infty} (S) \times (-\varepsilon , \varepsilon)$ is deleted)} \right]}_{{\mbox{\scriptsize this is a cell-complex (certainly not a 3-manifold,}} \atop \mbox{\scriptsize it has immortal singularities) and we will call it $[\Theta^3]$}} \cup \sum_{p_{\infty\infty} (S)} D^2 \left[ -\frac\varepsilon4 , \frac\varepsilon4 \right] \, ,
\end{equation}
where the two pieces are joined along $\underset{p_{\infty\infty} (S)}{\sum} C(p_{\infty\infty} (S)) \times \left[ -\frac\varepsilon4 , \frac\varepsilon4 \right]$. The object replacing now the $\Theta^4 (\Theta^3 (fX^2)$, ${\mathcal R})_{\rm II}$ from (\ref{eq1.2}), and which we called loosely $\Theta^4 (\Theta^3 ({\rm new}), {\mathcal R})$ in Lemma 1.5 is a train-track smooth 4-manifold (coming with its smooth triangulation, i.e. it is again a cell-complex, but less singular), with the following general structure
\begin{equation}
\label{eq2.2}
\Theta^4 \equiv \Theta^4 (\Theta^3 ({\rm new}), {\mathcal R}) = \underbrace{\Theta^4 [\Theta^3] , {\mathcal R})}_{{\mbox{\scriptsize this is a smooth}} \atop \mbox{\scriptsize 4-manifold $Y^4$}} \quad \bigcup_{\overbrace{X^3}} \quad \underbrace{\underset{p_{\infty\infty}(S)}{\sum} D^2 \times \left[ -\frac\varepsilon4 , \frac\varepsilon4 \right] \times I}_{{\mbox{\scriptsize this is a smooth}} \atop \mbox{\scriptsize manifold called $Z^4$}} \, ,
\end{equation}
where $X^3 \equiv \underset{p_{\infty\infty}(S)}{\sum} C(p_{\infty\infty}(S)) \times \left[ -\frac\varepsilon4 , \frac\varepsilon4 \right] \times I$, and where this $X^3$ comes with smooth embeddings
$$
Y^4 \supset \ring Y^4 \longleftarrow X^3 \longrightarrow \partial Z^4 \subset Z^4 \, . \eqno (2.2.1)
$$

\bigskip

\noindent {\bf Definition 2.1.} We will say that the $\Theta^4$ above is {\ibf 4$^{\mbox{\ibf d}}$ Dehn-exhaustible}, if for every compact sub-complex $K \subset \Theta^4$ we can find a compact simply connected $4^{\rm d}$ complex $M^4$ without cells of dimension $< 4$ which are NOT faces of some $4^{\rm d}$ cell, and a commutative diagram
\begin{equation}
\label{eq2.3}
\xymatrix{
K \ar[rr]^{j} \ar[dr]_-i &&M^4 \ar[dl]^-{g}  \\ 
&\Theta^4
}
\end{equation}
where $i$ is the canonical inclusion, $j$ is a simplicial injection, $g$ a simplicial  {\ibf immersion}, and where the following Dehn-type condition is fulfilled $jK \cap M_2 (g) = \emptyset$, inside $M^4$.

\smallskip

[We could as well ask that $M^4$ itself should be a $4^{\rm d}$ smooth train-track manifold, coming with a smooth $g$, but this will be unnecessary.]

\bigskip

\noindent {\bf Lemma 2.2.} {\it The $\Theta^4$ from {\rm (\ref{eq2.2})}  is $4^{\rm d}$ Dehn exhaustible.}

\bigskip

\noindent {\bf Proof.} We will adapt to the present situation the method of proof from \cite{23}, the hypothesis $V^3 \times B^n \in {\rm GSC}$ of the theorem in \cite{23} being replaced by the fact that $S_u ({\rm new}) \in {\rm GSC}$. Of course, in \cite{23} the $V^3$ which was then $3$-dimensional and smooth was an open manifold while now the non-compact non-smooth $4^{\rm d}$ $\Theta^4$ has $\partial \, \Theta^4 \ne \emptyset$. We consider, like in (\ref{eq1.3}) the projection
\begin{equation}
\label{eq2.4}
S_u ({\rm new}) \underset{\pi \, \equiv \, \pi_{N+4,4}}{-\!\!\!-\!\!\!-\!\!\!-\!\!\!-\!\!\!-\!\!\!-\!\!\!-\!\!\!\longrightarrow} \Theta^4
\end{equation}
and the zero-section
\begin{equation}
\label{eq2.5}
\Theta^4 \overset{\mathcal J}{-\!\!\!-\!\!\!-\!\!\!\longrightarrow} S_u ({\rm new}) = Y^4 \times B^N \underset{\overbrace{\mbox{\footnotesize$X^3 \times \frac12 B^N$}}}{\cup} Z^4 \times B^N \, ,
\end{equation}
gotten by sending the $Y^4 , X^3 , Z^4$ diffeomorphically into the respective $Y^4 \times \{0\} , X^3 \times \{0\} , Z^4 \times \{0\}$, when $0 = \bigl\{$the common center of $B^N$ and $\frac12 \, B^N \bigl\}$. Each of the two $Y^4$ and $Z^4$, which are smooth, will be endowed with a riemannian metric such that

\bigskip

\noindent (2.5.1) \quad On the $X^3$ which is contained both in ${\rm int} \, Y^4$ and in $\partial Z^4$, so that $X^3 = Y^4 \cap Z^4$ ($\Theta^4$ is train-track, remember), the two metrics coincide.

\bigskip

\noindent (2.5.2) \quad Each of the two $Y^4$ and $Z^4$ can be covered by small, geodesically convex charts, generically denoted by $U_i$.

\bigskip

\noindent (2.5.3) \quad In terms of both the metrics coming from $Y$ or from $Z$, the $X^3$ above is locally geodesically convex.

\bigskip

Since $\partial Y^4 \ne \emptyset \ne \partial Z^4$, the condition (2.5.2) is certainly not automorphic. One better starts with metrics on $\partial Y^4 , \partial Z^4$ and then one extends these carefully in the neighbourhood of the boundaries towards the interior, taking care of (2.5.3) among other things. Once we are far from the boundary, the extension becomes very easy.

\smallskip

Next, the $B^N$ itself is endowed with a standard euclidean metric. This will yield an atlas ${\mathcal U} = \bigl\{ U_i \times B^N \ {\rm OR} \ U_i \times \frac12 \, B^N$, according to the case$\bigl\}$, for $S_u ({\rm new})$.

\smallskip

Notice that $\pi \mid \pi^{-1} (\partial Y^4 + \partial Z^4 (\supset X^3))$ is violently degenerate, and a priori this is not compatible with the technology of \cite{23}. Our first step in the proof will be to change the geometry of (\ref{eq2.4}), without touching to the zero-section (\ref{eq2.5}), so as to demolish this  unwanted degeneracy. Figure 2.1 suggests how to achieve this goal for $\partial Y^4 + (\partial Z^4 - X^3)$.

$$
\includegraphics[width=13cm]{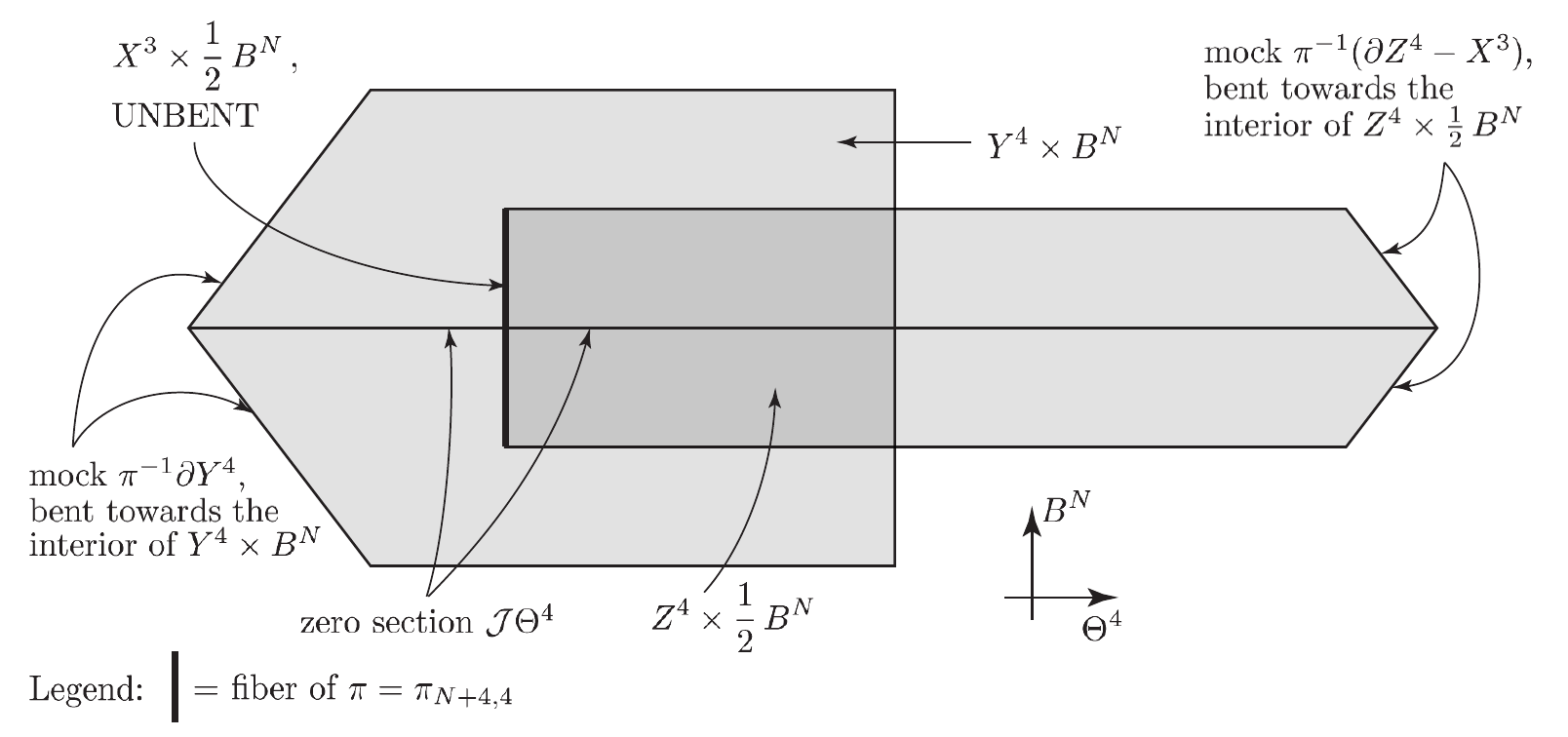}
$$
\label{fig2.1}

\centerline {\bf Figure 2.1.} 

\smallskip

\begin{quote}
A modification of the geometry of $S_u ({\rm new})$ in the neighbourhood of $\pi^{-1} (\partial Y^4 + (\partial Z^4 - X^3)$). When the shading is double, the two objects which (outside of $X^3 \times \frac12 \, B^N$) are disjoined, appear superposed by the obvious projection on $Y^4 \times B^N$; this is just an optical illusion.
\end{quote}

\bigskip

What happen to $X^3$ will be discussed afterwards. The trick here, is to change $\pi^{-1} (\partial Y^4 + (\partial Z^4 - X^3)) \subset \partial S_u ({\rm new})$, {\ibf without touching} to the zero-section (\ref{eq2.5}), into an object called
$$
{\rm mock} \, (\pi^{-1} (\partial Y^4 + (\partial Z^4 - X^3))) \, ,
$$
which is such that $\pi \mid {\rm mock} \, (\pi^{-1} (\partial Y^4 + (\partial Z^4 - X^3)))$ becomes {\ibf non-degenerate}. One bends $\pi^{-1} (\partial Y^4 + (\partial Z^4 - X^3))$ symmetrically around the zero-section, towards ${\rm int} \, S_u ({\rm new})$. We cannot apply this same treatment to $\pi^{-1} X^3 \subset \pi^{-1} (\partial Z^4)$, since we cannot mock around with the projection
\begin{equation}
\label{eq2.6}
X^3 \times \frac12 \, B^N \equiv \pi^{-1} X^3 \overset{\pi \mid \pi^{-1} X^3}{-\!\!\!-\!\!\!-\!\!\!-\!\!\!-\!\!\!-\!\!\!-\!\!\!-\!\!\!-\!\!\!\longrightarrow} X^3 
\end{equation}
which is common to the two pieces $Y,Z$, making that $\pi \mid X^3$ has to stay {\ibf degenerate}. How we will manage to live with this, will be soon explained.

\smallskip

Anyway, the modification of $S_u ({\rm new})$ just described, will not concern the (\ref{eq2.6}). As Figure 2.1 may suggest, we have an inclusion
$$
\{{\rm modified} \ S_u({\rm new})\} \subsetneqq \{{\rm original} \ S_u({\rm new})\} \, ,
$$
and the newly modified $S_u({\rm new})$ is endowed with the atlas ${\mathcal V} = {\mathcal U} \mid \{$modified $S_u({\rm new})\}$, coming with geodesically convex charts. From now on, $S_u({\rm new})$ should mean the $\{$modified $S_u({\rm new})\}$.

\bigskip

\noindent {\bf Sublemma 2.2.A.} {\it There is a smooth triangulation  $\tau$ of $S_u({\rm new})$ (meant now in its modified form), s.t.}

\medskip

1) {\it $\tau$ is GSC (actually asking only for WGSC would suffice for our present purposes) and $\pi^{-1} X^3$ is a subcomplex.}

\medskip

2) {\it The zero-section ${\mathcal J} \Theta^4 \subset S_u({\rm new})$ is a  subcomplex too.}

\medskip

3) {\it The simplexes of $\tau$ are ${\mathcal V}$-small, where ${\mathcal V}$ is the atlas above. We assume that the simplices
$$
\sigma \subset V \in {\mathcal V}
$$
of $\tau$ are geodesically convex. When the vertices of $\sigma$ are slightly perturbed, generically, into different positions inside ${\mathcal V}$, then this defines another, still geodesically convex, version of $\sigma$.}

\medskip

4) {\it Before we can actually state this new item, some preliminaries are necessary. We denote
\begin{equation}
\label{eq2.7}
\tau^{(4)} \diagup X^3 \equiv \{\mbox{the $4$-skeleton $\tau^{(4)}$ of $\tau$, from which all the ${\rm int} \, \sigma^4$ where $\sigma^4$ is a $4$-simplex}
\end{equation}
$$
\mbox{contained in $\pi^{-1} X^3$ are {\ibf deleted}}\} \, .
$$

This $\tau^{(4)} \diagup X^3$ continues to be GSC, just like $\tau$. With this, one can start by perturbing, in the manner explained at {\rm 3)} above, first the restriction $\tau^3 \cap \pi^{-1} X^3$, keeping each $3$-simplex inside its $\{\mbox{convex ${\mathcal V}$-chart}\} \cap \pi^{-1} X^3$, and next continue to perturb, in agreement with this, the rest of $\tau^{(4)} \diagup X^3$, keeping again each  $4$-simplex inside its convex ${\mathcal V}$-chart.

\smallskip

What such a perturbation can achieve are the following items: $\pi \mid \sigma^3$, where $\sigma^3 \subset \pi^{-1} X^3$ and each $\pi \mid \sigma^4$, where $\sigma^4 \subset \tau^{(4)} \diagup X^3$, should be an isomorphism on its image, and for these various $\sigma^3$'s, $\sigma^4$'s we should also have that
\begin{equation}
\label{eq2.8}
\mbox{$\pi \, \sigma_1^3 \cap \pi \, \sigma_2^3$ are in general position, modulo their incidences $\sigma_1^3 \cap \sigma_2^3$}
\end{equation}
$$
\mbox{relations, which should be respected. Similarly for $\pi \, \sigma_1^4 \cap \pi \, \sigma_2^4$} \, .
$$
}

\medskip

5) {\it Next, there is a {\ibf good subdivision} (and all one has to know right now about such subdivision is that they are plenty of them and that they preserve things like GSC and/or WGSC), call it $\tau \to \theta$, for which we denote $\theta^{(4)} \diagup X^3 = \{$the $\theta$-subdivision of $\tau^{(4)} \diagup X^3\}$, which is GSC, such that for $\lambda = 3$ or $4$ (see the context of $4)$ above) the intersection $\pi \, \sigma_1^{\lambda} \cap \pi \, \sigma_2^{\lambda}$ becomes a subcomplex of both $\pi \, \sigma_1^{\lambda}$ and of $\pi \, \sigma_2^{\lambda}$. Moreover, the following map
\begin{equation}
\label{eq2.9}
\theta^{(4)} \diagup X^3 \overset{\pi \mid (\theta^{(4)} \diagup X^3)}{-\!\!\!-\!\!\!-\!\!\!-\!\!\!-\!\!\!-\!\!\!-\!\!\!-\!\!\!-\!\!\!\longrightarrow} \Theta^{(4)} \, ,
\end{equation}
is {\ibf both simplicial and non-degenerate}.}

\bigskip

We will come back a bit later to the good subdivisions, which allow us to get from $S_u({\rm new}) \in {\rm GSC}$ to $\theta^{(4)} \diagup X^3 \in {\rm GSC}$, but right now I will offer some comments, in lieu of a formal proof for the Lemma 2.2. Notice, to begin with, that if our $\Theta^4$ would be replaced by a smooth open $4$-manifold $V^4$ s.t. $V^4 \times B^N \in {\rm GSC}$ then, the analogue of our lemma could be proved by a very simple-minded transposition of the arguments which have been used in \cite{23}.

\smallskip

As things actually stand, our $\Theta^4$ is only a train-track manifold, with the kind of structure which (\ref{eq2.2}) prescribes and one also has $\partial \, \Theta^4 \ne \emptyset$. The modification from Figure 2.1 suggests how one deals with $\partial \, \Theta^4 - X^3$, while the train-track locus $X^3$ is handled like in the point 4) of the Sublemma 2.2.A; see here, in particular, the formula (\ref{eq2.7}). In this context, I will offer here the following pedagogical toy-model.

\smallskip

Consider the linear spaces $A^{n-1} \subset B^n , C^N$ and the obvious projection $B^n \times C^N \overset{\pi}{-\!\!\!-\!\!\!\longrightarrow} B^n$. In this context, we consider simplexes $\sigma^n \subset A^{n-1} \times C^N$ and $\sigma^{n+1} \subset B^n \times C^N$ s.t.
$$
\sigma^{n+1} \cap (A^{n-1} \times C^N) = \sigma^n \subset \partial \, \sigma^{n+1} \, .
$$

In this generic set-up (where we should think in terms of $A^{n-1} \cong X^3 \subset Y^4 \cong B^N$ and $C^N \cong B^N$), after small admissible perturbations, both of the simplicial maps
$$
\partial \, \sigma^n  \overset{\pi}{-\!\!\!-\!\!\!\longrightarrow} A^{n-1} \, , \qquad \partial \, \sigma^{n+1} - {\rm int} \, \sigma^n \overset{\pi}{-\!\!\!-\!\!\!\longrightarrow} B^n \, ,
$$
can be rendered non-degenerate. End of the toy-model.

\bigskip

The point of this whole discussion is that with the items described, we get to (\ref{eq2.9}) and one can apply now, more or less directly, the arguments from the Section 4 of \cite{23}, and from our Lemma 2.2. \hfill $\Box$

\bigskip

But next, we will move from the relatively smooth $\Theta^4$ in (\ref{eq2.2}) to the much more singular $\Theta^3 ({\rm new})$ (from (\ref{eq2.2})) and prove the implication
$$
\{\Theta^4 \ \mbox{is Dehn-exhaustible}\} \Longrightarrow \{\Theta^3 ({\rm new}) \ \mbox{is Dehn-exhaustible}\} \, ,
$$
object of the next Lemma 2.4. The general idea is to adapt, once more, the technology from \cite{23}, but the road is now steeper than for the Lemma 2.2. Also, the initial input is now no longer the GSC property of $S_u ({\rm new})$, but the $4^{\rm d}$ Dehn-exhautibility of $\Theta^4$.

\smallskip

Before really proceeding further, I will open a LONGUISH PRENTICE concerning the good subdivisions which were mentioned in the statement of the Sublemma 2.2.A.

\smallskip

Whenever we talk about subdivisions for a simplicial complex, we will always mean linear subdivisions. Among these are the barycentric and stellar subdivisions, which clearly preserve the GSC feature, while the general linear ones might not. Concerning the stellar subdivisions there are also the old tricky results of Alexander and Newman.

\smallskip

For all these matters, there is a very nice and efficient approach due to Larry Siebenmann and, since his work is not available in print, at least not right now, I will briefly outline it here.

\smallskip

Siebenmann starts by introducing {\it cellulations}, which are an extension of simplicial complexes: instead of using simplexes we use now compact cells $D$ with a linear-convex structure. The notion of (linear) subdivision extends in an obvious way to cellulations and, also, instead of subcomplexes we can introduce now sub-cellulations. What we have gained with this approach is, among other things, the following useful fact: if $Y \subset Z$ is a sub-cellulation, then any subdivision $Y'$ extends canonically to a subdivision of $Z$, not affecting the open cells in $Z-Y$. An important class of subdivisions are the BISSECTIONS. These are localized at the level of an $i$-cell $D^i$ and are obtained by cutting $D^i$ with a hyperplane $H^{i-1} \subset D^i$ and splitting $D^i$ itself and any sub-cell of $D^i$ met by $H^{i-1}$ in the obvious way. Our ``useful fact'' above extends to bissections. No genericity conditions are required here fo $H^{i-1} \subset D^i$. If $X$ is a cellulation, then there is also a canonical way to subdivide $X$ to a simplicial complex $X$ (simplicial). We start by picking up for each $2$-cell $D^2 \subset X$ a point $q \in \ring D^2$ and then we subdivide $D^2$ in a way which should be obvious. Then we do the same for all the $3$-cells, next for the $4$-cells, a.s.o. The operation $X \Rightarrow X$ (simplical) will be called {\ibf stellation}, and quite obviously bissection and stellations preserve the GSC property.

\smallskip

Finally, there is also the following very useful fact, which is easy to prove, in Siebenmann's context. {\it If $X$ is a cellulation and $X'$ a (linear) subdivision of $X$, then there is a third cellulation $X_1$ such that one can go both from $X$ and from $X'$ to $X_1$ via bissections.} This is a nice elegant substitute for those old theorems of Alexander and Newman, the proofs of which was always a clumsy affair.

\smallskip

With this we close our prentice and our goal subdivisions {\ibf are} the bissections and stellations above.

\bigskip

What follows next is Definition 2.1 adapted now for $\Theta^3 ({\rm new})$.

\bigskip

\noindent {\bf Definition 2.3.} We define now, on the same lines as in the Definition 2.1 above, the {\ibf 3$^{\mbox{\ibf d}}$ Dehn-exhaustibility}. This definition makes sense for any $3^{\rm d}$ cell complex, in particular for $\Theta^3 ({\rm new})$. This is the only case where it will be needed, and we state it only for it. We will say that $\Theta^3 ({\rm new})$ is {\ibf 3$^{\mbox{\ibf d}}$ Dehn-exhaustible}, if for any compact subcomplex $k \subset \Theta^3 ({\rm new})$ we can find a compact simplicial complex $K^3$ with $\pi_1 K^3 = 0$, coming with a commutative diagram
\begin{equation}
\label{eq2.10}
\xymatrix{
k \ar[rr]^{j} \ar[dr]_-i &&K^3 \ar[dl]^-{\chi}  \\ 
&\Theta^3 ({\rm new})
}
\end{equation}
where $i$ is the canonical inclusion, $j$ a simplicial injection, $\chi$ a simplicial {\ibf immersion}, and where $jk \ \cap \ M_2 (\chi) = \emptyset$, inside $K^3$. 

\hfill $\Box$

\bigskip

\noindent {\bf Lemma 2.4.} {\it $\Theta^3 ({\rm new})$ is $3^{\rm d}$ Dehn-exhaustible.}

\bigskip

\noindent {\bf Proof.} For the convenience of the reader, we re-write schematically the formulae (\ref{eq2.1}), (\ref{eq2.2})
$$
\Theta^3 ({\rm new}) = [\Theta^3] \underset{\overbrace{\mbox{\footnotesize$X^2 \equiv C (p_{\infty\infty} (S)) \times \left[ -\frac\varepsilon4 , \frac\varepsilon4 \right]$}}}{\cup} D^2 \times \left[ -\frac\varepsilon4 , \frac\varepsilon4 \right] ,
$$
and
$$
\Theta^4 = \Theta^4 ([\Theta^3] , {\mathcal R}) \underset{\overbrace{\mbox{\footnotesize$X^2 \times I \equiv X^3$}}}{\cup} D^2 \times \left[ -\frac\varepsilon4 , \frac\varepsilon4 \right] \times I = Y^4 \underset{X^3}{\cup} Z^4 \, ,
$$
coming with $L_1^3 \equiv \partial Y^4$, $L_2^3 \equiv \partial Z^4$. We also have the (2.2.1). Like in (\ref{eq2.4}) $+$ (\ref{eq2.5}) we have again a natural projection
$$
\Theta^4 \overset{\pi_{4,3}}{-\!\!\!-\!\!\!-\!\!\!-\!\!\!\longrightarrow} \Theta^3 ({\rm new})
$$
and a natural zero-section
$$
\Theta^3 ({\rm new}) \overset{\mathcal J}{-\!\!\!-\!\!\!-\!\!\!-\!\!\!\longrightarrow} \Theta^4 \, ;
$$
we will not always distinguished, notationally, between $\Theta^3 ({\rm new})$ and ${\mathcal J} \Theta^3 ({\rm new})$.

\bigskip

\noindent {\bf Subemma 2.4.A.} 1) {\it Without any loss of generality, the following map, where now $\pi \equiv \pi_{4,3}$,
\begin{equation}
\label{eq2.11}
L_1^3 \overset{\pi \mid L_1^3}{-\!\!\!-\!\!\!-\!\!\!-\!\!\!\longrightarrow} [\Theta^3]
\end{equation}
is a submersion, except for simple fold singularities.}

\medskip

2) {\it Moreover, we have an isomorphism
\begin{equation}
\label{eq2.12}
(Y^4 , L_1^3) = \left( [\Theta^3] \cup (L_1^3 \times [0,1]) , L_1^3 \times \{1\} \right) \, ,
\end{equation}
where $L_1^3 \times [0,1]$ gets glued to $[\Theta^3]$ along $\pi \mid L_1^3 \times \{0\}$, and where for $t > 0$, each $\pi \mid L_1^3 \times t$ is isomorphic to $\pi \mid L_1^3$ in {\rm (\ref{eq2.8})}. So we have a foliation ${\mathcal F}_1$ with $3^{\rm d}$ leaves $L_1^3 \times t$, $t > 0$ of $Y^4 - {\mathcal J} [\Theta^3]$. When it comes to $X^3 = X^2 \times I \subset {\rm int} \, Y^4$, then ${\mathcal F}_1 \mid X^3$ is just the restriction of the standard foliation of $X^3$ by the $X^2 \times t$'s, where $t\in I$. We call this foliation ${\mathcal F}_3$. With the ${\mathcal F}_1 \mid X^3$, the ${\mathcal F}_3$ extends over the zero-section ${\mathcal J} X^2 \subset X^3$.}

\medskip

3) {\it Without loss of generality, the map
\begin{equation}
\label{eq2.13}
L_2^3 \overset{\pi \mid L_2^3}{-\!\!\!-\!\!\!-\!\!\!-\!\!\!\longrightarrow} D^2 \times \left[ -\frac\varepsilon4 , \frac\varepsilon4 \right] 
\end{equation}
is such that:}
\begin{enumerate}
\item[3.1)] {\it On the piece $X^3 \subset L_2^3$ it coincides with the canonical projection $X^3 = X^2 \times I \to X^2$}
\item[3.2)] {\it On $\overline{L_2^3 - X^3}$ it is, like the map {\rm (\ref{eq2.8})}, a submersion, except for simple fold singularities.}
\end{enumerate}

\medskip

4) {\it Moreover, we have an isomorphism
\newpage
\begin{equation}
\label{eq2.14}
(Z^4 ; L_2^3) = \Biggl( \left( D^2 \times \left[ -\frac\varepsilon4 , \frac\varepsilon4 \right] \right) \cup \left( \left(\overline{L_2^3 - X^3}\right) \times [0,1] \right) \, , \ \mbox{where the two pieces are glued along}
\end{equation}
$$
\pi \left(\overline{L_2^3 - X^3}\right) \times \{0\} = D^2 \times \left[ -\frac\varepsilon4 , \frac\varepsilon4 \right] ; \ \mbox{we have here} \ L_2^3 = \left(\left(\overline{L_2^3 - X^3}\right) \times \{1\} \right) \cup X^3 (=X^2 \times [0,1]) \, ,
$$
$$
\mbox{where the two pieces are glued along} \ X^2 \times \{0,1\} = \partial \left(\overline{L_2^3 - X^3}\right) \Biggl) .
$$

In the formula above, each $\left(\overline{L_2^3 - X^3}\right) \times (t > 0)$ is isomorphic to $\overline{L_2^3 - X^3}$, defining a foliation ${\mathcal F}_2$ of $D^2 \times \left[ -\frac\varepsilon4 , \frac\varepsilon4 \right] \times I - {\mathcal J} \left(D^2 \times \left[ -\frac\varepsilon4 , \frac\varepsilon4 \right] \right)$. The trace of this foliation on $X^3$ is, outside of the zero-section where ${\mathcal F}_2$ is undefined, the ${\mathcal F}_3$ from point {\rm 2)} above.}

\bigskip

\noindent {\bf Sketch of proof.} The $Z^4$-part of the lemma should be obvious. When it comes to $Y^4$, in particular to the (\ref{eq2.8}) from 1), one should notice that the only places where $[\Theta^3]$ fails to be a $3$-manifold are either branching points or undrawable singularities.

$$
\includegraphics[width=13cm]{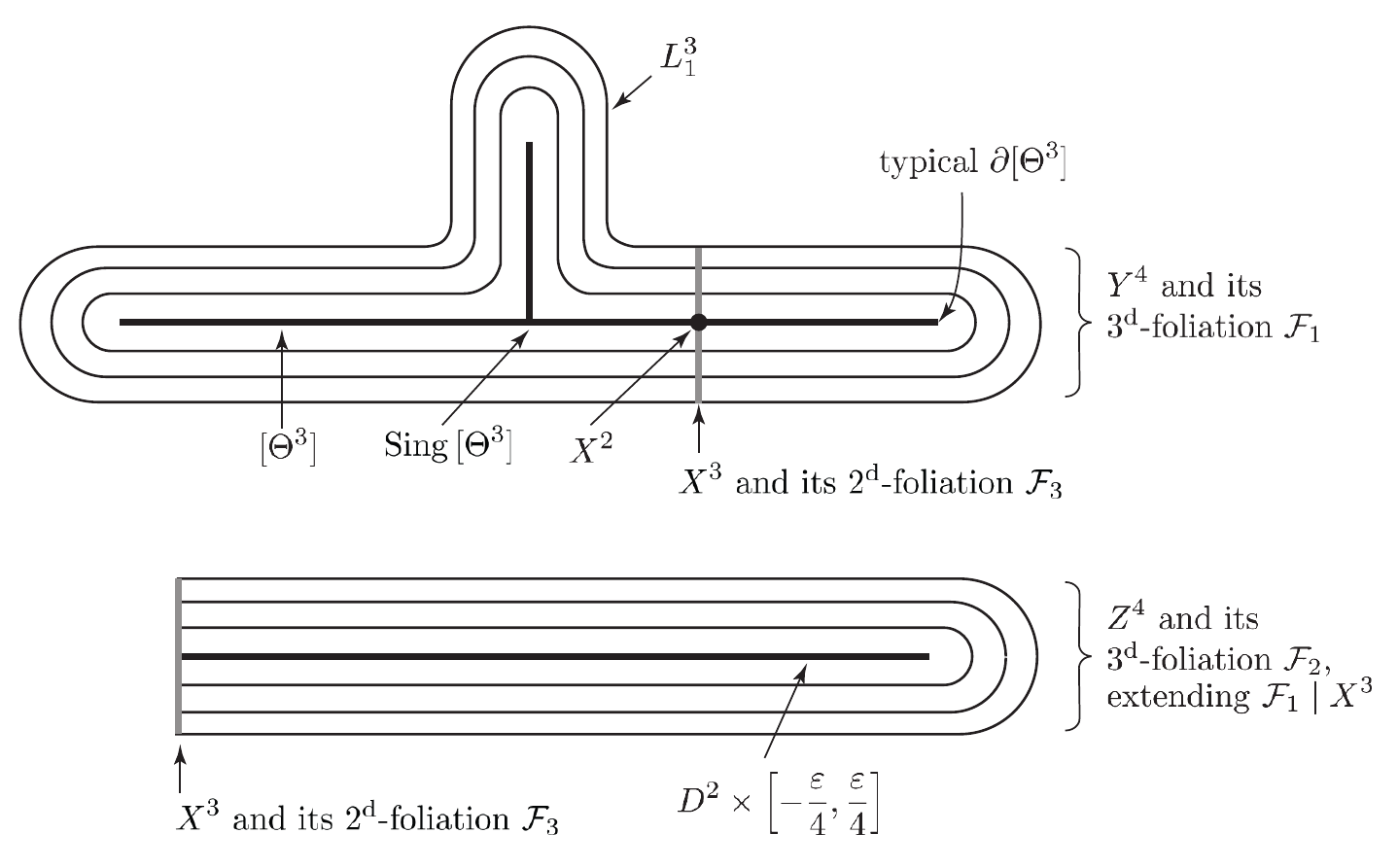}
$$
\label{fig2.21}

\centerline {\bf Figure 2.2.} 

\smallskip

\begin{quote}
Illustration for the formula ((\ref{eq2.12}) for $Y^4$) and ((\ref{eq2.14}) for $Z^4$). The two pieces we see here are glued together along $(X^3 , {\mathcal F}_3)$, so as to generate the full $(\Theta^4 , {\mathcal F})$.
\end{quote}

\bigskip

And, of course there are also the boundary points to be taken care of. For things like branchings and/or undrawable singularities, in proving the sublemma there is first a local issue to be faced and then a second issue of glueing together the local data. Details are left to the reader. \hfill $\Box$

\bigskip

When we glue together $(Y^4 , {\mathcal F}_1)$ and $(Z^4 , {\mathcal F}_2)$ along $(X^3 , {\mathcal F}_3)$, we get $\Theta^4$, endowed with a foliation, defined ouside of the zero-section. Call it ${\mathcal F}$; see here the Figure 2.2.

\bigskip

\noindent {\bf The PL Sublemma 2.4.B.} {\it There exists a smooth triangulation $\tau$ of $\Theta^4$ (which one should not mix up with the $\tau$ from the Sublemma {\rm 2.2.A}), such that:}
\begin{enumerate}
\item[1)] {\it ${\mathcal J} \Theta^3 ({\rm new}) \subset \Theta^4$ is a subcomplex of $\tau$.}
\item[2)] {\it $X^3 \subset \Theta^4$ is also a subcomplex of $\tau$ and we will introduce the notation (which should not be mixed up with {\rm (\ref{eq2.7})})
\begin{equation}
\label{eq2.15}
\tau^{(3)} \diagup X^3 \equiv \bigl\{\mbox{the $3$-skeleton $\tau^{(3)}$ of $\tau$, from which all the open cells ${\rm int} \, \sigma^3$,}
\end{equation}
$$
\mbox{where $\sigma^3 \subset X^3$ have been deleted}\bigl\}.
$$
From now on, $\pi \equiv \pi_{4,3} \mid (\tau^{(3)} \diagup X^3)$.}
\item[3)] {\it The map
\begin{equation}
\label{eq2.16}
\tau^{(3)} \diagup X^3 \overset{\pi}{-\!\!\!-\!\!\!\longrightarrow} \Theta^3 ({\rm new})
\end{equation}
is {\ibf simplicial and nondegenerate}.}
\end{enumerate}

\bigskip

\noindent {\bf Proof.} We start with a smooth triangulation $\tau_1$ of $\Theta^4$ having already the features 1) and 2) and which, moreover, is such that any $3$-simplex $\sigma^3$ not already  in ${\mathcal J} \Theta^3 ({\rm new}) \cup X^3$ is very close and almost parallel to, some leaf of ${\mathcal F}$. Similarly we ask for the 2-simplices $\sigma^2 \subset X^3$ to be parallel and very close to some leaf of ${\mathcal F} \mid X^3 \cong {\mathcal F}_3$. With this, when we are far from the fold singularities of (\ref{eq2.8}) $+$ (\ref{eq2.10}) and their counterparts on the leaves, then we may assume the $\pi \mid \sigma^3$, $\pi \mid \sigma^2$ injective already. From there on, our result is achieved by appropriate successive subdivision. Details are left to the reader. \hfill $\Box$

\bigskip

Like in (\ref{eq2.7}) we consider now $k \subset  \Theta^3 ({\rm new})$, which we assume subcomplex of the triangulation $\tau \mid  \Theta^3 ({\rm new})$. The $\pi^{-1} k \subset \Theta^4$ is ``$\pi$-closed'', meaning that $\pi^{-1} (\pi (\pi^{-1} k)) = \pi^{-1} k$; here $\pi$ is like in (\ref{eq2.13}) and clearly also $k \subset \pi^{-1} k$.

\smallskip

For each subcomplex $X \subset \tau^{(3)} \diagup X^3$ we have our basic equivalence relations (see \cite{22}, \cite{29})
\begin{equation}
\label{eq2.17}
\Psi (\pi \mid X) \subset \Phi (\pi \mid X) \, .
\end{equation}

\bigskip

\noindent {\bf Claim (2.18).} There exists a finite subcomplex $K_1$ with the feature $\tau^{(3)} \diagup X^3 \supset K_1 \supset \pi^{-1} k$, such that $K_1$ is $\pi$-closed $(\pi^{-1} \pi K_1 = K_1)$ and also that
$$
\Psi (\pi \mid K_1) \mid \pi^{-1} k = \Phi (\pi \mid \pi^{-1} k) \, .
$$

\bigskip

\noindent {\bf Proof.} One has to start by proving that
\setcounter{equation}{18}
\begin{equation}
\label{eq2.19}
\Psi (\pi) = \Phi (\pi) \ \mbox{for the map from (\ref{eq2.13})} \, ,
\end{equation}
which is done by the same arguments as for the formula (4.3) in the paper \cite{23}. From here on, the proof of our claim uses the same kind of compactness arguments as in the proof of Proposition B in \cite{23}. \hfill $\Box$

\bigskip

Since we know already that $\Theta^4$ is $4^{\rm d}$ Dehn-exhaustible we have an $M^4$, compact and simply-connected, which we may assume to be a subcomplex of $\tau$, s.t. like in the Definition 2.1,
\begin{equation}
\label{eq2.20}
\xymatrix{
K_1 \ar[rr]^{j} \ar[dr]_-i &&M^4 \ar[dl]^-{g}  \\ 
&\Theta^4
} , \ \mbox{with $j K_1 \cap M_2 (g) = \emptyset$.}
\end{equation}
Let $M^{(3)} \equiv \{$the $3$-skeleton of $M^4\} \supset M^{(3)} \diagup X^3 = \{$the obvious subcomplex of $\tau^{(3)} \diagup X^3\}$; since $\pi_1 M^4 = 0$, we also have $\pi_1 (M^{(3)} \diagup X^3) = 0$.

\smallskip

The following map $g$, restriction of the one from (\ref{eq2.20}), occurring below
$$
M^4 \supset M^{(3)} \diagup X^3 \underset{g}{-\!\!\!-\!\!\!\longrightarrow} \tau^{(3)} \diagup X^3 \subset \Theta^4 \, , \eqno (2.20.1)
$$
\vglue -6mm
$$
{\mid\mbox{\hglue -2mm}}_{-\!\!-\!\!-\!\!-\!\!-\!\!-\!\!-\!\!-\!\!-\!\!-\!\!-\!\!-\!\!-\!\!-\!\!-\!\!-\!\!-\!\!-\!\!-\!\!-\!\!-\!\!-\!\!-\!\!-\!\!-\!\!-\!\!-\!\!-\!\!-\!\!-\!\!-\!\!-\!\!-\!\!-\!\!-\!\!-\!\!-\!\!-\!\!-\!\!-\!\!-}{\mbox{\hglue -2mm}\uparrow}
$$

\smallskip

\noindent is simplicial nondegenerate; actually it is an immersion, just like the $g$ in (\ref{eq2.20}). With this, we extract the following commutative diagram from (\ref{eq2.20}) + (2.20.1)
\begin{equation}
\label{eq2.21}
\xymatrix{
K_1 \ar[rr]^{j} \ar[dr]_-i &&M^{(3)} \diagup X^3 \ar[dl]^-{g}  \\ 
&\Theta^4 \ar[rr]^{\pi} &&\Theta^3 ({\rm new}) \, ,
} 
\end{equation}
where $g$ is a {\ibf simplicial immersion} and where the following Dehn-type property gets inherited from (\ref{eq2.20}),
\begin{equation}
\label{eq2.22}
M^{(3)} \diagup X^3 \supset M_2 (g) \cap jK_1 = \emptyset \, .
\end{equation}
From (\ref{eq2.21}) we can pull out the composite map
$$
M^{(3)} \diagup X^3 \underset{g}{-\!\!\!-\!\!\!\longrightarrow} \Theta^4 \overset{\pi}{-\!\!\!-\!\!\!\longrightarrow} \Theta^3 ({\rm new}) \, .
$$
Here the map $g$ factors through $\tau^{(3)} \diagup X^3$, like in the (2.20.1) above, and the whole composite map $\pi \circ g$ is both simplicial and nondegenerate; see here the PL Sublemma 2.4.B too.

\smallskip

At this point, just like in \cite{23} we get an induced immersion
\begin{equation}
\label{eq2.23}
(M^{(3)} \diagup X^3) \diagup \Psi (\pi \circ g) \overset{g_1}{-\!\!\!-\!\!\!\longrightarrow} \Theta^3 ({\rm new})
\end{equation}
which comes with
$$
\pi_1 ((M^{(3)} \diagup X^3) \diagup \Psi (\pi \circ g)) = 0 \, . \eqno (2.23.1)
$$

\bigskip

\noindent {\bf Claim (2.24).} From the inclusion $k \subset \pi^{-1} k \subset K_1 \subset M^{(3)} \diagup X^3$, we can get a second inclusion
$$
k \subset (M^{(3)} \diagup X^3) \diagup \Psi (\pi \circ g) \, , \ \mbox{coming with} \ M_2 (g_1) \cap k = \emptyset \, .
$$

\bigskip

\noindent {\bf Proof.} We start by noticing that the situation marked $(*)$ below {\ibf cannot occur}
$$
x \in g \left( M^{(3)} \diagup X^3 - K_1\right) \subset \Theta^4 \, , \ y \in K_1 \subset \Theta^4 \ {\rm and} \ z \equiv \pi x = \pi y \in \Theta^3 ({\rm new}) \, . \eqno (*)
$$
Here is why $(*)$ cannot happen. Assume it does and denote $z \equiv \pi x = \pi y$. We then automatically get that $z = K_1$ and, since $K_1$ is $\pi$-closed, we also get that $\pi^{-1} z = \pi^{-1} \, \pi y \subset K_1$. Our $(*)$ above means that $x \in \pi^{-1} z \subset K_1$, which contradicts the Dehn property (\ref{eq2.22}).

\smallskip

So, by now we have proved that, at level $\Theta^3 ({\rm new})$, we have
\setcounter{equation}{24}
\begin{equation}
\label{eq2.25}
\pi \circ g \left(M^{(3)} \diagup X^3 - K_1 \right) \cap \pi K_1 = \emptyset \, .
\end{equation}
This (\ref{eq2.25}) implies that, when we restrict the equivalence relation $\Psi (\pi \circ g)$ which occurs in (\ref{eq2.23}), from $M^{(3)} \diagup X^3$ to the smaller set $K_1$, then this operates all the identifications $\Psi (\pi \mid K_1)$, but {\ibf nothing more}; hence $K_1 \diagup \Psi (\pi \mid K_1) \subset (M^{(3)} \mid X^3) \diagup \Psi (\pi \circ g)$.

\smallskip

We have now inclusions
\begin{equation}
\label{eq2.26}
\pi^{-1} k \diagup \Psi (\pi \mid K_1) \subset K_1 \diagup \Psi (\pi \mid K_1) \subset \left(M^{(3)} \diagup X^3 \right) \diagup \Psi ( \pi \circ g ) \, .
\end{equation}
Our Claim (2.18) tells us that, for the $\pi$-closed set $k$, we have
\begin{equation}
\label{eq2.27}
\pi^{-1} k \diagup \Psi (\pi \mid K_1) = \pi^{-1} k \diagup \Phi (\pi \mid \pi^{-1} k) = k \, .
\end{equation}

The combination of (\ref{eq2.26}) and (\ref{eq2.27})  gives us the desired inclusion occurring in our Claim (2.24), namely the
$$
k \subset \left(M^{(3)} \diagup X^3 \right) \diagup \Psi ( \pi \circ g ) \, ,
$$
and according to (\ref{eq2.25}) this inclusion factors through the $K_1 \diagup \Psi (\pi \mid K_1)$. In the Claim (2.24) there is also a Dehn-part, to the proof of which we turn now. For this purpose, in the context of
$$
k \subset \Theta^3 ({\rm new}) \overset{\mathcal J}{-\!\!\!-\!\!\!\longrightarrow} \Theta^4 \overset{\pi}{-\!\!\!-\!\!\!\longrightarrow} \Theta^3 ({\rm new}) \, , 
$$
\vglue -6mm
$$
{\mid\mbox{\hglue -2mm}}_{-\!\!-\!\!-\!\!-\!\!-\!\!-\!\!-\!\!-\!\!-\!\!-\!\!-\!\!-\!\!-\!\!-\!\!-\!\!-\!\!-\!\!-\!\!-\!\!-\!\!-\!\!-\!\!-\!\!-\!\!-\!\!-\!\!-\!\!-\!\!-\!\!-\!\!-\!\!-\!\!-\!\!-\!\!-\!\!-\!\!-\!\!-\!\!-\!\!-\!\!-}{\mbox{\hglue -2mm}\uparrow}
$$
\vglue -6mm
$${\rm id}$$
we make the identification $k = {\mathcal J} k \subset \pi^{-1} k$. Next, we go to the following big commutative diagram, which extends the (\ref{eq2.21})
\begin{equation}
\label{eq2.28}
\xymatrix{k \subset \pi^{-1} k \ar[d] \\
{ \ }}
\xymatrix{
\subset K_1 \ar[r]^{\!\!\!\!\!\!\!\!\!\!\!\!j} \ar[dr]^i  \ar[d]&M^{(3)} \diagup X^3 \ar[d]^g \ar@{>>}[r] &(M^{(3)} \diagup X^3) \diagup \Psi (\pi \circ g) \ar[d]^{g_1} \\ 
{ \ \atop \ }&\Theta^4 \ar[r]^{\pi} &\Theta^3 ({\rm new}) 
} 
\end{equation}
\vglue -2mm
$$
\quad\qquad\qquad\mbox{\hglue -25mm}\pi^{-1} \diagup \Psi (\pi \mid K_1) \ \subset \ K_1 \diagup \Psi (\pi \mid K_1) \quad \overset{\rm inclusion}{-\!\!\!-\!\!\!-\!\!\!-\!\!\!-\!\!\!-\!\!\!-\!\!\!-\!\!\!-\!\!\!-\!\!\!-\!\!\!\longrightarrow} \quad (M^3 \diagup X^3) \diagup \Psi (\pi \circ g) \, .
$$
$\bigl[$Remember here that the lower inclusion follows from the following fact, itself a consequence of (\ref{eq2.25}), namely that
$$
K_1 \diagup \Psi (\pi \circ g) = K_1 \diagup \Psi (\pi \mid K_1) \, . \bigl]
$$
Inside $\Theta^3 ({\rm new})$, we have
\begin{eqnarray}
\pi \circ g \, (M^{(3)} \diagup X^3 - K_1) &= &g_1 \left[ (M^{(3)} \diagup X^3) \diagup \Psi (\pi \circ g) - K_1 \diagup \Psi (\pi \mid K_1) \right] \nonumber \\
&= &\pi \left\{ [g(M^{(3)} \diagup X^3)] \diagup \Psi (\pi \circ g) - K_1 \diagup \Psi (\pi \mid K_1) \right\} \, . \nonumber
\end{eqnarray}
Here, the map $K_1 \diagup \Psi (\pi \mid K_1) \longrightarrow \Theta^3 ({\rm new})$, clearly factors through $\pi K_1 \subset \Theta^3 ({\rm new})$ and invoking (\ref{eq2.25}) we can see that we also have, inside $(M^{(3)} \diagup X^3) \diagup \Psi (\pi \circ g)$, the following
\begin{equation}
\label{eq2.29}
M_2 (g_1) \cap \left( K_1 \diagup \Psi (\pi \mid K_1) \right) = \emptyset \, .
\end{equation}
By (\ref{eq2.26}) $+$ (\ref{eq2.27}) the inclusion $k \subset (M^{(3)} \diagup X^3) \diagup \Psi (\pi \circ g)$ from the proved part of the Claim (2.24), factors through $K_1 \diagup \Psi (\pi \mid K_1)$. Hence, the (\ref{eq2.29}) implies the desired Dehn property $M_2 (g_1) \cap k = \emptyset$.

\smallskip

Our Claim (2.24) has been completely proved, and Lemma 2.4 follows now from (\ref{eq2.23}) $+$ (2.23.1) + (2.24). In the diagram (\ref{eq2.7}) for $k$, take now
$$
K^3 = (M^{(3)} \diagup X^3) \diagup \Psi (\pi \circ g) \, , \quad {\rm and} \quad \chi = g_1 \, .
$$

\newpage

\section{The proof of the multigame Lemma 1.5}\label{sec3}
\setcounter{equation}{0}

In the context of \cite{39} a set of
\setcounter{equation}{-1}
\begin{eqnarray}
\label{eq3.0}
\{\mbox{Holes}\} &= &\{\mbox{completely normal Holes, contained inside walls $W$ (non-complementary)}\} \\
&+ &\{\mbox{BLACK Holes}\} + \{H(p_{\infty\infty})\}  \nonumber
\end{eqnarray}
has been defined, see the beginning of Section IV of \cite{39}. We will introduce now a rather similar (but not quite) set of holes, with a very different utility than the one of the Holes in \cite{39}, namely the
\begin{eqnarray}
\label{eq3.1}
\{\mbox{New Holes}\} &= &\{\mbox{All the completely normal Holes from (\ref{eq3.0}), {\ibf NOT} concerning walls in $Y^2$ (\ref{eq1.23})}\} \nonumber \\
&+ &\{\mbox{one Hole for each $W$ (BLACK complete), {\ibf NOT} in $Y^2$}\} \, . 
\end{eqnarray}
The last item in (\ref{eq3.1}) is independent of the BLACK Holes of \cite{39}, although the two items have a large common intersection.

\bigskip

\noindent {\bf Lemma 3.1.} 1) {\it The multigame from {\rm (\ref{eq1.28})} consists of the following two kinds of steps: we delete all the new Holes above and we also add the} $\underset{n}{\sum} \, {\mathcal B}_n \times [0,\infty)$.

\medskip

2) {\it There is a big $2^{\rm d}$ collapse (to be made explicit later)}
\begin{equation}
\label{eq3.2}
\left\{ \Theta^3 ({\rm new}) - \sum_n {\mathcal B}_n \times [0,\infty)\right\} (\mbox{{\it see} (\ref{eq1.31})}) \supset \Theta^3 ({\rm provisional}) \cup (fX^2 - \{\mbox{new Holes}\}) \overset{\pi (2)}{-\!\!\!-\!\!\!-\!\!\!\longrightarrow} \Theta^3 ({\rm provisional}).
\end{equation}

3) {\it We have natural isomorphisms, at $\pi_1$-level}
$$
\pi_1 \left[ \Theta^3 ({\rm provisional}) \cup (fX^2 - \{\mbox{new Holes}\})\right] \approx \pi_1 (fX^2 - \{\mbox{new Holes}\}) \approx \pi_1 \, fX^2 = 0 \, .
$$

4) {\it There exists also a $3^{\rm d}$ collapse
$$
g(\infty) \, Y(\infty) \overset{\pi (3)}{-\!\!\!-\!\!\!-\!\!\!\longrightarrow} \Theta^3 ({\rm provisional}) \, ,
$$
and} $\{\mbox{new Holes} \} \approx \{\mbox{the 2-cells killed by $\pi(3)$}\}$.

\bigskip

The collapse $\pi(3)$ is essentially the following
$$
\Theta^3 (fX^2)_{\rm II} \cup \sum {\mathcal B} \times [0,\infty) + \sum \mbox{all $H_n^3$'s} \longrightarrow \Theta^3 ({\rm new}) \overset{(1.29)}{-\!\!\!-\!\!\!-\!\!\!\longrightarrow} \Theta^3 ({\rm provisional}) \, .
$$
It is the first half of the collapse above which kills the New Holes. Notice that the combination of (\ref{eq3.1}) with 1) in our lemma already prescribes all the individual elementary games, so it only remains to determine their order.

\smallskip

Now, enough has been said in Section I, to make it clear that, provided we choose the order correctly, namely use the preliminary cleaning before each elementary game, then
$$
{\rm BLUE} < {\rm RED} < {\rm BLACK} \, ,
$$
and provided we only perform one single elementary game, at a given time, we can realize the kind of things which are stated in Lemma 3.1. Figure 3.1 provides a toy-model which should suggest what the multigame does.

$$
\includegraphics[width=145mm]{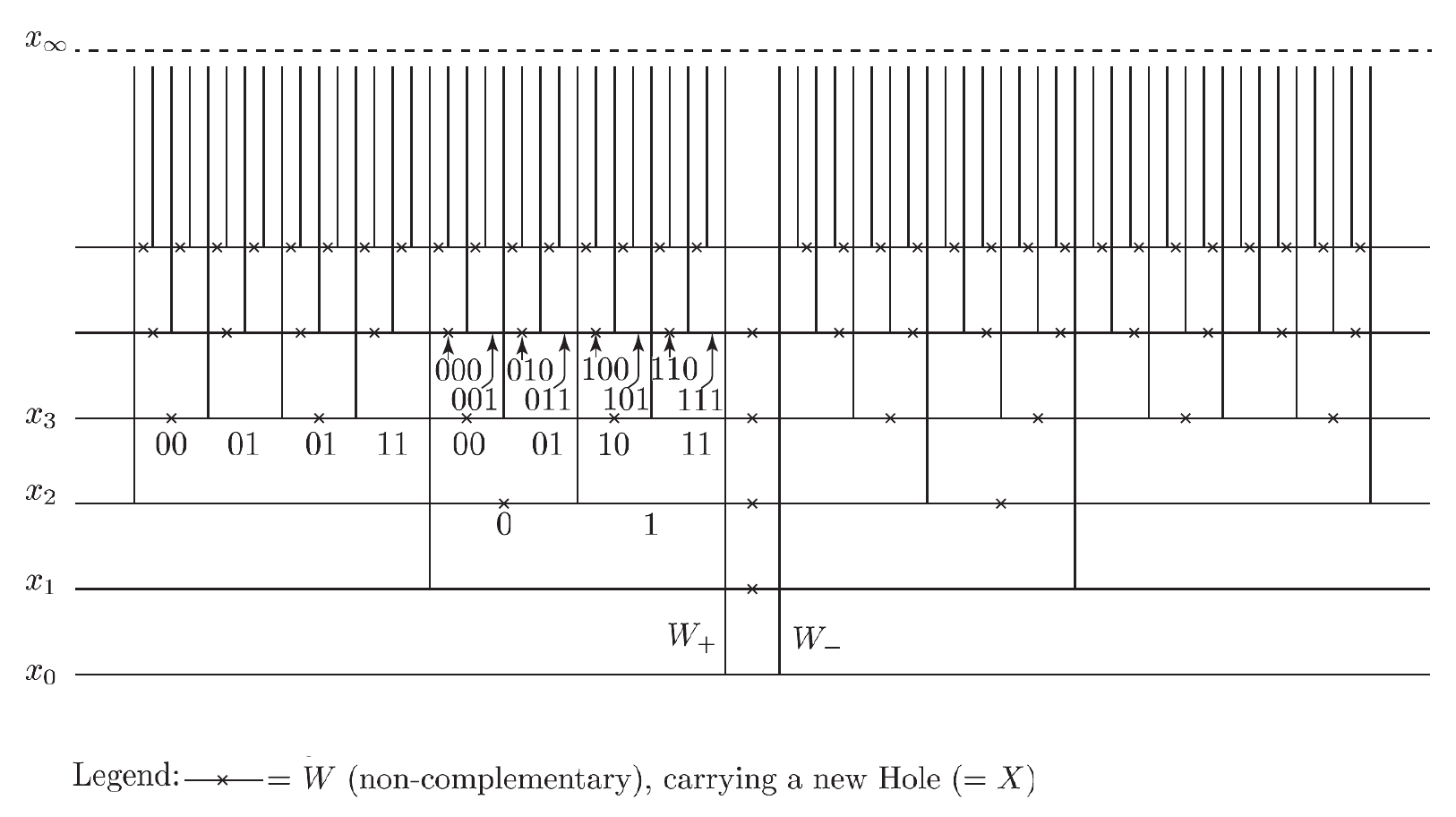}
$$
\label{fig3.1}
\centerline {\bf Figure 3.1.} 

\smallskip

\begin{quote}
Schematical, toy-model representation of the multigame.
\end{quote}

\bigskip

In Figure 3.1, the $Y^2$ is located at $x_0 \geq x$, at $x = x_{\infty}$ we have a limit wall and, moreover the following pattern has been set up. We see various non-complementary walls $W$ (see here (4.16) in \cite{39} for the distinction between complementary and non-complementary walls) labelled by $\{ \varepsilon_1 , \varepsilon_2 , \ldots , \varepsilon_n \}$, when each $\varepsilon_i$ is $0$ or $1$. It is exactly the $W (\varepsilon_1 , \varepsilon_2 , \ldots , \varepsilon_{n-1} , \varepsilon_n = 0)$ which carry a $\{$New Hole$\}$ while, remember, each non-complementary wall carries a completely normal Hole. All this is schematical, of course.  The figure may also suggest the collapse (\ref{eq3.2}).

\smallskip

The main item of this section is the following statement which completes the Lemma 1.5, and the proof of which will also yield the proof of our multigame Lemma 1.5.

\bigskip

\noindent {\bf Lemma 3.2.} 1) {\it There is an equivariant codimension one space, which is a surface with branching lines (locally like $Y \times R$), PROPERLY and properly embedded
\begin{equation}
\label{eq3.3}
(S_0 , \partial S_0) \subset \left(\Theta^3 ({\rm new}) - \sum_{{\mathcal B}_n} {\mathcal B}_n \times (0,\infty) , \partial \, \Theta^3 ({\rm new}) \right) \, ,
\end{equation}
which meets the $\underset{n}{\sum} \, {\mathcal B}_n$ transversally, and which induces the following {\ibf splitting}, via which $\Theta^3 (\mbox{\rm co-compact})$ from {\rm (\ref{eq1.32})} is {\ibf defined} (once the surface $S_0$ has been explicitly specified)
\begin{equation}
\label{eq3.4}
\Theta^3 ({\rm new}) - \sum_n {\mathcal B}_n \times (0,\infty) = \Theta^3 (\mbox{\rm co-compact}) \underset{\overbrace{\mbox{\footnotesize$S_0$}}}{\cup} \Theta^3_0 \, (\mbox{\rm residual space}) \, ,
\end{equation}
coming with the inclusions, prescribed by Lemma {\rm 1.5},
\begin{equation}
\label{eq3.5}
\Theta^3 (\mbox{\rm co-compact}) \subset \Theta^3 ({\rm provisional}) \subset \Theta^3 \, ({\rm new}) - \sum_n {\mathcal B}_n \times (0,\infty) 
\end{equation}
which also determines on which side of $S_0$ the $\Theta^3 (\mbox{\rm co-compact})$ is located. Figure {\rm 3.3} can serve as a first approximate description of the first inclusion in} (\ref{eq3.5}).

\medskip

2) {\it For each ${\mathcal B}$, each of the non-void intersections ${\mathcal B} \cap \Theta^3 (\mbox{\rm co-compact})$ is a compact {\ibf connected} surface, with non-empty boundary.}

\medskip

3) {\it Let $L = {\mathcal B} ({\rm BLACK}) \cap {\mathcal B} ({\rm BLACK})$ be one of the lines of tranversal intersection from $1)$ in Lemma {\rm  1.5.} Then, we also have}
$$
L \cap \Theta^3 (\mbox{co-compact}) = \emptyset \, . \eqno (3.5.1)
$$

4) {\it The $\partial \, S_0$ has sufficiently many connected components so that we can find a proper and PROPER embedding of disjoined finite trees
\begin{equation}
\label{eq3.6}
\sum_i A_i \subset S_0 \, , \ \mbox{s.t. the following things happen:}
\end{equation}
}

4.1. {\it The ramifications of each $A_i$ reflect exactly the intersections of $A_i$ with the ramifications of $S_0$, so that $\underset{i}{\sum} \, A_i \subset S_0$ induces a clean codimension one splitting of $S_0$.}

\smallskip

4.2. {\it This splitting break $S_0$ into a disjoined union of compact collapsible pieces}
\begin{equation}
\label{eq3.7}
S_0 = \sum_j B_j \, .
\end{equation}

5) {\it There is a $3^{\rm d}$ collapse}
\begin{equation}
\label{eq3.8}
\Theta_0^3 \overset{\pi}{-\!\!\!-\!\!\!-\!\!\!\longrightarrow} S_0 \, .
\end{equation}

\bigskip

The big collapse (\ref{eq1.29}) from Lemma 1.5 reduces, essentially, to the collapse from (\ref{eq3.8}) above.

\bigskip

\noindent {\bf Complements to Lemma 3.2.} 1) {\it Generically (meaning when outside things like the immortal singularities), along a ${\mathcal B} \times \{ 0 \}$, the $\Theta^3 ({\rm new})$ is like a figure $Y$, see Figure {\rm 1.6} for an illustration. The $\Theta^3 (\mbox{\rm co-compact})$ goes through ${\mathcal B} \times \{ 0 \}$, without entering the ${\mathcal B} \times [0,\infty)$ arm of the $Y$ in question.}

\medskip

2) {\it From now on, $\Sigma (\infty)$ will be like in {\rm (2.13.1)} in {\rm \cite{39}}, with all the  contribution of $p_{\infty\infty} (S)$ deleted. This comes with ${\rm int} \, \Sigma (\infty) \subset \Sigma (\infty)$, which occurs in {\rm (1.1.bis)}. The $\Sigma (\infty)$ comes with a second surface}
$$
\Sigma (\infty) (\mbox{co-compact}) \equiv \Sigma (\infty) \cap \Theta^3 (\mbox{co-compact}) \, ,
\eqno (3.8.1)
$$
{\it and, modulo $\Gamma$, this is a surface of finite type, except that finitely many arcs, properly embedded
$$
\sum_{\{p_{\infty\infty} (S)\} \diagup \Gamma} p_{\infty\infty} (S) \times [-\varepsilon , \varepsilon]
$$
have been deleted (leaving us with punctures).

\smallskip

We also have}
\begin{equation}
\label{eq3.9}
\partial \, \Sigma (\infty) \cap \Sigma (\infty) (\mbox{co-compact}) = \emptyset \, ,
\end{equation}
{\it and it is the {\ibf ramification} of $\Theta^3 (\mbox{\rm co-compact})$, actually making that the natural free action
$$
\Gamma \times \Theta^3 (\mbox{\rm co-compact}) \longrightarrow \Theta^3 (\mbox{\rm co-compact})
$$
{\ibf is} co-compact (i.e. has a compact fundamental domain), which is the main reason for the splitting} (\ref{eq3.4}). {\it Let us say that $\Theta^3 (\mbox{\rm co-compact})$ looks very much like $\Theta^3 ({\rm provisional})$, {\ibf BUT} with the big difference that, while the free action of $\Gamma$ on $\Theta^3 ({\rm provisional})$ is {\ibf not} co-compact, the free action of $\Gamma$ on $\Theta^3 (\mbox{\rm co-compact})$ {\ibf is}.}

\medskip

3) {\it Remember, at this point, that we find
$$
{\rm int} \, \Sigma (\infty) \times [0,\infty) \subset \Theta^3 (fX^2)_{\rm II} \cap \Theta^3 ({\rm new}) \, ,
$$
and independently of things like immortal singularities or ${\mathcal B}$'s, the two $\Theta^3$'s, new and provisional, fail already to be $3$-manifolds along the ${\rm int} \, \Sigma (\infty) \times [0,\infty)$. Figure {\rm 3.2} describes the interaction of $S_0$ with ${\rm int} \, \Sigma (\infty) \times \{0\}$ outside the immortal singularities and the ${\mathcal B}$'s.}

$$
\includegraphics[width=7cm]{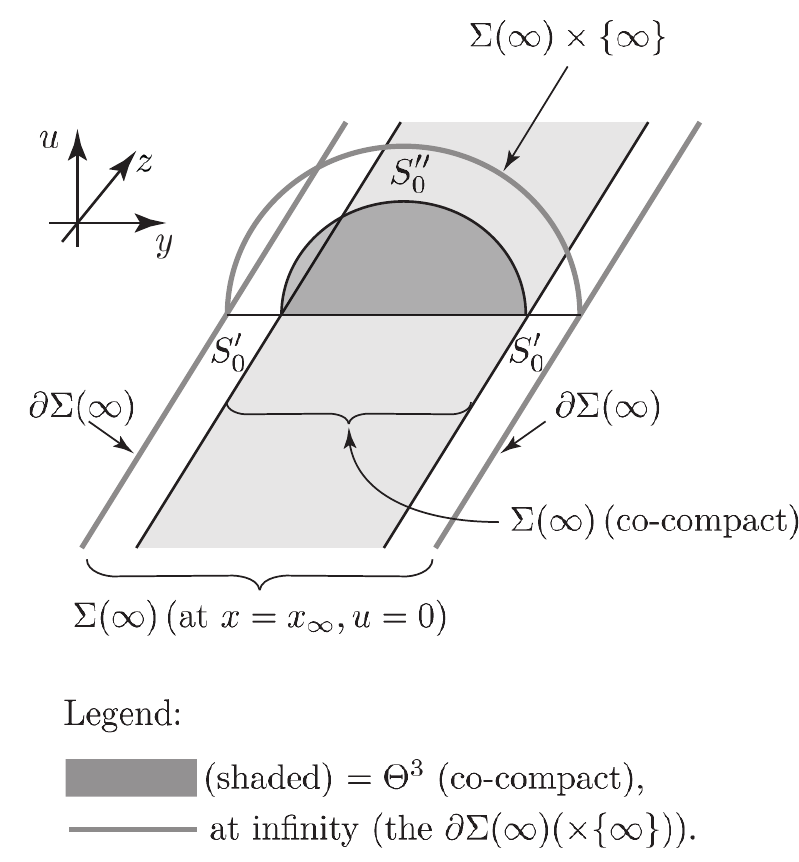}
$$
\label{fig3.2}
\centerline {\bf Figure 3.2.} 

\smallskip

\begin{quote}
We see here a small detail of ${\rm int} \, \Sigma (\infty) \times [0,\infty) \subset \Theta^3 ({\rm new})$, in a simple location far from immortal singularities, ${\mathcal B}$'s or bifurcations of $\Sigma (\infty)$. The splitting is $S_0 = S'_0 \cup S''_0$ from the formula (\ref{eq3.10}) below. 

\smallskip

The bare local coordinate system concerns $\Theta^3 (fX^2)_{\rm II}$. This restricts to $(y,z)$ along $\Sigma (\infty)$ and the axis $u$ has been added for ${\rm int} \, \Sigma (\infty) \times [0,\infty)$. The $\Sigma (\infty)$ lives at $x=x_{\infty}$ and, while $S'_0$ continues along $-M \leq x-x_{\infty} \leq M$, the $S''_0$ continues along $-N \leq z \leq N$.
\end{quote}

\bigskip

4) {\it Figure {\rm 3.2} should give an idea about $S_0 \cap (\Sigma (\infty) \times [0,\infty))$ and it should also suggest a decomposition}
\begin{equation}
\label{eq3.10}
S_0 = S'_0 \underset{\overbrace{\mbox{\footnotesize$\partial \, \Sigma (\infty)$ (co-compact)}}}{\cup} S''_0 \, , \quad {\rm with} \quad S''_0 \equiv S_0 \cap \left( \Sigma (\infty) \times [0,\infty) \right) \, .
\end{equation}

5) {\it Concerning the collapse $\pi$ from {\rm (\ref{eq3.8})}, each $\pi^{-1} A_i$ collapses into $A_i$, the $\underset{i}{\sum} \, \pi^{-1} A_i$ breaks $\Theta_0^3$ into $\underset{j}{\sum} \, \pi^{-1} B_j$ and each $\pi^{-1} B_j$ collapses into $B_j$.}

\bigskip

Figure 3.3 should help understand the splitting (\ref{eq3.4}) in the neighbourhood of $\Sigma (\infty)$. This figure is supposed to have a good fit with Figure 3.2.

$$
\includegraphics[width=13cm]{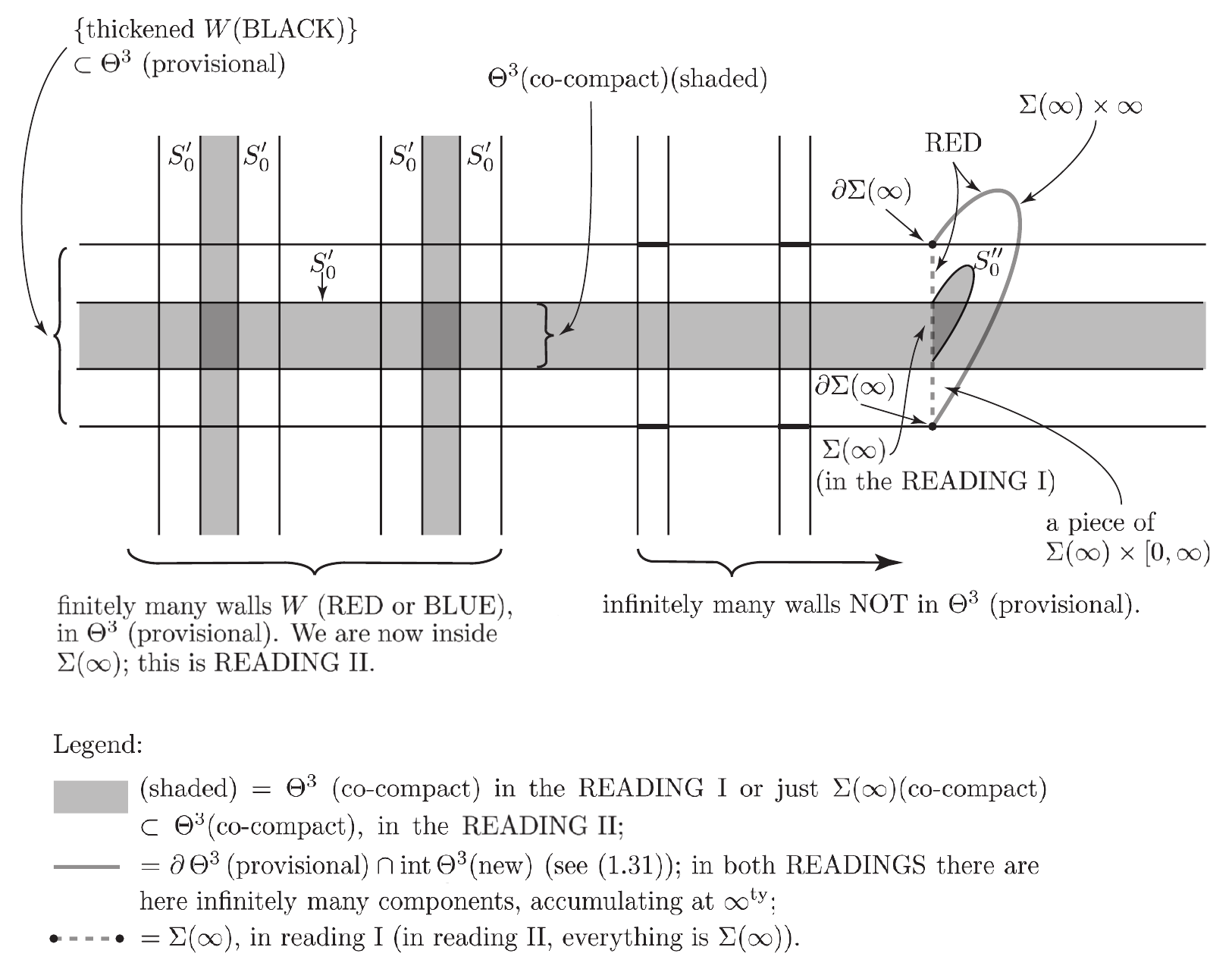}
$$
\label{fig3.3}
\centerline {\bf Figure 3.3.} 

\smallskip

\begin{quote}
This figure, which is very much in the style of Figure 3.2 should help understand the articulation between $\Theta^3 ({\rm provisional}) \supset \Theta^3 (\mbox{co-compact})$ and $\Sigma (\infty)$. There are two readings for the present figure, a READING I where the plane of the figure meets transversally the $W({\rm BLACK}) \subset \Theta^3 ({\rm provisional})$ and then also the $\Sigma (\infty)$. This last item is then the red detail at the right side of our figure. We are here far from the $p_{\infty\infty}$'s.

There is also a READING II for our figure, where what we see is $\Sigma (\infty)$ itself, the plane of the figure being then $x=x_{\infty}$.
\end{quote}

\bigskip

The $S'_0 \subset S_0$ in (\ref{eq3.10}), is essentially a copy of (see (\ref{eq1.27}))
$$
\partial \, \Theta^3 (Y^2) = \partial \left(\Theta^3 ({\rm provisional}) - \ring\Sigma \, (\infty) \times [0,\infty)\right) \, ,
$$
pushed towards the interior of $\Theta^3 ({\rm provisional})$. The $S''_0$ is another copy of $\Sigma (\infty) (\mbox{co-compact})$.

\smallskip

To fully understand the $S_0$ which defines our $\Theta^3 (\mbox{co-compact})$ we need to make precise its structure when localized at the immortal singularities and at the ${\mathcal B} ({\rm BLACK})$'s.

\smallskip

I will explain now how the double points of the map ${\mathcal J}$ (\ref{eq1.30}) appear. When one compares the Figures 1.4 and 1.6, one sees the following. The ${\mathcal B} (W) - D^2 (p_{\infty\infty} (S))$ is a surface of boundary $C(p_{\infty\infty} (S))$ in Figure 1.4, occupying the rest of $W({\rm BLACK})$ and then climbing up the $z$-axis (which is perpendicular to the plane of Figure 1.4 and looks towards the observer, like in the generic Figure 1.2), towards the $S_{\infty}^1$, living at the infinity of ${\mathcal B} (W)$. Similarly, the
$$
{\mathcal B} (W_1^*) - D^2 (p_{\infty\infty} (S^*)) ({\rm of} \ W_1^*) \quad (\mbox{Figure 1.6})
$$
leaves $C(p_{\infty\infty} (S^*))$ and cuts through the green line marked ${\mathcal B} (W_1^*)$ in Figure 1.4, the ${\mathcal B} (W)$. The green line in question is actually the beginning of the
$$
L = {\mathcal B} (W) \cap {\mathcal B} (W_1^*) \underset{\rm TOP}{=} R_+ \, ,
$$
which starts at the point $\partial L \in C(p_{\infty\infty} (S^*)) \times \left[ -\frac\varepsilon4 , \frac\varepsilon4 \right]$ in the Figure 1.4, and then goes in the direction $-x$ in that figure and also in the Figures 3.4 and 3.5. These figures will be our starting point for explaining what $S_0$ and the $\Theta^3 (\mbox{co-compact})$ which cobounds it look like in the neighbourhood of the immortal singularities. In the figure in question the $\Theta^3 (\mbox{co-compact})$ has been shaded, and the
$$
D^2 (p_{\infty\infty} (S)) \subset \Theta^3 (\mbox{co-compact})
$$
doubly shaded. When we say that Figures 3.4 $+$ 3.5 are only a ``starting point'' for defining $S_0$ ($=$ the splitting hypersurface), we have in mind the following fact, which is actually a flaw:
\begin{equation}
\label{eq3.11}
\mbox{The $S_0$ which, simple-mindedly, they may suggest, does {\ibf NOT} split.}
\end{equation}
Here for the clarity of our exposition, the immortal $S$ from Figure 3.5 has been re-drawn in Figure 3.6. Let us introduce the notation
$$
\sigma \equiv \{\mbox{the shaded area} \ [a_1 , b_1 , b_2 , c_1 , c_2 , d_1 , d_2 , a_2] = \{\mbox{singularity of} \ \Theta^3 (Y^2)\} , \ \mbox{in the Figure 3.6}\} \, .
$$
As things have been drawn so far, the branch (see (\ref{eq1.31}))
$$
\left(\ring\Sigma \, (\infty) \times [0,\infty)\right) \mid \sigma \subset \Theta^3 ({\rm provisional}) \subset \Theta^3 ({\rm new})
$$
has not yet been properly taken care of. The cure is to combine the treatment of the Figures 3.4 $+$ 3.5 with the one in the Figure 3.2. This entails enriching our $S_0$ with a piece $S''_0$ like in that figure, cutting through ${\rm int} \, \Sigma (\infty) \times [0,\infty)$ and isolating a piece which goes together with $\Theta^3 (\mbox{co-compact})$. Notice that this will also create for the immortal singularities of $\Theta^3 (\mbox{co-compact})$ a number of branches superior to the canonical two for the undrawable singularities in \cite{8}, \cite{19}, \cite{36}, and then accordingly, more complicated desingularizations ${\mathcal R}$. But no harm comes with this. We just have to live with a $\Theta^3 (\mbox{co-compact})$ which is train-track, generically with three $R_+^3$ branches, which may become four when in the presence of $\Sigma(\infty) \times [0,\infty)$.

$$
\includegraphics[width=13cm]{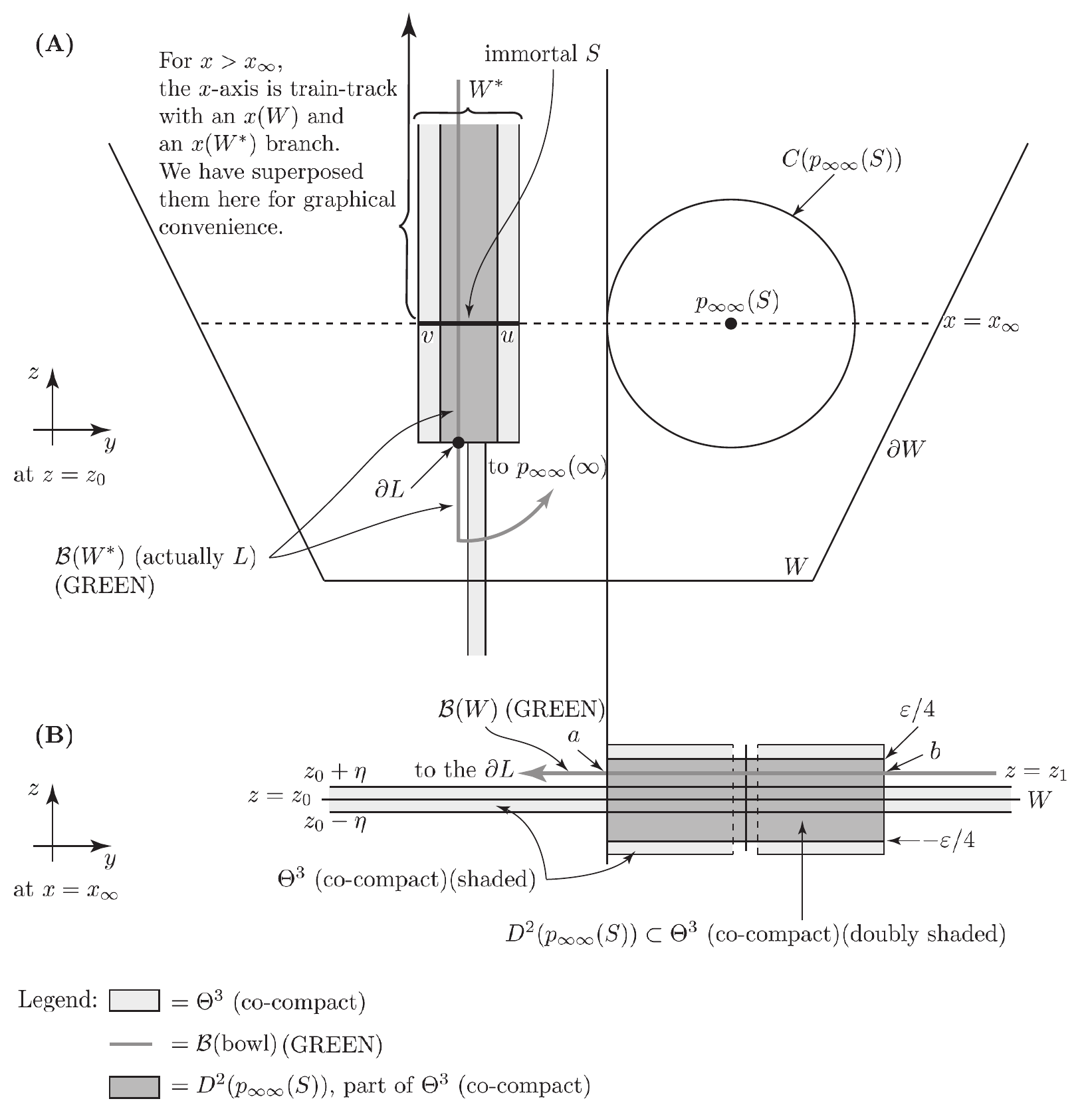}
$$
\label{fig3.4}
\centerline {\bf Figure 3.4.} 

\smallskip

\begin{quote}
This figure should be compared with Figure 1.4 and it concerns the situation when $W+W^* \subset Y^2$ and when at the surviving $S$ (see here the (\ref{eq1.35})), $W^*$ overflows while $W$ is subdued. We have shaded the $\Theta^3 (\mbox{co-compact})$. The $D^2 (p_{\infty\infty} (S,S^*) (W \ {\rm and} \ W^*))$ which are both inside $\Theta^3 (\mbox{co-compact})$ are doubly shaded. All the Figure (A), when outside the $C(p_{\infty\infty} (S))$ is covered by ${\mathcal B} (W)$. The figure should explain how, in the considered situation, the double line
$$
{\mathcal J} {\mathcal B} (W) \cap {\mathcal J} {\mathcal B} (W^*) \subset M^2({\mathcal J})\, ,
$$
with ${\mathcal J}$ like in  (\ref{eq1.30}) starts at $\partial L$ (in (A)) staying {\ibf disjoined} from $\Theta^3 (\mbox{co-compact})$, as it is stipulated by (3.5.1). This conclusion will remain valid when $W^*$ is subdued and $W$ overflows; see here the Figure 3.5.

At $z=z_0$, in (A) the $W$ is shaded (i.e. it is in $\Theta^3 (\mbox{co-compact})$), while at $z=z_1$, the same (A), when outside $C(p_{\infty\infty} (S))$ is green, like ${\mathcal B}$; then ${\mathcal B} (W)$ covers it.
\end{quote}

\bigskip

\noindent {\bf Additional explanations for the Figures 3.4 and 3.5.} Here $z_0$ stands for the $\{z=0$ of the wall $W({\rm BLACK})\}$. The thickness of $\Theta^3 (fX^2)_{\rm II} \mid W$ is $2\varepsilon$ and the thickness of $\Theta^3 (\mbox{co-compact}) \mid W$ is $2\eta$. We have
$$
z_0 + \varepsilon > z_0 + \frac\varepsilon\eta \gg z_1 (\mbox{level of ${\mathcal B} (W)$}) > z_0 + \eta > z_0 (\mbox{level of $W$}) > z_0 - \eta \, .
$$
In (A) we have shaded only $\Theta^3 (\mbox{co-compact}) \mid W^*$ but, in real life, all the thickness $z_0 - \eta \leq z \leq z_0 + \eta$ is shaded, i.e. {\ibf all} of (A). Our inequalities above make that $L = {\mathcal J} {\mathcal B} (W) \pitchfork {\mathcal J} {\mathcal B} (W^*)$ lives at $z=z_1$, outside of $\Theta^3 (\mbox{co-compact})$, at least  in the seable region. Along the line $[ab]$, the ${\mathcal B} (W)$ rides on top of the $D^2 (p_{\infty\infty} (S))$. In both this figure and in 3.5, the $u,v$ refer to the immortal singularity. Here $[uv] = \{$an immortal $S$ projected on $W_{\rm subdued}\}$. In Figure 3.5,
$$
[uv] \subset S \subset W_{\rm overflowing} \cap \Theta^3 (\mbox{co-compact}) \, .
$$
End of explanations.

$$
\includegraphics[width=125mm]{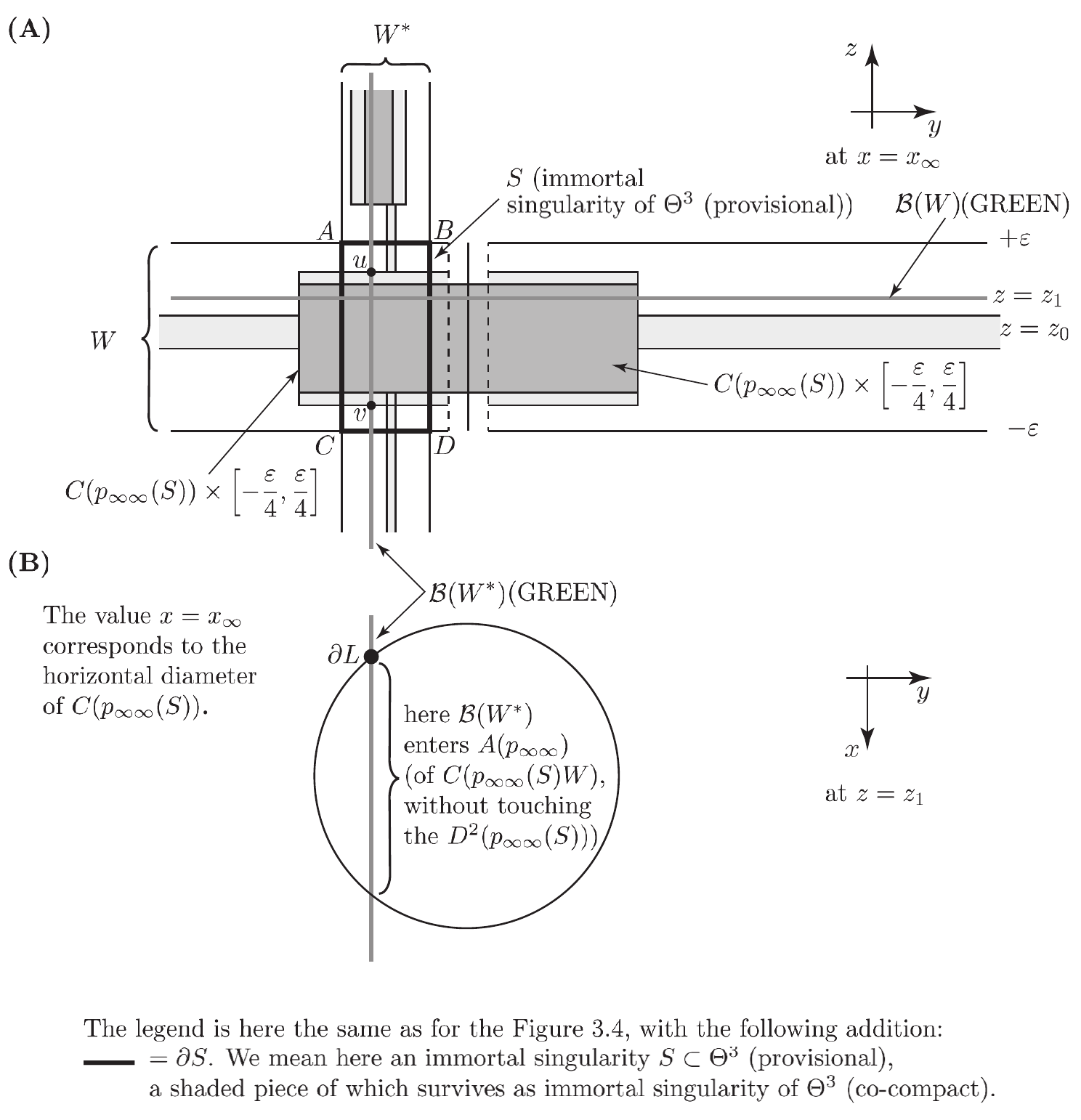}
$$
\label{fig3.5}
\centerline {\bf Figure 3.5.} 

\smallskip

\begin{quote}
This figure is in the same style as 3.5, but it refers now to Figure 1.6. Here $W$ overflows and $W^*$ is subdued. Inside the (red) contour marked $S$, we have the unique surviving $S \subset \Theta^3 ({\rm new})$, from some $\overline S$. The whole of (B), outside $C(p_{\infty\infty} (S))$ is covered by ${\mathcal B} (W)$. This explains, in the present situation, the $L \subset M^2({\mathcal J})$, $L = {\mathcal B} (W) \cap {\mathcal B} (W^*)$, which is starting at $\partial L$ and is {\ibf not} touching the shaded $\Theta^3 (\mbox{co-compact})$.

The Figures 3.4, 3.5 do not focus on the same items but they both refer to a single situation, an immortal $S$ like in (\ref{eq1.34}), (\ref{eq1.35}), and a pair of black walls in duality overflowing/subdued. There is a perfect duality between $W$ and $W^*$.
\end{quote}

\bigskip

Concerning now the Figure 3.6, notice the four white corners, like the $(A,a_1,a_2)$. There are immortal singularities for the $\Theta_0^3$ in (\ref{eq3.4}).

$$
\includegraphics[width=12cm]{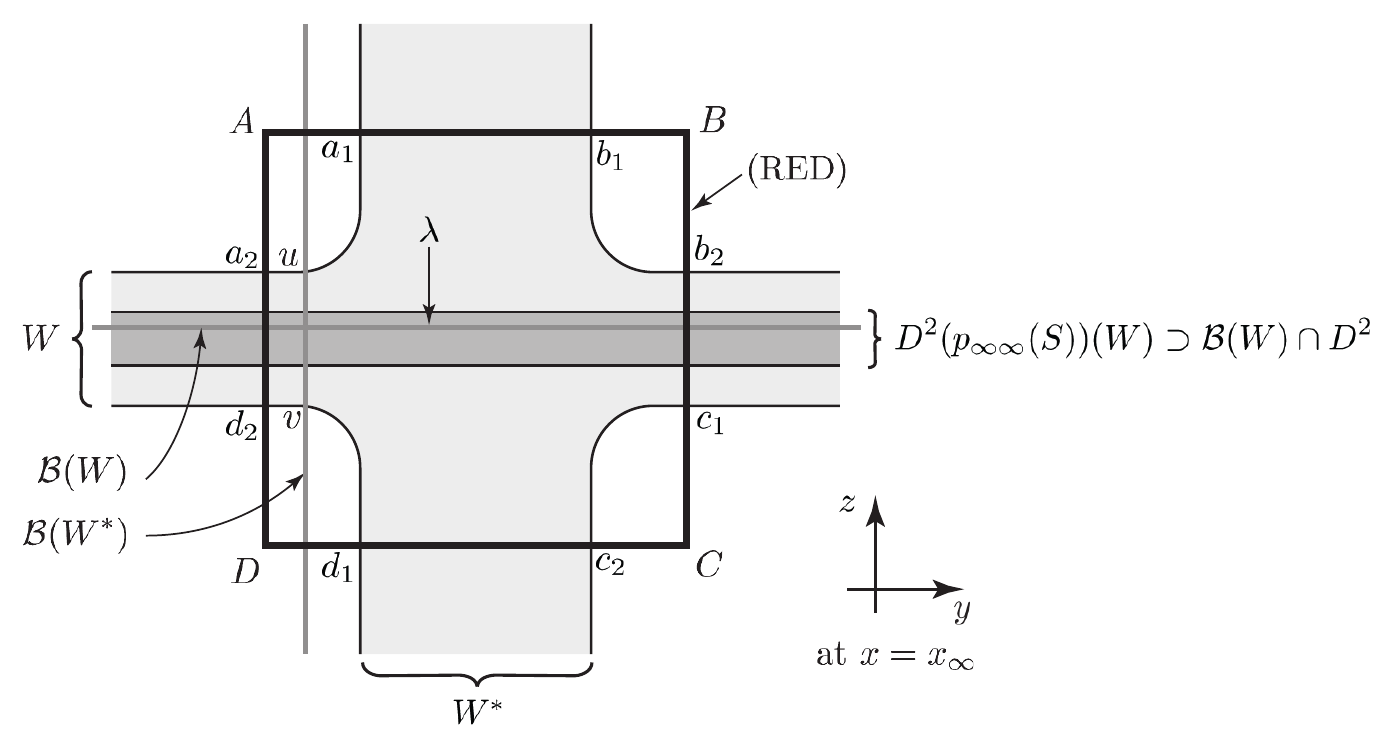}
$$
\label{fig3.6}
\centerline {\bf Figure 3.6.} 

\smallskip

\begin{quote}
A detailed view of the singularity $S$ from the Figure 3.5. We see here $\sigma = [a_1 , b_1 , b_2 , c_1 , c_2 , d_1 , d_2 , a_2]$, immortal singularity of $\Theta^3 (\mbox{co-compact})$.

The ${\mathcal B} (W) , {\mathcal B} (W^*)$, which are GREEN, never meet at the level of this figure, which lives at $x = x_{\infty}$. The point $\lambda$ is ficticious.
\end{quote}

\bigskip

In order to describe, abstractly, a bowl ${\mathcal B} = {\mathcal B} (W({\rm BLACK}))$, we start by thinking in terms of the following decomposition for $W = W({\rm BLACK})$ itself, $W = \{$a central $D^2$ which, in terms of the Figure 1.1 in \cite{39} is bounded by the dotted hexagon with vertices $p_{\infty\infty}\} \cup \{$a collar piece $\partial W \times [1,0) = W-D^2\}$, with $\partial D^2 = \partial W \times \{0\}$, $\partial W = \partial W \times \{1\}$. With this, here is the abstract description of ${\mathcal B} = {\mathcal B}(W)$
\begin{eqnarray}
\label{eq3.12}
{\mathcal B}(W) &\cong &D^2 \underset{\overbrace{\mbox{\footnotesize$\partial W \times \{0\}$}}}{\cup} \partial W \times [0,1] \underset{\overbrace{\mbox{\footnotesize$\partial W \times \{1\} = \partial W$}}}{\cup} \partial W \times [1,0) \\
&= &W \cup \partial W \times [1,0) = D^2 \cup \{\mbox{a {\ibf collar} piece} \ \partial D^2 \times [1,0)\} \, , \nonumber
\end{eqnarray}
with the last $\partial D^2 \times \{0\}$ living at infinity.

\smallskip

The Figure 3.7 below, where the $D^2 (W) , D^2 (W^*)$ live on the other side of the square $[ABCD] \subset S_{\infty}^2$ should help understand the following two features from the Lemma 3.2: we have both a {\ibf connected} ${\mathcal B} \cap \Theta^3 (\mbox{co-compact})$ and then also, like in (3.5.1), the $L = {\mathcal J} {\mathcal B} ({\rm BLACK}) \pitchfork {\mathcal J} {\mathcal B} ({\rm BLACK})$ is far from $\Theta^3 (\mbox{co-compact})$. On the other hand, the fact that ${\mathcal B} \cap \Theta^3 (\mbox{co-compact})$ {\ibf is} compact, is an immediate consequence of the fact that $\Theta^3 (\mbox{co-compact})$ stays far from the infinity of ${\mathcal B}$.

\smallskip

In order to simplify our discussion, we will ignore, provisionally at least, the fact that the ${\mathcal B}$'s ride on top of the $D^2 (p_{\infty\infty} (S))$'s, like in the Figures 1.6, 3.4, 3.5, 3.6. We will also pretend that $W \cap W^* \subset \Theta^3 (fX^2)_{\rm II}$ is generic with one transversal intersection point for $\partial W \cap {\rm int} \, W^*$, like in the Figures 1.6, 3.4.

\smallskip

We will also ignore, provisionally, the fact that rendering the 2) and the 4) in Lemma 3.2 compatible with each other will require some doctoring, which we leave for later on.

$$
\includegraphics[width=12cm]{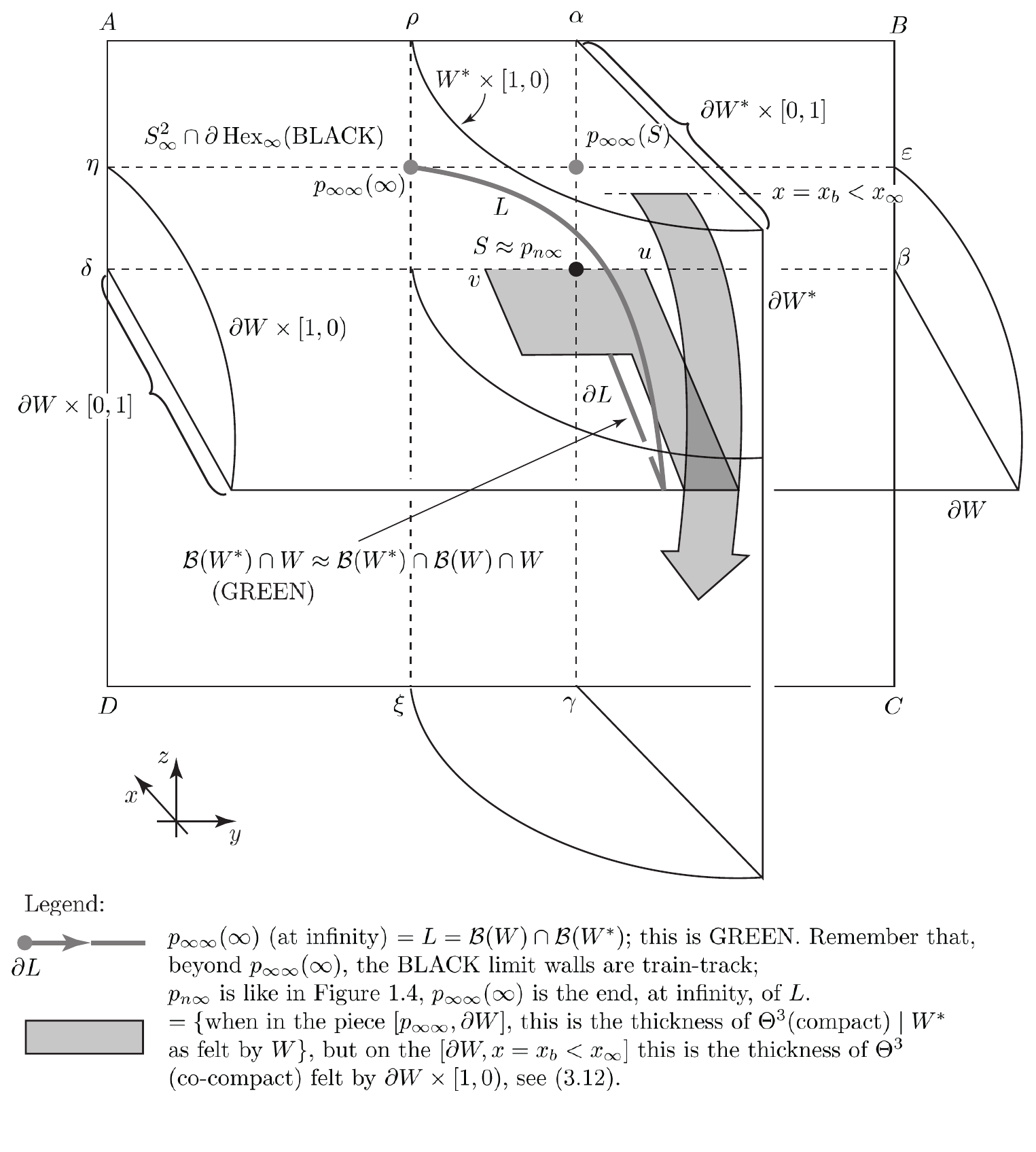}
$$
\label{fig3.7}

\centerline {\bf Figure 3.7.} 

\smallskip

\begin{quote}
The geometry of ${\mathcal B}(W) , {\mathcal B}(W^*)$ and of $L = {\mathcal B}(W) \cap {\mathcal B}(W^*)$. This figure is completely concentrated in the {\ibf NON}-traintrack region $x \leq x_{\infty}$. The $[A,B,C,D] \subset S_{\infty}^2$ (at $x=x_{\infty}$). Beyond $x_{\infty}$, the $x$-axis becomes train-track, branching into $x(W) \geq x_{\infty}$ and $x(W^*) \geq x_{\infty}$. The figure should help understand the contact between the pieces of ${\mathcal B}(W) , {\mathcal B}(W^*)$, each ${\mathcal B}$ being decomposed as in (\ref{eq3.12}). The $[\rho , \xi]$ and $[\delta, \beta]$ are in $\partial {\rm Hex}_{\infty} ({\rm BLACK})$ with $p_{\infty\infty} (\infty)$ an immortal singularity of limit walls, just like $p_{\infty n}$ is an immortal singularity involving $W,W^*$. The two ${\rm Hex}_{\infty} ({\rm BLACK})$'s and $D^2 (W,W^*)$ (see (\ref{eq3.12})) continue beyond $x_{\infty}$ in a train-track manner. To be concretely explicit, think here in terms of the $W^* , W$ in Figure 3.4.
\end{quote}

\bigskip

$$
\includegraphics[width=12cm]{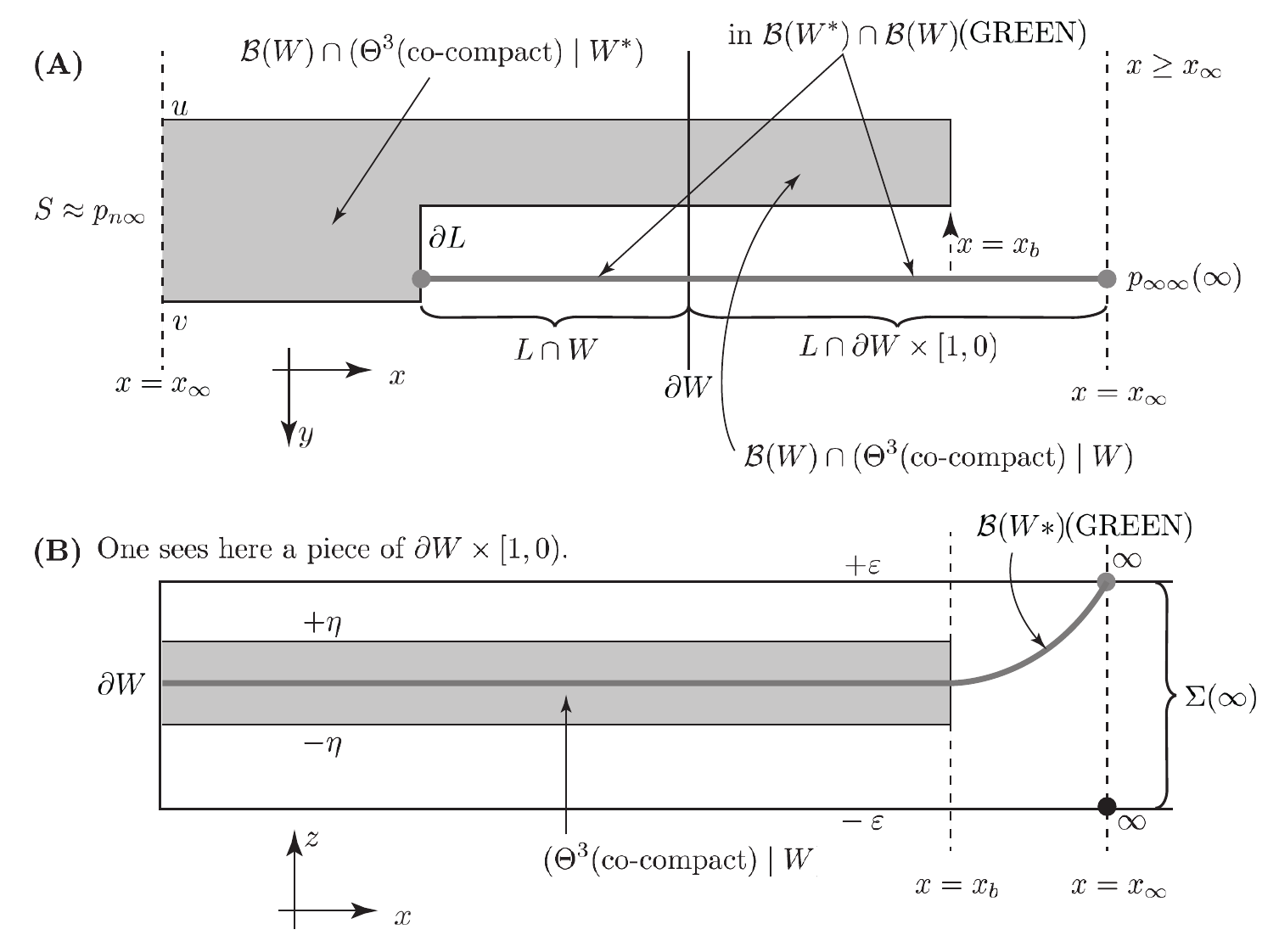}
$$
\label{fig3.7.bis}
\centerline {\bf Figure 3.7.bis.} 

\smallskip

\begin{quote}
This figure continues and completes 3.7. Imagine that (A) represents a piece of ${\mathcal B} (W)$ following very closely $(\partial W \times [0,1]) \cup (\partial W \times [1,0))$. Very importantly, as far as $x \leq x_{\infty}$ and $x \geq x_{\infty}$ in this figure are concerned, at the line $\partial W$ the sign of $x$ switches, and $x \leq x_{\infty}$ (respectively $x \geq x_{\infty}$) becomes $x \geq x_{\infty}$ (respectively $x \leq x_{\infty}$).
\end{quote}

\bigskip

In the Figures 3.7 $+$ 3.7.bis, it is the sanitizing $S_0$, the splitting surface from (\ref{eq3.3}) $+$ (\ref{eq3.4}), more precisely its main branch $S'_0$ from the Figure 3.2, which stops the $\Theta^3 (\mbox{co-compact})$ at $x = x_b < x_{\infty}$, keeping it away from the $\Sigma (\infty)$ at $x=x_{\infty}$.

\smallskip

Look now at the Figures 1.6 and 3.4, which are supposed to account for the $W,W^*$ in the Figures 3.7 $+$ 3.7.bis. To be very concrete, we assume that the $p_{n\infty}$ from Figures 3.7 $+$ 3.7.bis is actually the $p_{1\infty}$ in the Figure 1.6, so that our $S \approx p_{n\infty}$ (Figure 3.7) is actually the $[\alpha,\beta,\gamma,\delta]$ in the LHS of Figure 1.6. Then, when $W \times [-\varepsilon , \varepsilon]$ (Figure 1.6) is collapsed down to $W$, like in Figure 3.4, then $[\delta\gamma]$ become $v$ and $[\alpha\beta]$ become $u$, accounting for the $[u,v]$ which occurs in the Figures 3.4 and 3.7. With this, $\partial L$ (Figures 3.4 and 3.7) rests actually on the circle $C (p_{\infty\infty} (S))(W^*) = \partial (D^2 (p_{\infty\infty} (S)) (W^*))$, more explicitly on $C(p_{\infty\infty} (S) (W^*)) \times (z_1 =$ level of ${\mathcal B} (W^*))$. To the left of $\partial L$, towards $x \geq x_{\infty}$, in (A) Figure 3.7, our ${\mathcal B} (W^*)$ rides on $D^2 (p_{\infty\infty} (S)) (W^*)$, staying at level $z=z_1$, and when it leaves $D^2 (p_{\infty\infty} (S) (W^*))$, then it is outside of $\Theta^3 (\mbox{co-compact})$.

\smallskip

With all this enough has been said concerning the points 2) and 3) in Lemma 3.2, and we turn now to the point 4). Notice, to begin with, that if for $\partial S_0$ (actually for $\partial S'_0$) we manage to create sufficiently many, well-located components, then 4) is true for our $S_0$, essentially for the same reasons which make it true in the realm of smooth surfaces. But then, $S'_0$ runs very closely parallel to $\partial \Theta^3 ({\rm provisional})$, so we can create more $\partial S'_0$ by sending feelers from $\Theta^3 (\mbox{co-compact}) \subset {\rm int} \, \Theta^3 ({\rm provisional})$, to $\partial \Theta^3 ({\rm provisional})$. Figure 3.8 suggests how to do this, without violating the connectivity of each
$$
{\mathcal B}_n \cap \Theta^3 (\mbox{co-compact}) \, .
$$

\bigskip

$$
\includegraphics[width=13cm]{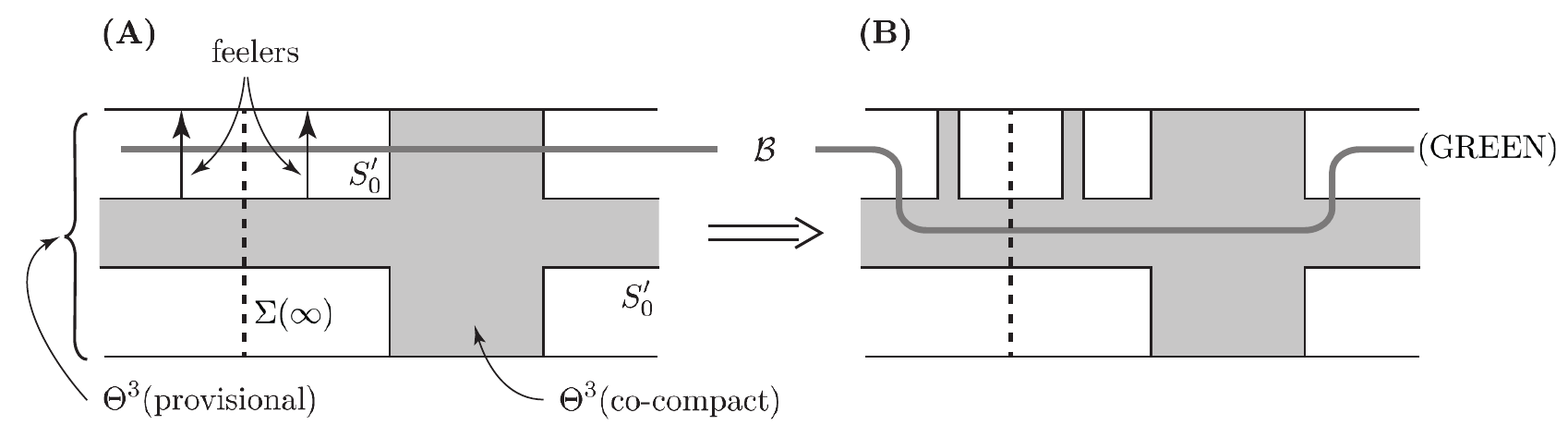}
$$
\label{fig3.8}

\centerline {\bf Figure 3.8.} 

\smallskip

\begin{quote}
How to create more $\partial S_0$, without violating the connectivity of ${\mathcal B}_n \cap \Theta^3 (\mbox{co-compact})$.
\end{quote}

\newpage

\section{From Dehn-exhaustibility to QSF; final arguments}\label{sec4}
\setcounter{equation}{0}

The main object of the present section is to prove the following.

\bigskip

\noindent {\bf Lemma 4.1.} {\it We have the following implication
\begin{equation}
\label{eq4.1}
\{\Theta^3 ({\rm new}) \ \mbox{is $3^{\rm d}$ Dehn-exhaustible}\} \Longrightarrow \{\Theta^3 (\mbox{\rm co-compact}) \ (\mbox{see {\rm (\ref{eq1.32})}) is QSF} \, \} \, .
\end{equation}
Since Lemma {\rm 2.4} proves that $\Theta^3 ({\rm new}) \in {\rm DE}$, our present lemma proves that $\Theta^3 (\mbox{\rm co-compact}) \in {\rm QSF}$, hence $\forall \, \Gamma \in {\rm QSF}$.}

\bigskip

\noindent {\bf Proof of Lemma 4.1.} The proof in question will occupy the rest of this section. We pick up some finite simplicial complex $k \subset \Theta^3 (\mbox{co-compact})$ and our aim will be to show that there exists a commutative diagram
\begin{equation}
\label{eq4.2}
\xymatrix{
k \ar[rr]^{j_0} \ar[dr]_-{{\rm canonical} \atop {\rm inclusion}} &&K_0 \ar[dl]^-{\chi_0}  \\ 
&\Theta^3 (\mbox{co-compact})
}
\end{equation}
where $K_0$ is an (abstract) compact simply-connected simplicial complex, $j_0$ a simplicial injection, $\chi_0$ a simplicial map and where the Dehn-type condition below is satisfied
$$
M_2 (\chi_0) \cap j_0 \, k = \emptyset, \ \mbox{inside $K_0$}.
\eqno (4.2.1)
$$
This expresses of course that $\Theta^3 (\mbox{co-compact}) \in {\rm QSF}$ and the rule of the game should be that here $k$ is arbitrary. But, clearly, once $k$ has been chosen, enlarging it comes with no harm. So, we can assume to begin with that $k$ is {\ibf connected}. A further extension is presented below. The fact that $\Theta^3 ({\rm new}) \in {\rm DE}$ provides us with an abstract compact simplicial complex $K$ with $\pi_1 K = 0$, coming with a commutative diagram
\begin{equation}
\label{eq4.3}
\xymatrix{
k \subset \Theta^3(\mbox{co-compact}) \ar[dr]_-{i} &\subset &\Theta^3({\rm new}) \supset \underset{n}{\sum} \ {\mathcal B}_n \times [0,\infty) \, ,  \\ 
&K \ar[ur]_{\chi}
}
\end{equation}
when the inclusion $\Theta^3(\mbox{co-compact}) \subset \Theta^3({\rm new}) - \sum {\mathcal B}_n \times (0,\infty)$ is the composition of (\ref{eq1.32}) with (\ref{eq1.31}), where $i$ is a simplicial inclusion, $\chi$ a simplicial {\ibf immersion}, and where the following Dehn-type condition is satisfied
$$
ik \cap M_2 (\chi) = \emptyset \, . \eqno (4.3.1)
$$
Next, keeping $k$ connected and compact, we extend it until the  following conditions are satisfied too

\bigskip

\noindent (4.4.1) \quad If ${\mathcal B}_n$ is any bowl such that $k \cap ({\mathcal B}_n \cap \Theta^3(\mbox{co-compact})) \ne \emptyset$, then ${\mathcal B}_n \cap \Theta^3(\mbox{co-compact}) \subset k$.

\bigskip

\noindent (4.4.2) \quad Let $A_i , B_j$ be like in Lemma 3.2. Then, if $k \, \cap \, A_i \ne 0$, we also have $A_i \subset k$, and also, if $ ({\rm int} \, B_j) \cap k \ne \emptyset$, then $B_j \subset k$.

\bigskip

Since ${\mathcal B}_n \cap \Theta^3(\mbox{co-compact})$, $A_i$, $B_j$ are all three compact, and since $\Theta^3(\mbox{co-compact})$ can only touch finitely many ${\mathcal B}_n$'s [remember that $\Theta^3(\mbox{co-compact}) \subset \Theta^3 ({\rm new})$ and ${\mathcal B}_n \to \infty$ in $\Theta^3 ({\rm new})$, when $n \to \infty$], there is no problem in implementing (4.4.1) and (4.4.2).

\smallskip

Only after $k$ has been extended so that all these conditions are fulfilled, do we fix the $K$, $\chi$ in (\ref{eq4.3}), for the time being, at least. In the arguments which we will develop later, $K$ may change but it will remain compact and simply-connected, $\chi$ may loose its feature of being immersive, but the sacro-sancted Dehn condition
$$
M_2 (\ldots) \cap k = \emptyset
$$
will never be violated.

\smallskip

We certainly have $k \cap \underset{n}{\sum} \, {\mathcal B}_n \times [0,\infty) = \emptyset$, and our FIRST STEP towards Lemma 4.1 will be to demolish $\chi (K) \cap \underset{n}{\sum} \, {\mathcal B}_n \times (0,\infty)$.

\smallskip

We can extend quite naturally the ${\mathcal J}$ (\ref{eq1.30}) to a PROPER
\setcounter{equation}{4}
\begin{equation}
\label{eq4.5}
\sum_n {\mathcal B}_n \times [0,\infty) \overset{{\mathcal J}_1}{-\!\!\!-\!\!\!-\!\!\!-\!\!\!\longrightarrow} \Theta^3 ({\rm new}) \, ,
\end{equation}
such that $M^2 ({\mathcal J}_1) = M^2 ({\mathcal J})$.

\smallskip

According to (4.4.1), there are finitely many ${\mathcal B}_n$'s such that ${\mathcal B}_n \cap \Theta^3 (\mbox{pre-compact}) \subset k$, the others being disjoined from $k$. We have $k \cap {\mathcal B}_n \times [0,\infty) = k \cap {\mathcal B}_n \times \{0\}$ (see (\ref{eq1.31})), but generally speaking $K \cap {\mathcal B}_n \times (0,\infty) \ne \emptyset$, and these are the sets which we want to destroy now, after which we will forget about the $\underset{n}{\sum} \, {\mathcal B}_n \times (0,\infty)$ altogether. We fix a precise ${\mathcal B}_n$ which we may call ${\mathcal B}$. Let
\begin{equation}
\label{eq4.6}
G^3 \mid {\mathcal B}_n \equiv \left\{\mbox{the germ} \ \left( \Theta^3 ({\rm new}) - \sum {\mathcal B} \times (0,\infty) \right) \mid {\mathcal B}_n \times \{0\} \right\} \, .
\end{equation}
In the generic situation of Figure 1.6 we will assume that exactly the $W({\rm BLACK}) + W_1^* + W_2^*$ are there, in such a way that the ${\mathcal B} (W({\rm BLACK}))$ cuts transversally through the immortal singularity $p_{1\infty}$ and not through $p_{i\infty}$, $i > 1$. Also, outside of the plane of Figure 1.6, we have transversal intersection lines
\begin{equation}
\label{eq4.7}
L_1 = {\mathcal B} (W({\rm BLACK})) \cap {\mathcal B} (W_1^*) \, , \qquad L_2 = {\mathcal B} (W({\rm BLACK})) \cap {\mathcal B} (W_{i > 1}^*) \, .
\end{equation}
When we consider $G^3 \mid {\mathcal B} (W({\rm BLACK}))$ (Figure 1.6), then outside of $p_{1\infty}$ and of (\ref{eq4.7}), the $G^3 \mid {\mathcal B} (W({\rm BLACK}))$ is a smooth 3-manifold, just like the ${\mathcal B} \times [0,\infty)$'s.

\smallskip

The $\chi$ in (\ref{eq4.3}) is just a simplicial immersion, but the following things may be assumed, without loss of generality,

\bigskip

\noindent (4.8) \quad Both $K \cap ({\mathcal B} \times [0,\infty)) \equiv \chi^{-1} (\chi \, K \cap ({\mathcal B} \times [0,\infty))$ and $K \cap ({\mathcal B} \times \{0\}) \equiv \chi^{-1} (\chi \, K \cap {\mathcal B} \times \{0\})$ are smooth manifolds of dimensions three and two respectively, on which the restriction of $\chi$ (into ${\mathcal B} \times [0,\infty)$, respectively into ${\mathcal B} \times \{0\}$) is smooth.

\bigskip

\noindent (4.9) \quad When we move from $G^3 \mid {\mathcal B}_n$ to the larger
$$
\overline G^3 \mid {\mathcal B}_n \equiv \{\mbox{the germ of} \ \Theta^3 ({\rm new}) \ {\rm at} \ {\mathcal B}_n \times \{0\}\} \, ,
$$
then, generically, the local structure of $K \cap \overline G^3 \mid {\mathcal B}_n$ is the union along $K \cap ({\mathcal B}_n \times \{0\})$ of $K \cap G^3 \underset{\rm TOP}{=} (K \cap ({\mathcal B}_n \times \{0\})) \times [-\varepsilon , \varepsilon]$ with $K \cap ({\mathcal B}_n \times \{0\}) \times [0,\varepsilon]$, producing a structure $\{$figure $Y$ of vertices $-\varepsilon,+\varepsilon,+\varepsilon\} \times R^2$.

\bigskip

This picture may have to become slightly more complicated when in the presence of the $p_{1\infty} + L_1 + L_2$ mentioned above, but we will ignore this, at least for the time being. With this, we want now to eliminate, successively, all the finitely many ${\mathcal B} \times (0,\infty)$'s which touch $K$. According to (4.8), the $K \cap {\mathcal B}_n \times [0,\infty)$ is a finite union of disjoined components each a smooth 3-manifold generically called $M^3$. The $\partial M^3 \cap ({\mathcal B}_n \times \{0\})$ is a, not necessarily connected, codimension zero submanifold of $\partial M^3$. Call its generic connected component $N$. There is a connected component of $K \cap {\mathcal B}_n \times \{0\}$, call it $N_0$, which is such that
\setcounter{equation}{8}
\begin{equation}
\label{eq4.9}
\{{\mathcal B}_n \cap \Theta^3 (\mbox{co-compact})(\mbox{which by (4.2.1) is contained in $k$}\} \subset N_0 \, , \ \mbox{and this $N_0$ is necessarily {\ibf UNIQUE}}.
\end{equation}
The reason for the uniqueness above is the following. According to the Lemma 3.2, the ${\mathcal B}_n \cap \Theta^3 (\mbox{co-compact})$ is {\ibf connected}. So, imagine now that there are $N'_0 \ne N''_0$ with
$$
N''_0 \supset {\mathcal B}_n \cap \Theta^3 (\mbox{co-compact}) \subset N'_0 \, .
$$
This would contradict then the Dehn property $k \cap M_2 (\chi) = \emptyset$. {\ibf Forgetting temporarily about $k$}, we consider the natural immersion
\begin{equation}
\label{eq4.10}
N_0 \overset{\chi}{-\!\!\!-\!\!\!-\!\!\!-\!\!\!\longrightarrow} {\mathcal B}_n \times \{0\} \, .
\end{equation}
We choose a very dense skeleton $M_0 \subset N_0$ and restrict (\ref{eq4.10}) to it
\begin{equation}
\label{eq4.11}
M_0 \overset{\chi}{-\!\!\!-\!\!\!-\!\!\!-\!\!\!\longrightarrow} {\mathcal B}_n \times \{0\} \, .
\end{equation}

According to our convenience, we may think of the $M_0$ in (\ref{eq4.11}) as being an immersed connected graph, OR as an immersed surface, thin regular neighbourhood of the same graph. From the surface $M_0$ one may get back our initial $N_0$ by adding the small $2$-cells $D_1^2 + D_2^2+ \ldots + D_p^2$. Assuming $M_0$ very dense in $N_0$, these disks are individually embedded by $\chi$ in ${\mathcal B}_n \times \{0\}$.

\smallskip

Continuing to ignore $k$, we replace $K$ by the smaller, still simply-connected object, where small open ``$3$-cells'' get deleted
\begin{equation}
\label{eq4.12}
K_0 \equiv K - \sum_1^p {\rm int} \, D^2_i \times \left( -\frac\varepsilon2 , \frac\varepsilon2 \right) \, .
\end{equation}

We have written here ``$3$-cells'', with quotation marks, since our present $\left( -\frac\varepsilon2 , \frac\varepsilon2 \right)$ is rather a $\{$figure $Y\} - \partial Y$, but this will not change the little argument which will follow next.

\smallskip

For pairs like $(M_0 , K_0)$, immersed into $({\mathcal B}_n \times \{0\} , \Theta^3 ({\rm new}))$, we will consider {\ibf elementary moves}, each consisting of several successive steps

\bigskip

\noindent (4.13.I) \quad Find inside the graph $M_0$ an arc $I = [0,1] \subset M_0$, with $\chi \mid (0,1)$ injective, s.t. $\chi_0 I$ closes to an embedded circle bounding an embedded disk $\delta^2 \subset {\mathcal B}_n \times \{0\}$. This can take one of the three forms displayed in the Figure 4.1, where the $\delta^2$ has been shaded.

\bigskip

We ignore here the other pieces of $M_0$, which are not connected at the source with what is displayed in our Figure 4.1; these pieces may be superposed to it, at the target.

\bigskip

\noindent (4.13.II) \quad We start by adding $\delta^2$ to $\chi (M_0)$ and, at the same time the $3$-cell $\delta^2 \times \left[-\frac\varepsilon4 , \frac\varepsilon4 \right]$, considered here as a {\ibf$2$-handle}, to $\chi \, K_0$; this will require some EXPLANATIONS, which are following now. Let us start with the attachment to $\chi (M_0)$. What this means is the following. To begin with, we consider the abstract object
$$
M_0 \cup \delta^2 \equiv (M_0 + \delta^2) \diagup \{\mbox{the equivalence relation which performs the identifications $\chi (0) = \chi (1)$,}
$$
$$
\mbox{and next, $\partial \delta^2 = \chi (I)$}\} \, .
$$

This object comes endowed with a natural nondegenerate map
$$ M_0 \cup \delta^2 \overset{\chi'}{-\!\!\!-\!\!\!-\!\!\!-\!\!\!\longrightarrow} {\mathcal B}_n \times \{0\} \, ,
$$
which fails to be immersive at some mortal singularities. The $\Psi/\Phi$ abstract nonsense theory, \`a la \cite{22} and \cite{39} can be afterwards applied. When it comes to $K_0$, one has a $\delta^2 \times \left[-\frac\varepsilon4 , \frac\varepsilon4 \right]$ which is actually something like $\delta^2 \times \{{\rm figure} \ Y\}$.

$$
\includegraphics[width=13cm]{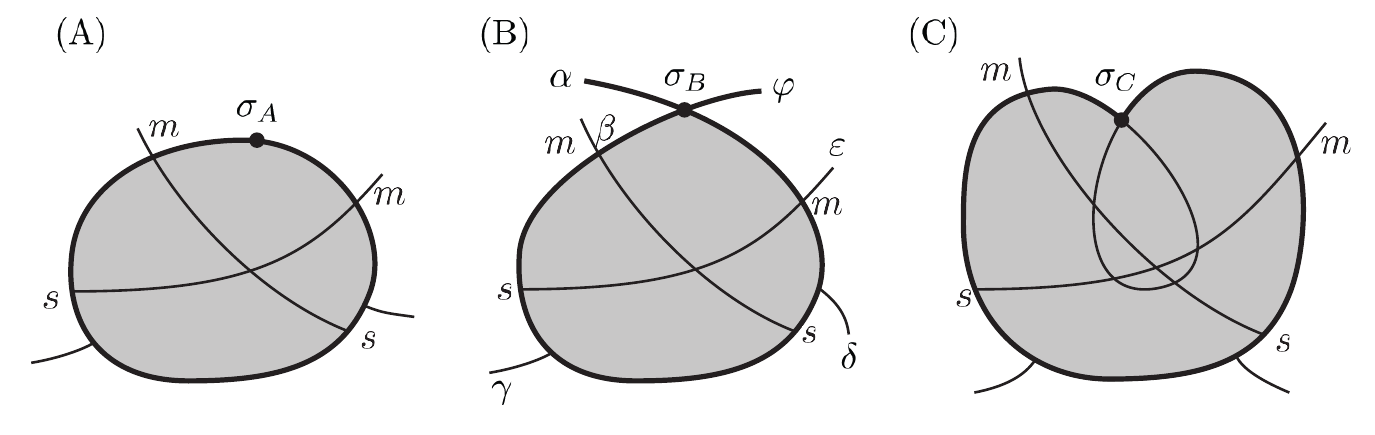}
$$
\label{fig4.1}

\centerline {\bf Figure 4.1.} 

\smallskip

\begin{quote}
We are here inside ${\mathcal B}_n \underset{\rm TOP}{=} R^2$. The circle $\chi (I)$ is drawn in thick lines. At the points marked $m$ or $\sigma$ we see points in $\chi (M_2 (\chi))$, with $\sigma \equiv \{\chi_0 (0) = \chi_0(1)\}$. At the points marked $s$, things are glued to $I \subset M_0$, at the level of the source $M_0$ of $\chi$. The present figure, presents not only the $\delta^2 ({\rm shaded})$, with $S^1 = \chi I = \partial \delta^2$, but also the typical continuation of $\chi M_0$ outside of $\chi (I)$. Sometimes, we write $\chi_0$ for $\chi$.
\end{quote}

\bigskip

Consider now, to begin with the abstract $K_0 \cup \delta^2$ defined by noticing that $M_0 \subset \partial K_0 \times \{0\}$, hence we have $I \to K_0$ and then force the identification $\chi (0) = \chi (1)$ at level $K_0$, after which $\delta^2 = \delta^2 \times \{0\}$ can be attached to get the $K_0 \cup \delta^2$. Without loss of generality, not only do we have $I \subset \partial K_0 \times \{0\}$, but from $I$ start three strata of type $I \times [0,\varepsilon] \subset \partial K_0$. At the level of our $K_0 \cup \delta^2$ we have a singularity $\sigma$ involving three double lines of $K_0 \cup \delta^2 \overset{\chi}{-\!\!\!-\!\!\!\longrightarrow} \Theta^3 ({\rm new})$, along three stata. We zip them along $\left[0,\frac\varepsilon4 \right]$, after which we fill in the missing pieces of the $\delta^2 \times \left[0,\frac\varepsilon4 \right]$ via three dilatations. The result is our $K_0 \cup \delta^2 \times \left[-\frac\varepsilon4 , \frac\varepsilon4 \right]$. End of EXPLANATIONS.

\bigskip

For the newly created objects $M_0 \cup \delta^2$, $K_0 \cup \left( \delta^2 \times \left[-\frac\varepsilon4 , \frac\varepsilon4 \right]\right)$ the points $s$, and $\sigma_0$ in Figure 4.1 are now mortal singularities. When it comes to $K_0 \cup \delta^2 \times \left[-\frac\varepsilon4 , \frac\varepsilon4 \right]$ alone, such singularities also occur at $\sigma_B \times \left\{ \pm \frac\varepsilon4 \right\}$ and $\sigma_C \times \left\{ \pm \frac\varepsilon4 \right\}$. When we talk here about singularities, we have in mind the nondegenerate maps from (4.14.1).

\smallskip

Our step (4.13.II) continues with the commutative diagram below, where all the vertical arrows, except the upper left one, are the obvious inclusions
\setcounter{equation}{13}
\begin{equation}
\label{eq4.14}
\xymatrix{
M'_1 \equiv M_0 \cup \delta^2 \diagup \, \widehat{\rm Cl}_Z^2 (s+\sigma_C) \ar[rr]_{\qquad\qquad\chi'_1} \ar[d]&&{\mathcal B}_n \times \{0\} \ar[d] \\ 
K_1 \equiv K_0 \cup \delta^2 \diagup \, \widehat{\rm Cl}_Z^3 (s+\sigma_B + \sigma_C) \ar[rr]_{\qquad\qquad\chi_1} &&\Theta^3 ({\rm new}) \\ 
M_1 \equiv M'_1 \diagup \, \widehat{\rm Cl}_Z^3 (s+ \sigma_B + \sigma_C) \ar[u] \ar[rr]_{\qquad\qquad\chi_1} &&{\mathcal B}_n \times \{0\}. \ar[u]
}
\end{equation}
Here the $\widehat{\rm Cl}_Z$ are the equivalence relations defined in \cite{29}, the first of the three papers in this series (see, in particular, formula (2.6) in \cite{29}), for the maps
$$
M_0 \cup \delta^2 \longrightarrow {\mathcal B}_n \times \{0\} \quad {\rm and} \quad K_0 \cup \delta^2 \longrightarrow \Theta^3 ({\rm new}) , \quad {\rm respectively}.
\eqno (4.14.1)
$$
The subscript ``$1$'' occurring in the horizontal arrows of (\ref{eq4.14}) is like in formula (2.4) of \cite{39}. The $K^2 \cup \delta^2$ means, of course $K^2 \cup \left( \delta^2 \times \left[-\frac\varepsilon4 , \frac\varepsilon4 \right]\right)$ and, in the middle line of (\ref{eq4.14}) each of the $\sigma_B , \sigma_C$ accounts for two, respectively for four immortal singularities; the $\pm \frac\varepsilon4$ have to be taken into account at $\sigma_C$ too. We have $M_1 \cup \delta^2 \subset K_0 \cup \delta^2$, so the identifications of the second line, affect the third line too. We have now
$$
\pi_1 K_1 = 0 \, , \ \mbox{since we have already $\pi_1 K_0 = \pi_1 K = 0$, and $\pi_1 M'_1 \leq \pi_1 M_0$,}
$$
but the $\widehat{\rm CL}_Z^3$ induces additional identification at the $2^{\rm d}$ level of ${\mathcal B}_n \times \{0\}$, affecting the third line in (\ref{eq4.14}) and the $\pi_1 M_1$ is no longer controlled.

\smallskip

What we have just done is, by definition the {\ibf elementary move}
$$
(K_0 , M_0) \Longrightarrow (K_1,M_1) \, ,
$$
and $(K_1,M_1)$ is just ready for the iteration of the process. We replace now the initial data $\Bigl\{ K_0 \supset M_0 \, \overset{\chi \, \equiv \, \chi_0}{-\!\!\!-\!\!\!-\!\!\!-\!\!\!\longrightarrow}$ ${\mathcal B}_n \times \{0\} \Bigl\}$ by
\begin{equation}
\label{eq4.15}
K_1 \supset M_1 \overset{\chi_1}{-\!\!\!-\!\!\!-\!\!\!-\!\!\!\longrightarrow} {\mathcal B}_n \times \{0\} \, .
\end{equation}
Notice that, in moving from (\ref{eq4.11}) to (\ref{eq4.15}), what we have gained is that
\begin{equation}
\label{eq4.16}
\# \ M^2 (\chi_1 \ (\mbox{from (\ref{eq4.15})}))< \# \ M^2 (\chi_0 \equiv \chi \ (\mbox{from (\ref{eq4.11})})) \, .
\end{equation}

\smallskip

Let me explain this. Look at $(M'_1 , \chi'_1)$ (\ref{eq4.14}) and, for the sake of the argument, we will consider the case (B) of Figure 4.1. Then, when we go to (\ref{eq4.16}) the $M'_1$ gets replaced at the level of our figure, by the {\ibf blub} $\delta^2$ with the six outcoming branches at $\alpha , \beta , \gamma , \delta , \varepsilon , \varphi$, which really is now a 6-valued vertex. Some double points {\ibf do die} in this process, and the phenomenon is quite general.

\bigskip

\noindent {\bf Sublemma 4.2.} 1) {\it After a sufficiently long sequence of elementary moves, and no kind of special strategy is required here
\begin{equation}
\label{eq4.17}
(M_0 , K_0) \Longrightarrow (M_1 , K_1) \Longrightarrow \ldots \Longrightarrow (M_{\omega - 1} , K_{\omega - 1}) \Longrightarrow (M_{\omega} , K_{\omega}) \, ,
\end{equation}
we get a final $M_{\omega} \underset{\chi_{\omega}}{-\!\!\!-\!\!\!-\!\!\!\longrightarrow} {\mathcal B}_n \times \{0\}$, which is such that
$$
\pi_1 M_{\omega} = 0 \quad \mbox{and} \quad \chi_{\omega} \quad \mbox{injects (i.e. $M^2 (\chi_{\omega}) = \emptyset$).}
$$

The $\pi_1 K_i$ continues to stay zero through the whole process.}

\medskip

2) {\it In the end, we get a commutative diagram
\begin{equation}
\label{eq4.18}
\xymatrix{
M_{\omega} \ar[rr]_{\chi_{\omega}} \ar[d] &&{\mathcal B}_n \times \{0\} \ar[d] \\
K_{\omega} \ar[rr]_{\chi_{\omega}} &&\Theta^3 ({\rm new})
}
\end{equation}
where the lower $\chi_{\omega}$ is an immersion and all the other three maps inject.}

\medskip

3) {\it In going from} (\ref{eq4.10}) {\it to} (\ref{eq4.11}), {\it we have an induced map
$$
\sum_{i=1}^P D_i^2 \overset{\chi}{-\!\!\!-\!\!\!-\!\!\!\longrightarrow} {\mathcal B}_n \times \{0\} \, .
$$
We can assume that our} (\ref{eq4.17}) {\it includes enough degenerate elementary moves where $\chi_i \mid [0,1]$ closes already at the source, so that $\chi \mid \underset{i}{\sum} \, D_i^2$ factors through $\chi_{\omega} \, M_{\omega} \subset {\mathcal B}_n \times \{0\}$.

\smallskip

The original $3$-handles of $K-K_0$ have gotten fragmented, each into many mini $3$-handles for $K_{\omega}$, each living either in $\left[ -\frac\varepsilon2 , 0 \right]$ or in $\left[ 0,\frac\varepsilon2 \right]$. Since we continue to find
$$
K_{\omega} \mid \left[ -\frac\varepsilon2 , \frac\varepsilon2 \right] = \chi_{\omega} \, M_{\omega} \times \left[ -\frac\varepsilon2 , \frac\varepsilon2 \right] \, ,
$$
it is possible to add to $K_{\omega}$ the $\Bigl\{${\ibf fragmented} $3$-handles of $\underset{i=1}{\overset{p}{\sum}} \, D_i^2 \times \left[ -\frac\varepsilon2 , \frac\varepsilon2 \right]\Bigl\}$ and they change $K_{\omega}$ into a larger object we call $\overline K_{\omega}$, which stays with $\pi_1 \, \overline K_{\omega} = 0$. All these operations do not touch the existing double points, nor do they create new ones. Diagram {\rm (\ref{eq4.18})} changes now into}
\begin{equation}
\label{eq4.19}
\xymatrix{
\chi_{\omega} \, M_{\omega} = \overline K_{\omega} \cap ({\mathcal B}_n \times 0) \ar[rr] \ar[d] &&{\mathcal B}_n \times \{0\} \, , \ar[d] \\
\overline K_{\omega} \ar[rr]^{\overline\chi_{\omega}} &&\Theta^3 ({\rm new})
}
\end{equation}
{\it when $\overline\chi_{\omega}$ is immersive, all the other arrows are injective, and when $\pi_1 (\chi_{\omega} \, M_{\omega}) = 0$. Also $\pi_1 \, \overline K_{\omega} = 0$, as already said.}

\medskip

4) {\it We finally can {\ibf put back the $k$}, now into the new context. More explicitly, we have a factorization
$$
k \subset \overline K_{\omega} - M_2 (\overline\chi_{\omega}) \subset \Theta^3 ({\rm new})
$$
for the $k \subset \Theta^3 ({\rm new})$ from} (\ref{eq4.3}).

\bigskip

\noindent {\bf Proof.} Via an iterated number of elementary moves one can kill all the double points of $\chi \mid M_0$, being possibly stranded with some $\pi_1 M \ne 0$. This can be killed by some additional disk filling move like in a Figure 4.1.(A), devoid now of any $s$ or $m$. This proves 1) in our lemma and 2) $+$ 3) are left to the reader.

\smallskip

As far as 4) is concerned, for the region $\left[ -\frac\varepsilon2 , \frac\varepsilon2 \right]$ where everything embedds, there is clearly no problem. When we go outside it, one has to notice the following basic fact: we have only applied equivalence relations compatible with $\chi \mid K_0$ and so, because of the Dehn property $k \cap M_2 (\chi) = \emptyset$ of (\ref{eq4.3}) our $k \overset{i}{\hookrightarrow} K$ does not feel the change $K \to \overline K_{\omega}$. \hfill $\Box$

\bigskip

Going back now to (\ref{eq4.9}), what we have gained by our Sublemma 4.2 is that we can also assume now that we also have $\pi_1 N_0 = 0$. Also, in terms of the decomposition into connected components
\begin{equation}
\label{eq4.20}
K \cap ({\mathcal B}_n \times \{0\}) = N_0 \, (4.9) + N_1 + N_2 + \ldots + N_Q \, ,
\end{equation}
so far we have only dealt with $N_0$. But the $N_{i \geq 1}$'s can be treated similarly, things are then even easier, since there is no longer the $k$ to be worried about. So, we may assume, with a possibly new $K$, that in the context of (\ref{eq4.20}) we have $\pi_1 N_i = 0$ for all $i$'s. Here ${\mathcal B}_n$ is generic and we can split away from $K$ all the $K \cap {\mathcal B} \times (0,1)$'s retaining a smaller $K$ which, by Van Kampen, continues to be {\ibf simply-connected}, with $k \subset K$.

\smallskip

The next result is that we have replaced (\ref{eq4.3}) by a new diagram, which still retains $\pi_1 K = 0$, $\chi$ immersive and the Dehn property from (4.3.1), with the following form
\begin{equation}
\label{eq4.21}
\xymatrix{
k \subset \Theta^3(\mbox{co-compact}) \ar[dr]_-{i} &\subset &\Theta^3({\rm new}) - \underset{n}{\sum} \ {\mathcal B}_n \times (0,\infty) \, .  \\ 
&K \ar[ur]_{\chi}
}
\end{equation}

Remember now, that we have gotten so far pretending all the time that we do not encounter the specific difficulties $p_{1\infty} + L_1 + L_2$. I claim that their presence does not change our conclusions, and here is the reason why. To begin with, as we know from (3.5.1), $\Theta^3 (\mbox{co-compact}) \cap (L_1 + L_2) = \emptyset$, and so the presence of the lines $L_1 + L_2$ leaves $k$ untouched.

\smallskip

Next, we certainly have the (4.8), and
$$
({\mathcal B} \times [0,\infty) , {\mathcal B} \times \{0\}) \underset{\rm TOP}{=} (R_+^3 , R^2) \, .
$$
The ${\mathcal B}$'s ride, of course on top of the compensating $2$-handles too. In the process which modifies $(K,\chi)$, when we go from (\ref{eq4.3}) to the (\ref{eq4.21}), each mini step only deals with one individual pair $({\mathcal B} \times [0,\infty) , {\mathcal B} \times \{0\})$ at a time, and $p_{1\infty} + L_1 + L_2$ is not in the way. 

\bigskip

\noindent [EXPLANATIONS. What ``$p_{1\infty}$'' actually stands for, is a contact
$$
{\mathcal B} (W({\rm BLACK})) \cap \{ S (W \cap W^*) , \ \mbox{our $p_{1\infty}$ in Figure 1.6}\} \, , \eqno (*_1)
$$
and $L_1 , L_2$ mean contacts of the type
$$
{\mathcal B} \pitchfork {\mathcal B} \, . \eqno (*_2)
$$
All of $(*_1)$, $(*_2)$ happen far from our compensating $2$-handles $D^2 (p_{\infty\infty} (S))$, on which the ${\mathcal B}$'s may ride.

\smallskip

When dealing with $(*_1)$ we deal essentially with $U^3 ({\rm BLACK})$, Figure 1.6, with the $D^2 (p_{\infty\infty} (S))$ not part of this $U^3$. When we deal with a ${\mathcal B}_n \times \{0\}$ partaking in a context $(*_2)$ we start by replacing ${\mathcal B}_n \times \{0\}$ with ${\mathcal B}_n \times \{\eta > 0\} \subset {\mathcal B}_n \times [0,\infty)$. This leads to a diagram like (\ref{eq4.21}) where our specific ${\mathcal B}_n \times (0,\infty)$ gets replaced by ${\mathcal B}_n \times (\eta,\infty)$ and which is such that
$$
K \cap ({\mathcal B}_n \times [0,\infty)) = K \cap ({\mathcal B}_n \times [0,\eta]) = (K \cap {\mathcal B}_n \times \{0\}) \times [0,\eta] \, ,
$$
and which continues to come with $\pi_1 K = 0$. It is not hard, afterwards, to delete ${\mathcal B}_n \times (0,\eta)$, without creating any harm, and get back exactly the (\ref{eq4.21}).]

\bigskip

So, by now the $\chi (K) \cap \underset{n}{\sum} \, {\mathcal B}_n \times (0,\infty)$ has been demolished, and our FIRST STEP is finished.

\bigskip

\noindent {\bf Second step.} In terms of (\ref{eq4.3}) what we have managed to do, so far, was to take $K$ off the $\underset{n}{\sum} \, {\mathcal B}_n \times (0,\infty)$ and now we want to take it off the $\Theta^3_0$ (see (\ref{eq3.4})) too.

\smallskip

The first ministep will be to demolish the intersections
$$
K \cap \sum_i (\pi^{-1} A_i - A_i) \ \mbox{with $A_i$ like in (\ref{eq3.6}) and $\pi$ in (\ref{eq3.8})}.
$$
The (4.4.2) is with us, and we consider first the case when $A_1 \subset k$; for the corresponding collapsible space $\pi^{-1} A_1$ we have, of course that, $\pi^{-1} A_1 \subset \Theta^3_0$. There is a decomposition into finitely many connected components
$$
K \cap \pi^{-1} A_1 \equiv \chi^{-1} (\chi \, K \cap \pi^{-1} A_1) = C_1 + C_2 + \ldots + C_{\lambda} \, . \eqno (4.21.1)
$$
The $C_i$'s are connected, codimension one subcomplexes of $K$. We have assumed that $A_1 \subset k$; then, up to a notational change one may assume that $A_1 \subset C_1$ and $C_i \cap A_i = \emptyset$ if $i > 1$. Inside $\Theta^3 ({\rm new}) - \underset{n}{\sum} \, {\mathcal B}_n \times (0,\infty)$ we have here $A_1 \subset \partial \, \pi^{-1} A_1$, while at least at the level of $K$, each of the
$$
(\pi^{-1} A_1 - A_1) \mid C_{i > 1}
$$
induces a clear splitting. For the $A_1 \subset \partial \, \pi^{-1} A_1$, the $A_1 \subset K$ is of codimension two. Here it is only along
$$
(\pi^{-1} A_1 - A_1) \mid (C_1 - A_1)
$$
that there is a clean splitting, while along $A_1$ things stick.

\smallskip

Generally speaking, also, $\pi_1 C_i \ne 0$.

\smallskip

Each of the immersion
\begin{equation}
\label{eq4.22}
C_i \overset{\chi \, \mid \, C_i}{-\!\!\!-\!\!\!-\!\!\!-\!\!\!-\!\!\!-\!\!\!-\!\!\!-\!\!\!\longrightarrow} \pi^{-1} A_1 \, , \quad i > 1
\end{equation}
will be replaced by the following {\ibf simplicial map}, no longer immersive (in general), but which has a simply-connected source
\begin{equation}
\label{eq4.23}
{\rm Cone} \, (C_i) \overset{\chi_i \, \equiv \, {\rm Cone} \, (\chi \, \mid \, C_i)}{-\!\!\!-\!\!\!-\!\!\!-\!\!\!-\!\!\!-\!\!\!-\!\!\!-\!\!\!-\!\!\!-\!\!\!-\!\!\!-\!\!\!-\!\!\!-\!\!\!\longrightarrow} \pi^{-1} A_1 - A_1 \, .
\end{equation}

All this was for $C_{i > 1}$. For $C_1$, we can define just like above
$$
{\rm Cone} \, (C_1) \overset{\chi_1}{-\!\!\!-\!\!\!-\!\!\!\longrightarrow} \pi^{-1} A_1 \, , \ \mbox{coming now with} \ \chi_1 ({\rm Cone} \, (C_1)) \cap A_1 \ne \emptyset \, . \eqno (4.23.1)
$$

But, once $A_1 \subset \partial \, \pi^{-1} A_1$ we can certainly ask, in the context of (4.23.1) that we should also have
\begin{equation}
\label{eq4.24}
\chi_1 \, M^2 (\chi_1) \cap A_1 = \emptyset \, .
\end{equation}

One should notice that at the level of these last moves we started moving from Dehn-exhaustibility to the weaker QSF property, where the map $\chi$ (\ref{eq4.21}) looses its immersive property, remaining a mere simplicial (i.e. continuous) map.

\smallskip

Inside $\Theta_0^3$ we define now the following simplicial complexes
\begin{equation}
\label{eq4.25}
K_1 \equiv \{K \ \mbox{split along} \ \pi^{-1} A_1 - A_1\} \, .
\end{equation}

For each $C_{i>1}$ we find now two copies  of $C_i$, each of them with $C_i^{\pm} \subset K_1$. We also find two copies of $C_i^{\pm} - A_1 \subset K_1$.

\bigskip

\noindent (4.26) \quad From $K_1$ we go to $K_2$ by adding, to begin with, for each $C_{i>1}^{\pm}$ a copy of $\chi$ (cone $C_i^{\pm}$), defined like in (\ref{eq4.23}), on the corresponding side of $K_1$, with respect to the split. Then, proceeding like in (4.23.1), we also add two copies of $\chi_1 ({\rm cone} \, C_1^{\pm}) \mid (C_1^{\pm} - A_1)$.

\bigskip

With the unique, obvious, $A_1 \subset K_1$, then extend to two complete copies of $\chi_1 ({\rm cone} \, C_1^{\pm}) \subset K_2$.

\bigskip

\noindent {\bf Claim (4.27).} 1) $\pi_1 K_2 = 0$.

\medskip

2) {\it There is a {\ibf simplicial map} (no longer an immersion!)}
\setcounter{equation}{27}
\begin{equation}
\label{eq4.28}
K_2 \overset{\chi(2)}{-\!\!\!-\!\!\!-\!\!\!\longrightarrow} \Theta^3 ({\rm new}) - \sum_n {\mathcal B}_n \times (0,\infty) \, .
\end{equation}

\medskip

3) {\it We have}
\begin{equation}
\label{eq4.29}
\chi (2) \, M^2 (\chi (2)) \cap A_1 = \emptyset \, .
\end{equation}

\medskip

3.bis) {\it The original inclusion $k \subset K$ from {\rm (\ref{eq4.21})} induces an inclusion $k \subset K_2$, and since $k \cap \pi^{-1} A_1 = A_1$, we also have}
\begin{equation}
\label{eq4.30}
k \cap M_2 (\chi (2)) = \emptyset \, .
\end{equation}

\medskip

4) {\it By a small perturbation of {\rm (\ref{eq4.28})}, localized inside ${\rm int} \, \Theta^3_0$, we can disentangle completely $K_2$ from $\pi^{-1} A_1$, thereby replacing {\rm (4.21.1)} by the following formula, where the LHS should be read in the manner of {\rm (4.21.1)}}
\begin{equation}
\label{eq4.31}
K_2 \cap \pi^{-1} A_1 = A_1 \subset k \, .
\end{equation}

\bigskip

By iterating the process
$$
(K,\chi) \Longrightarrow (K_2 , \chi (2))
$$
we can get a $(K_m , \chi (m))$ which is now disentangled completely from $\underset{j}{\sum} \, \pi^{-1} A_j$, and which we continue to call $(K,\chi)$.

\smallskip

We still have to deal with
$$
K \cap \Theta^3_0 = \left\{ K \cap \sum_j \pi^{-1} B_j \, , \ \mbox{with all the contribution of} \ K \cap \sum_i (\pi^{-1} A_i - A_i) \ \mbox{by now already removed}\right\} .
$$
Each $B_j \subset S_0$ is collapsible and either $B_j \subset k$ or $({\rm int} \, B_j) \cap k = \emptyset$. Also, as just explained, all the contribution $(\pi^{-1} A_i - A_i) \cap K \subset \pi^{-1} B_j$ has already been dealt with.

\smallskip

So, let us move to the most complicated case when $B_j \subset k \subset K$. The codimension one space
$$
\widehat B_j \equiv B_j \cup \sum_{\overbrace{\mbox{\footnotesize$A_i \cap \partial B_j$}}} \pi^{-1} A_i \subset \Theta_0^3 \subset \Theta^3({\rm new}) - \sum_n {\mathcal B}_n \times (0,\infty)
$$
splits. Also, $\pi_1 \, \widehat B_j = 0$, and $\widehat B_j$ does not touch the double points of the map $K \longrightarrow \Theta^3 ({\rm new}) - \underset{n}{\sum} \, {\mathcal B}_n \times (0,\infty)$.

\smallskip

We are now in a context similar with the one of our previous dealings with ${\mathcal B} \times (0,\infty)$ or with $\pi^{-1} A_i$, but easier.

\smallskip

After an appropriate cone-construction, in the style of (\ref{eq4.23}), or (4.23.1), the $\widehat B_j$ splits $K$ into two simply-connected pieces only one of which fully contains $k$. So, we can happily replace $K$ by $K - (\pi^{-1} B_j - \widehat B_j)$. By a finite iteration we realize
$$
K \cap \Theta_0^3 = \left( \sum_j \pi^{-1} B_j \right) \cap K = S_0 \cap k \, .
$$

This, finally, replaces (\ref{eq4.21}) with a diagram of the desired form (\ref{eq4.2}).

\end{document}